\documentclass[12pt]{article}

\usepackage{setspace}
\usepackage{lscape}
\usepackage{pdflscape}
\usepackage{ulem}
\usepackage{amsmath,amssymb,amsthm} 

\usepackage{wrapfig} 

\usepackage{mathrsfs} 
\usepackage{url} 
\usepackage{latexsym} 
\usepackage{bm} 

\usepackage{mathtools}

\usepackage{booktabs,multirow}
\usepackage{array}
\newcolumntype{P}[1]{>{\centering\arraybackslash}p{#1}}
\usepackage{diagbox}
\usepackage{threeparttable}
\usepackage{siunitx} 

\usepackage{algorithm}
\usepackage[noend]{algpseudocode}

\usepackage[sort&compress,numbers]{natbib}

\usepackage{hyperref}
\usepackage{cleveref}
\usepackage{autonum} 

\usepackage{thmtools,thm-restate}

\usepackage{tikz,tikz-3dplot}
\usetikzlibrary{positioning,shapes,shapes.callouts,shapes.multipart,arrows.meta,fit,calc}
\usepackage{pgfplots}
\pgfplotsset{compat=newest}
\usepackage{pxpgfmark} 
\usetikzlibrary{bayesnet} 

\usepackage{multicol}
\usepackage{wasysym} 
\usepackage{adjustbox}
\usepackage{xparse} 
\usepackage{pbox} 
\usepackage{authblk} 

\usepackage[subrefformat=parens]{subcaption}
\usepackage{enumitem} 
\hypersetup{
  colorlinks,
  linkcolor={green!35!black},
  citecolor={green!35!black},
  urlcolor={blue!50!black}
}

\DeclareMathOperator{\diag}{diag}

\DeclareMathOperator{\rank}{rank}

\DeclareMathOperator{\proj}{proj}

\DeclareMathOperator{\dist}{dist}

\DeclareMathOperator*{\argmin}{argmin}

\DeclareMathOperator{\EE}{\mathbb{E}}

\DeclarePairedDelimiter\ceil{\lceil}{\rceil}

\DeclarePairedDelimiterX\inner[2]{\langle}{\rangle}{{#1},{#2}}

\DeclarePairedDelimiter\norm{\|}{\|}
\DeclarePairedDelimiter\set{\{}{\}}
\DeclarePairedDelimiter\prn{(}{)}
\DeclarePairedDelimiter\bra{[}{]}
\DeclarePairedDelimiterX\Set[2]{\{}{\}}{\mspace{2mu}{#1}\;\delimsize|\;{#2}\mspace{2mu}}
\DeclarePairedDelimiterX\Prn[2]{(}{)}{\mspace{2mu}{#1}\;\delimsize|\;{#2}\mspace{2mu}}
\DeclarePairedDelimiterX\Bra[2]{[}{]}{\mspace{2mu}{#1}\;\delimsize|\;{#2}\mspace{2mu}}

\newcommand{\N}{\mathbb N}
\newcommand{\Z}{\mathbb Z}
\newcommand{\R}{\mathbb R}

\newcommand{\0}{\mathbf 0}

\newcommand{\subjectto}{\mathrm{subject\ to}}

\renewcommand{\epsilon}{\varepsilon}

\NewDocumentCommand{\exsub}{s m O{} m}{%
  \IfBooleanT{#1}{\EE_{#2}\nolimits\bra*{#4}}%
  \IfBooleanF{#1}{\EE_{#2}\nolimits\bra[#3]{#4}}%
}

\algnewcommand{\Break}{\textbf{break}}
\algnewcommand{\Continue}{\textbf{continue}}
\algnewcommand{\algorithmicgoto}{\textbf{go to} Line}
\algnewcommand{\Goto}[1]{\algorithmicgoto~\ref{#1}}
\algblockdefx{MRepeat}{EndRepeat}{\textbf{repeat}}{}
\algnotext{EndRepeat}

\usepackage{CJKutf8}

\newcommand{\mathInd}{\hphantom{{}={}}}
\newcommand{\by}[2][]{\text{\pbox[c]{\textwidth}{(by \pbox[t]{\textwidth}{\,\!#2)#1}}}}

\newcommand{\tred}{\textcolor{red!60!black}}
\newcommand{\cstEB}{\rho}
\newcommand{\cstSIL}{c}
\newcommand{\myparagraph}[1]{\medskip\noindent\textbf{#1}\ }
\renewcommand{\emph}[1]{\textit{#1}}

\declaretheoremstyle[
shaded={bgcolor=gray!15},
]{thmsty}

\declaretheorem[
name=Theorem,
style=thmsty,
]{theorem}

\declaretheorem[
name=Corollary,
style=thmsty,
]{corollary}

\declaretheorem[
name=Proposition,
style=thmsty,
]{proposition}

\declaretheorem[
name=Lemma,
style=thmsty,
]{lemma}

\declaretheorem[
name=Definition,
style=thmsty,
]{definition}

\declaretheorem[
name=Assumption,
style=thmsty,
]{assumption}

\declaretheorem[
name=Condition,
style=thmsty,
]{condition}

\declaretheorem[
name=Remark,
style=thmsty,
]{remark}

\crefname{theorem}{Theorem}{Theorems}
\crefname{corollary}{Corollary}{Corollaries}
\crefname{proposition}{Proposition}{Propositions}
\crefname{lemma}{Lemma}{Lemmas}
\crefname{definition}{Definition}{Definitions}
\crefname{assumption}{Assumption}{Assumptions}
\crefname{remark}{Remark}{Remarks}
\crefname{algorithm}{Algorithm}{Algorithms}
\crefname{line}{Line}{Lines}
\crefname{section}{Section}{Sections}
\crefname{appendix}{Section}{Sections}
\crefname{condition}{Condition}{Conditions}
\crefname{proof}{Proof}{Proofs}
\crefname{table}{Table}{Tables}
\crefname{figure}{Figure}{Figures}
\crefname{equation}{}{}
\Crefname{equation}{Eq.}{Eqs.}

\newlist{enuminasm}{enumerate}{1} 
\setlist[enuminasm]{
  label=(\roman*), ref=\theassumption(\roman*),
  leftmargin=0.06\textwidth,
  topsep=0.2\baselineskip,
  itemsep=0.2\baselineskip,
  font=\normalfont}
\crefalias{enuminasmi}{assumption}

\newlist{enumincond}{enumerate}{1} 
\setlist[enumincond]{
  label=(\roman*), ref=\theassumption(\roman*),
  leftmargin=0.06\textwidth,
  topsep=0.2\baselineskip,
  itemsep=0.2\baselineskip,
  font=\normalfont}
\crefalias{enumincondi}{condition}

\newlist{enuminthm}{enumerate}{1}
\setlist[enuminthm]{
  label=(\roman*), ref=\thetheorem(\roman*),
  leftmargin=0.06\textwidth,
  topsep=0.2\baselineskip,
  itemsep=0.2\baselineskip,
  font=\normalfont}
\crefalias{enuminthmi}{theorem}

\newlist{enuminlem}{enumerate}{1}
\setlist[enuminlem]{
  label=(\roman*), ref=\thelemma(\roman*),
  leftmargin=0.06\textwidth,
  topsep=0.2\baselineskip,
  itemsep=0.2\baselineskip,
  font=\normalfont}
\crefalias{enuminlemi}{lemma}

\usepackage[margin=1.0truein]{geometry}
\providecommand{\keywords}[1]{{\small \textbf{Keywords:} \pbox[t]{0.8\linewidth}{#1}}}

\title{
  Majorization-minimization-based Levenberg--Marquardt method for constrained nonlinear least squares%
}

\author[1]{Naoki Marumo\footnote{Corresponding author. E-mail: \url{naoki_marumo@mist.i.u-tokyo.ac.jp}}}
\author[2,3]{Takayuki Okuno}
\author[1,3]{Akiko Takeda}

\affil[1]{Graduate School of Information Science and Technology, University of Tokyo, Tokyo, Japan}
\affil[2]{Faculty of Science and Technology, Seikei University, Tokyo, Japan}
\affil[3]{Center for Advanced Intelligence Project, RIKEN, Tokyo, Japan}

\begin{document}
\maketitle

\begin{abstract}
  A new Levenberg--Marquardt (LM) method for solving nonlinear least squares problems with convex constraints is described. 
  Various versions of the LM method have been proposed, their main differences being in the choice of a damping parameter.
  In this paper, we propose a new rule for updating the parameter so as to achieve both global and local convergence even under the presence of a convex constraint set. 
  The key to our results is a new perspective of the LM method from majorization-minimization methods. 
  Specifically, we show that if the damping parameter is set in a specific way, the objective function of the standard subproblem in LM methods becomes an upper bound on the original objective function under certain standard assumptions.

  Our method solves a sequence of the subproblems approximately using an (accelerated) projected gradient method.
  It finds an $\epsilon$-stationary point after $O(\epsilon^{-2})$ computation and achieves local quadratic convergence for zero-residual problems under a local error bound condition.
  Numerical results on compressed sensing and matrix factorization show that our method converges faster in many cases than existing methods.
\end{abstract}


\section{Introduction}
\label{sec: introduction}
In this study, we consider the constrained nonlinear least-squares problem:
\begin{equation}
  \min_{x \in \R^d} \
  f(x)
  \coloneqq
  \frac{1}{2}\|F(x)\|^2
  \quad
  \subjectto\quad x \in \mathcal C,
  \label{eq:mainproblem}
\end{equation}
where $\norm{\cdot}$ denotes the $\ell_2$-norm, $F: \R^d \to \R^n$ is a continuously differentiable function, and $\mathcal C \subseteq \R^d$ is a closed convex set. 
If there exists a point $x \in \mathcal C$ such that $F(x) = \0$, 
the problem is said to be \emph{zero-residual}, 
and is reduced to the constrained nonlinear equation:
\begin{equation}
  \text{find} \quad
  x \in \mathcal C\quad
  \text{such that}\quad
  F(x) = \0.
  \label{eq:mainproblem equation}
\end{equation}
Such problems cover a wide range of applications, 
including chemical equilibrium systems~\citep{meintjes1987methodology}, 
economic equilibrium problems~\citep{dirkse1995mcplib}, 
power flow equations~\citep{wood2013power}, 
nonnegative matrix factorization~\citep{berry2007algorithms,lee2000algorithms}, 
phase retrieval~\citep{candes2015phase,zhang2017nonconvex}, 
nonlinear compressed sensing~\citep{blumensath2013compressed},
and learning constrained neural networks~\citep{chorowski2014learning}.

\emph{Levenberg--Marquardt (LM) methods}~\citep{levenberg1944method, marquardt1963algorithm} are efficient iterative algorithms for solving problem \eqref{eq:mainproblem}; 
they were originally developed for unconstrained cases (i.e., $\mathcal C = \R^d$) 
and later extended to constrained cases by \citep{kanzow2004levenberg}.
Given a current point $x_k \in \mathcal C$, 
an LM method defines a model function $m^k_\lambda: \R^d \to \R$ 
with a damping parameter $\lambda > 0$:
\begin{equation}
  m^k_\lambda(x)
  \coloneqq
  \frac{1}{2}\|F_k + J_k(x-x_k)\|^2 + \frac{\lambda}{2}\|x-x_k\|^2,
  \label{eq: def of model}
\end{equation}
where $F_k \coloneqq F(x_k) \in \R^d$ {and} $J_k \coloneqq J(x_k) \in \R^{n \times d}$ {with} 
$J: \R^d \to \R^{n \times d}$ {being} the Jacobian matrix function of $F$. 
The next point $x_{k+1} \in \mathcal C$ is set to an exact or approximate solution to the convex subproblem:
\begin{equation}
  \min_{x \in \R^d}\
  m^k_\lambda(x)\quad
  \subjectto \quad
  x \in \mathcal C
  \label{eq:subproblem}
\end{equation}
for some $\lambda = \lambda_k$. 
Various versions of this method have been proposed, 
and their theoretical and practical performances largely depend on 
how the damping parameter $\lambda_k$ is updated.

\subsection{Our contribution}
We propose an LM method with a new rule for updating $\lambda_k$.
Our method is based on \emph{majorization-minimization} (MM) methods, which successively minimize a \emph{majorization} or, in other words, an upper bound on the objective function. 
The key to our method is the fact that the model $m^k_\lambda$ defined in \cref{eq: def of model} is a majorization of the objective $f$ under certain standard assumptions. 
This MM perspective enables us to create an LM method with desirable properties, including global and local convergence guarantees.
Although there exist several MM methods for problem~\cref{eq:mainproblem} and relevant problems~\citep{nesterov2007modified,nesterov2006cubic,griewank1981modification,bellavia2015strong,bellavia2010convergence}, as far as we know, no studies have elucidated that the model in \cref{eq: def of model} is a majorization of $f$. 
Another feature of our LM method is the way of generating an approximate solution of subproblem~\eqref{eq:subproblem}.
It is sufficient to apply one iteration of a projected gradient method to \eqref{eq:subproblem} for deriving the iteration complexity of our LM method, which leads to an overall complexity bound.

Our contributions are summarized as follows:
\begin{enumerate}[
    label=(\roman*), ref=(\roman*),
    leftmargin=0.06\textwidth,
    topsep=0.3\baselineskip,
    itemsep=0.1\baselineskip
  ]
  \item
  \label{item:mmperspective}
  \textbf{A new MM-based LM method:}
  We prove that the LM model defined in \cref{eq: def of model} is a majorization of $f$ if the damping parameter $\lambda$ is sufficiently large.
  See \cref{lem: majorization LM} for a precise statement. This result provides us with a new update rule of $\lambda$, bringing about a new LM method for solving problem~\eqref{eq:mainproblem}.

  \item
  \label{item:property_global}
  \textbf{Iteration and overall complexity for finding a stationary point:}
  The iteration complexity of our LM method for finding an $\epsilon$-stationary point (see \cref{def:eps-stationary-point}) is proved to be $O(\epsilon^{-2})$ under mild assumptions on the Jacobian.
  Because the computational complexity per iteration of our method does not depend on $\epsilon$, the overall complexity is also evaluated as $O(\varepsilon^{-2})$ through
  \[
    (\text{Overall complexity}) = 
    (\text{Iteration complexity}) \times (\text{Complexity per iteration}).
  \]
  See \cref{cor:iteration_complexity,cor: global overall complexity} for a precise statement.

  \item
  \label{item:property_local}
  \textbf{Local quadratic convergence:}
  For zero-residual problems, assume that a starting point $x_0 \in \mathcal C$ is sufficiently close to an optimal solution, and assume standard conditions, including a local error bound condition.
  Then, if the subproblems are solved with sufficient accuracy, a solution sequence $(x_k)$ generated by our method converges quadratically to an optimal solution.
  See \cref{thm: local superliner convergence and unsuccessful} for a precise statement.

  \item
  \label{item:ext_const}
  \textbf{Improved convergence results even for unconstrained problems:}
  Our method achieves both the $O(\epsilon^{-2})$ iteration complexity bound and local quadratic convergence.
  An LM method having such global and local convergence results is new for unconstrained and constrained problems, as shown in \cref{table: comparison of LM methods}.
\end{enumerate}


Numerical results show that our method converges faster and is more robust than existing LM-type methods~\citep{fan2013levenberg,kanzow2004levenberg,facchinei2013family,goncalves2021inexact}, a projected gradient method, and a trust-region reflective method~\citep{branch1999subspace,virtanen2020scipy}.

\subsection{Oracle model for overall complexity bounds}
\label{sec:oracle_model}
To evaluate the overall complexity of LM methods, we count the number of \emph{basic operations}---evaluation of $F(x)$, Jacobian-vector multiplications $J(x) u$ and $J(x)^\top v$, and projection onto $\mathcal C$---required to find an $\epsilon$-stationary point, following \citep[Section 6]{drusvyatskiy2019efficiency}.
The important point is that we do not assume an evaluation of $J_k \coloneqq J(x_k)$ but access the Jacobian only through products $J_k u$ and $J_k^\top v$ to solve subproblem~\cref{eq:subproblem}.
Computing vectors $J_k u$ and $J_k^\top v$ for given $u \in \R^d$ and $v \in \R^n$ is much cheaper than evaluating the matrix $J_k$.\footnote{
  Automatic differentiation libraries such as JAX~\citep{jax2018github} compute the Jacobian-vector products at several times the cost of evaluating $F(x)$.
  See, e.g., the JAX documentation \citep{jax2022autodiff}.
}
Avoiding the computation of the $n \times d$ matrix $J_k$ makes algorithms practical for large-scale problems where $n$ and $d$ amount to thousands or millions.
We note that some existing LM-type methods \citep{bellavia2010convergence,bellavia2015strong,cartis2011evaluation,cartis2011adaptive,cartis2011adaptive2,cartis2012adaptive,cartis2020strong,goncalves2021inexact} compute the Jacobian explicitly.

\subsection{Paper organization}
In \cref{sec:related_work}, we review LM methods and related algorithms for problem~\cref{eq:mainproblem}.
In \cref{sec:proposed_method}, a key lemma is presented and the LM method (\cref{alg: proposed LM-GD}) is derived based on the lemma.
\cref{sec: global convergence analysis,sec: local convergence analysis} show theoretical results for \cref{alg: proposed LM-GD}: iteration complexity, overall complexity, and local quadratic convergence.
In \cref{sec:generalization}, we generalize \cref{alg: proposed LM-GD} and present a more practical variant of \cref{alg: proposed LM-GD}.
This variant also achieves the theoretical guarantees given for \cref{alg: proposed LM-GD} in \cref{sec: global convergence analysis,sec: local convergence analysis}.
\cref{sec:experiments} provides some numerical results and \cref{sec:conclusion} concludes the paper.

\subsection{Notation}
Let $\R^d$ denote a $d$-dimensional Euclidean space equipped with the $\ell_2$-norm $\norm{\cdot}$ and the standard inner product $\inner{\cdot}{\cdot}$.
For a matrix $A \in \R^{m \times n}$, let $\norm{A}$ denote its spectral norm, or its largest singular value. 
For $a \in \R$, let $\ceil{a}$ denote the least integer greater than or equal to $a$.


\section{Comparison with related works}
\label{sec:related_work}
We review existing methods for problem~\cref{eq:mainproblem} and compare them with our work.

\begin{figure}
  \centering
  \includegraphics[width=0.5\linewidth]{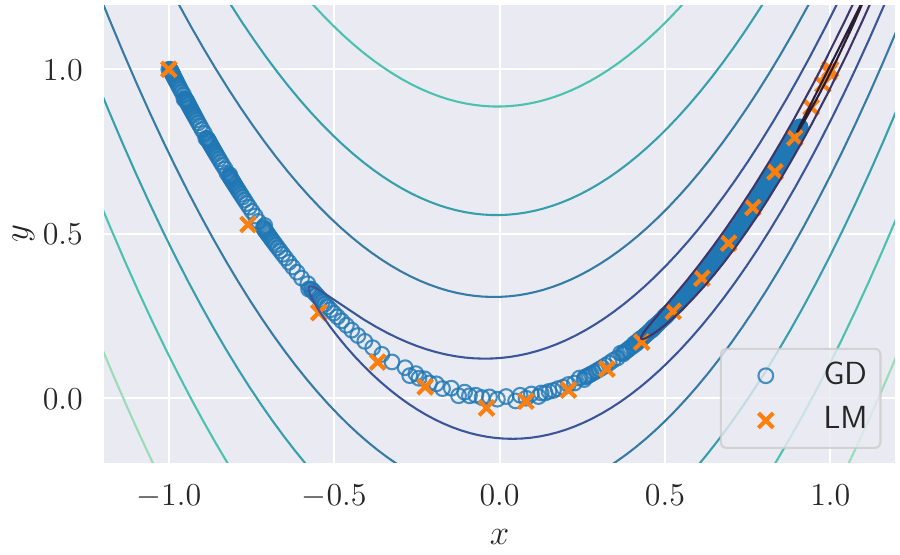}
  \caption{
    Minimization of the Rosenbrock function~\citep{rosenbrock1960automatic}, $f(x, y) = (x - 1)^2 + 100 (y - x^2)^2$.
    Both the gradient descent (GD) and our LM start from $(-1, 1)$ and converge to the optimal solution, $(1, 1)$.
    One marker corresponds to one iteration, and the GD and LM are truncated after 1000 and 20 iterations, respectively.
  }
  \label{fig:rosenbrock}
\end{figure}

\subsection{General methods}
Algorithms for general nonconvex optimization problems, not just for least-squares problems, also solve problem~\cref{eq:mainproblem}.
For example, the projected gradient method have an overall complexity bound of $O(\epsilon^{-2})$; our LM method enjoys local quadratic convergence in addition to that bound, which seems difficult to achieve with general first-order methods.
\cref{fig:rosenbrock} illustrates that our LM successfully minimizes the Rosenbrock function, a valley-like function that is notoriously difficult to minimize numerically.
Although quadratic convergence is proved only locally around an optimal solution, in practice, the LM method may perform considerably better than general first-order methods, even when started far from the optimum.

Some methods, such as the Newton method, achieve local quadratic convergence using second-order or higher-order derivatives of $f$; our LM achieves it without the second-order derivative.
Besides the fact that our LM does not require a computationally demanding Hessian matrix, it has another advantage: subproblem~\cref{eq:subproblem} is very tractable.
Whereas our subproblem is smooth and strongly convex, those in second- or higher-order methods are nonconvex in general.
The matter becomes more severe under the presence of constraints because the subproblems may be NP-hard, as pointed out in~\citep{cartis2012adaptive}.

\subsection{Specialized methods for least squares}
Several methods, including the LM method, utilize the least-squares structure of problem~\cref{eq:mainproblem}.
Focusing on those algorithms without second-order derivatives, we review them from three points of view: (i) subproblem, (ii) complexity for finding a stationary point, and (iii) local superlinear convergence.
Most of the methods discussed in this section are summarized in \cref{table: comparison of LM methods}.
The table shows the following:
\begin{itemize}
  \item
  Our method can achieve an overall computational complexity bound, $O(\epsilon^{-2}) \times O(1) = O(\epsilon^{-2})$, for finding an $\epsilon$-stationary point for constrained problems.
  \item
  To the best of our knowledge, this is the first LM that achieves such a complexity bound with local quadratic convergence, even for unconstrained problems.
\end{itemize}

\begin{table}
  \centering
  \begin{threeparttable}
    \caption{
      Comparison of methods for problem~\cref{eq:mainproblem}.
    }
    \label{table: comparison of LM methods}
    \def\arraystretch{0.95}
    \footnotesize
    \begin{tabular}{@{}ccccccc@{}}\toprule
      \multirow{2}{*}{Subproblem}                                                       & \multirow{2}{*}{References}                                                                       & \multirow{2}{*}{Constr.}    & \multicolumn{2}{c}{Complexity}                                                & \multicolumn{2}{c}{Local conv.} \\[-0.4\baselineskip]\cmidrule(lr){4-5}\cmidrule(lr){6-7}
                                                                                        &                                                                                                   &                             & \#iterations                           & complexity/iter.                          & order            & inexact           \\[-0.2\baselineskip]\midrule
      \cref{eq:subproblem} (LM)                                                         & \citep{osborne1976nonlinear,ueda2010global,zhao2016global}                                        &                             & $O(\epsilon^{-2})$                     &                                      &                  &                   \\\cmidrule{2-7}
                                                                                        & \citep{bellavia2018levenberg}                                                                     &                             & $O(\epsilon^{-2})$                     & $O(1)$                               &                  &                   \\\cmidrule{2-7}
                                                                                        & \citep{bergou2020convergence}\footnotemark[1]                                                     &                             & $O(\epsilon^{-2} \log \epsilon^{-1})$  & $O(1)$                               & $2$              &                   \\\cmidrule{2-7}
                                                                                        & \citep{yamashita2001rate,fan2001convergence,fan2005quadratic,fan2003modified,fan2006convergence}  &                             &                                        &                                      & $2$              &                   \\\cmidrule{2-7}
                                                                                        & \citep{dan2002convergence,fan2004inexact,fischer2010inexactness,fan2011convergence}               &                             &                                        &                                      & $2$              & \checkmark        \\\cmidrule{2-7}
                                                                                        & \citep{kanzow2004levenberg,fan2013levenberg}                                                      & \checkmark                  &                                        &                                      & $2$              &                   \\\cmidrule{2-7}
                                                                                        & \citep{behling2012unified,facchinei2013family}                                                    & \checkmark                  &                                        &                                      & $2$              & \checkmark        \\\cmidrule{2-7}
                                                                                        & \textbf{This work}                                                                                & \checkmark                  & $O(\epsilon^{-2})$                     & $O(1)$                               & $2$              & \checkmark        \\\cmidrule{1-7}
      \cref{eq:subproblem_1_2}                                                          & \citep{nesterov2007modified}\footnotemark[2]                                                      &                             & $O(\epsilon^{-2})$                     &                                      & $2$              &                   \\\cmidrule{2-7}
                                                                                        & \citep{cartis2011evaluation}                                                                      &                             & $O(\epsilon^{-2})$                     &                                      &                  &                   \\\cmidrule{2-7}
                                                                                        & \citep{cartis2020strong}                                                                          & \checkmark                  & $O(\epsilon^{-2})$                     &                                      &                  &                   \\\cmidrule{2-7}
                                                                                        & \citep{bellavia2015strong,bellavia2010convergence}                                                &                             &                                        &                                      & $2$              & \checkmark        \\\cmidrule{1-7}
      \multirow{2}{*}{\pbox[b]{12em}{\cref{eq:subproblem_2_3} and its\\generalization}} & \citep{cartis2011adaptive,cartis2011adaptive2}                                                    &                             & $O(\epsilon^{-2})$                     & $O(1)$                               &                  &                   \\\cmidrule{2-7}
                                                                                        & \citep{cartis2012adaptive}                                                                        & \checkmark                  & $O(\epsilon^{-2})$                     & $O(1)$                               &                  &                   \\\cmidrule{2-7}
                                                                                        & \citep{bellavia2015strong}                                                                        &                             &                                        &                                      & $2$              & \checkmark        \\\cmidrule{1-7}
      \cref{eq:subproblem_trustregion}                                                  & \citep{cartis2011evaluation}                                                                      &                             & $O(\epsilon^{-2})$                     &                                      &                  &                   \\\cmidrule{2-7}
                                                                                        & \citep{zhang2003new,fan2006convergencerate}                                                       &                             &                                        &                                      & $< 2$            &                   \\\cmidrule{2-7}
                                                                                        & \citep{fan2011improved}                                                                           &                             &                                        &                                      & $2$              &                   \\
                                                                                        \bottomrule
    \end{tabular}
    \begin{tablenotes}
      \item
      \footnotemark[1]%
      The complexity analysis in \citep{bergou2020convergence} assumes the iterates not to converge to a zero-residual solution.
      If the solution sequence converges to a zero-residual solution, then $\bar f$ defined in \citep[Section~3]{bergou2020convergence} is $\bar f = 0$.
      Then, $\mu_{\max}$ defined in \citep[Lemma~3.2]{bergou2020convergence} becomes $\mu_{\max} = \Theta(\epsilon^{-2})$, resulting in the iteration complexity of $O(\epsilon^{-4} \log \epsilon^{-1})$.
      \item
      \footnotemark[2]%
      The complexity analysis in \citep{nesterov2007modified} assumes that $\rank J(x) = n$ for all $x$, which is quite restrictive because such an assumption implies that all stationary points are global optima.
      The local convergence analysis in \citep{nesterov2007modified} assumes that the solution sequence $(x_k)$ is in the neighborhood of a solution $x^*$ such that $F(x^*) = \0$ and $\rank J(x^*) = n$.
    \end{tablenotes}
  \end{threeparttable}
\end{table}

\subsubsection{Subproblems}
Most algorithms for the nonlinear least-squares problem~\cref{eq:mainproblem} generate a solution sequence $(x_k)_{k \in \N}$ by repeatedly solving convex subproblems, and we focus on such algorithms.
There are three popular subproblems, in addition to the LM subproblem~\cref{eq:subproblem}:
\begin{alignat}{2}
  \min_{x \in \R^d}\
  &\norm{F_k + J_k(x - x_k)} + \frac{\lambda}{2} \norm{x - x_k}^2\quad
  &&\subjectto \quad
  x \in \mathcal C,
  \label{eq:subproblem_1_2}\\
  \min_{x \in \R^d}\
  &\norm{F_k + J_k(x - x_k)}^2 + \frac{\lambda}{2} \norm{x - x_k}^3\quad
  &&\subjectto \quad
  x \in \mathcal C,
  \label{eq:subproblem_2_3}\\
  \min_{x \in \R^d}\
  &\norm{F_k + J_k(x - x_k)}^2\quad
  &&\subjectto \quad
  x \in \mathcal C, \  \norm*{x - x_k} \leq \Delta,
  \label{eq:subproblem_trustregion}
\end{alignat}
where $\lambda, \Delta > 0$ are properly defined constants.
Methods using subproblems \cref{eq:subproblem_1_2}, \cref{eq:subproblem_2_3}, and \cref{eq:subproblem_trustregion} have been proposed and analyzed in
\citep{bellavia2010convergence,cartis2020strong,nesterov2007modified,bellavia2015strong},
\citep{bellavia2015strong}, and
\citep{cartis2011evaluation,fan2006convergencerate,fan2011improved,zhang2003new},
respectively.
Other works~\citep{cartis2011adaptive,cartis2011adaptive2,cartis2012adaptive} propose methods with a more general version of \cref{eq:subproblem_2_3}.
These four subproblems~\cref{eq:subproblem,eq:subproblem_1_2,eq:subproblem_2_3,eq:subproblem_trustregion} are closely related in theory; one subproblem becomes equivalent to the others with specific choices of the parameters $\lambda$ and $\Delta$.

In practice, these four subproblems are quite different, and the LM subproblem~\cref{eq:subproblem} is the most tractable one because the objective function $m^k_\lambda$ is smooth and strongly convex.
Thanks to smoothness and strong convexity, we can efficiently solve subproblem~\cref{eq:subproblem} with linearly convergent methods such as the projected gradient method.
Note that the objective function of \cref{eq:subproblem_1_2} is nonsmooth, and \cref{eq:subproblem_2_3,eq:subproblem_trustregion} are not necessarily strongly convex.
Although some algorithms for subproblems~\cref{eq:subproblem_1_2,eq:subproblem_2_3,eq:subproblem_trustregion} \emph{without} constraints have been proposed \citep{bellavia2010convergence,cartis2011adaptive,zhang2003new}, 
efficient algorithms are nontrivial under the presence of constraints.
Hence, the LM method is more practical than methods using other subproblems.

\subsubsection{Complexity for finding a stationary point}
For unconstrained zero-residual problems, \citet{nesterov2007modified} proposed a method with subproblem~\cref{eq:subproblem_1_2} and showed that the method finds an $\epsilon$-stationary point after $O(\epsilon^{-2})$ iterations under a strong assumption (see footnote 2 of \cref{table: comparison of LM methods} for details).
After that, for unconstrained (possibly) nonzero-residual problems, several methods with subproblems~\cref{eq:subproblem,eq:subproblem_1_2,eq:subproblem_trustregion} have been proposed \citep{ueda2010global,zhao2016global,cartis2011evaluation}, and they achieve the same iteration complexity bound under weaker assumptions such as the Lipschitz continuity of $J$ or $\nabla f$.
The method of \citep{cartis2011evaluation} has been extended for constrained problems~\citep{cartis2020strong}.\footnote{
  More precisely, \citet{cartis2020strong} proposed a framework for arbitrary-order methods and it includes a method with subproblem~\cref{eq:subproblem_1_2} as a special case.
}
These methods~\citep{nesterov2007modified,ueda2010global,zhao2016global,cartis2011evaluation,cartis2020strong} have the iteration complexity bound, but computational complexity per iteration, i.e., complexity for a subproblem, is unclear.

The key to bounding complexity per iteration is that we do not need to solve subproblems so accurately to derive the iterative complexity bound.
Several algorithms have been proposed based on this fact for both unconstrained \citep{bellavia2018levenberg,bergou2020convergence,cartis2011adaptive,cartis2011adaptive2} and constrained \citep{cartis2012adaptive} problems.
They use a point that decreases the model function value sufficiently compared to the value at the current iterate $x_k$.
Such a point can be computed with an $\epsilon$-independent number of basic operations: evaluation of $F(x)$, Jacobian-vector multiplications $J(x) u$ and $J(x)^\top v$, and projection onto $\mathcal C$.
Thus, the methods in \citep{cartis2012adaptive,bellavia2018levenberg,cartis2011adaptive,cartis2011adaptive2} achieve the overall complexity $O(\epsilon^{-2}) \times O(1) = O(\epsilon^{-2})$.

Our LM method also finds an $\epsilon$-stationary point within $O(\epsilon^{-2})$ iterations, and the complexity per iteration is $O(1)$ when subproblems are solved approximately like \citep{bellavia2018levenberg,bergou2020convergence,cartis2011adaptive,cartis2011adaptive2,cartis2012adaptive}.
Thus, the overall complexity amounts to $O(\epsilon^{-2})$ same as \citep{cartis2012adaptive,bellavia2018levenberg,cartis2011adaptive,cartis2011adaptive2}.

\subsubsection{Local superlinear convergence}
For unconstrained zero-residual problems, many methods with subproblems~\cref{eq:subproblem,eq:subproblem_1_2,eq:subproblem_2_3,eq:subproblem_trustregion} have achieved local quadratic convergence under a local error bound condition~\citep{yamashita2001rate,fan2001convergence,fan2005quadratic,fan2003modified,fan2006convergence,bellavia2015strong,fan2011improved,dan2002convergence,fan2004inexact,fischer2010inexactness,fan2011convergence}.
These local convergence results have been extended to constrained problems~\citep{kanzow2004levenberg,fan2013levenberg,behling2012unified,facchinei2013family}.
Some methods~\citep{zhang2003new,fan2006convergencerate} have local convergence of an arbitrarily order less than 2.
Other methods~\citep{fan2012modified,fan2014accelerating,fan2015modified} achieve local (nearly) cubic convergence by solving two subproblems in one iteration.
We note that the local convergence analyses in \citep{nesterov2007modified,bellavia2010convergence} assume the solution sequence $(x_k)$ is in the neighborhood of a solution $x^*$ such that $F(x^*) = \0$ and $\rank J(x^*) = n$, which is a stronger assumption than the local error bound.

Among these methods, some~\citep{dan2002convergence,fan2004inexact,fischer2010inexactness,fan2011convergence,behling2012unified,facchinei2013family,bellavia2015strong,bellavia2010convergence} use an approximate solution to subproblems while preserving local quadratic convergence.
The approximate solution
is more accurate than that used to derive the global complexity mentioned in the previous section.
We also use the same kind of approximate solution as \citep{dan2002convergence,fan2004inexact,fischer2010inexactness,fan2011convergence,behling2012unified} to prove local quadratic convergence.
See \cref{cond:local_convergence} in \cref{sec:generalization} for the details of the approximate solution.

\section{Majorization lemma and proposed method}
\label{sec:proposed_method}

Here, we will prove a \emph{majorization lemma} that shows that the LM model $m^k_\lambda$ defined in \cref{eq: def of model} is an upper bound on the objective function. 
In view of this lemma, we can characterize our LM method as a majorization-minimization (MM) method.

For $a, b \in \R^d$, we denote the sublevel set and the line segment by
\begin{align}
 \mathcal S(a)
 &\coloneqq
 \Set{x \in \R^d}{f(x) \leq f(a)},
 \label{eq: def of S}\\
 \mathcal L(a, b)
 &\coloneqq
 \Set{ (1-\theta) a + \theta b \in \R^d }{ \theta \in [0, 1]}.
 \label{eq: def of L}
\end{align}

\subsection{LM method as majorization-minimization}
MM is a framework for nonconvex optimization that successively performs (approximate) minimization of an upper bound on the objective function. 
The following lemma, a majorization lemma, shows that the model $m^k_\lambda$ defined in \cref{eq: def of model} is an upper bound on the objective $f$ over some region under certain assumptions.

\begin{lemma}
  \label{lem: majorization LM}
  Let $\mathcal X \subseteq \R^d$ be any closed convex set, and suppose $x_k \in \mathcal X$.
  Moreover, assume that for some constant $L > 0$,
  \begin{equation}
    \|J(y) - J(x)\| \leq L\|y - x\|,\quad
    \forall x,y \in \mathcal X \ \ \text{s.t.} \ \
    \mathcal L(x, y) \subseteq \mathcal S(x_k).
    \label{eq: J Lip on X}
  \end{equation}
  Then for any $\lambda > 0$ and $x \in \mathcal X$ such that
  \begin{align}
  &
  \lambda \geq L\|F_k\|\quad\text{and}
  \label{eq: lam >= L||F_k||}\\
  &
  m^k_\lambda(x)
  \leq m^k_\lambda(x_k),
  \label{eq: m^k(x) <= m^k(x_k)}
  \end{align}
  the following bound holds:
  \begin{equation}
    f(x)
    \leq
    m^k_\lambda(x).
    \label{eq: f(x_k+u) <= m(u)}
  \end{equation}
\end{lemma}
The proof is given in \cref{sec: proof of MM lemma}.

The assumption in \cref{eq: J Lip on X} is the Lipschitz continuity of $J$ and is analogous to the Lipschitz continuity of $\nabla f$, which is often used in the analysis of first-order methods.
\Cref{eq: lam >= L||F_k||} requires a sufficiently large damping parameter, which corresponds to a sufficiently small step-size for first-order methods.
\Cref{eq: m^k(x) <= m^k(x_k)} requires the point $x \in \mathcal X$ to be a solution that is at least as good as the current point $x_k \in \mathcal X$ in terms of the model function value.

\subsection{Proposed LM method}
\label{sec: proposed LM method}
\begin{algorithm}[t]
  \small
  \caption{Proposed LM method for solving \cref{eq:mainproblem}}
  \label{alg: proposed LM-GD}
  \begin{algorithmic}[1]
    \Require{
      \pbox[t]{40em}{
        $x_0 \in \mathcal C$;
        $M_0, \eta_0 > 0$; $\alpha, \alpha_{\mathrm{in}} > 1$;
        $T \in \Z_{> 0} \cup \set{\infty}$;
        $c > 0$
      }
      \vspace{0.2\baselineskip}
    }
    \State{$M \gets M_0$, $\eta \gets \eta_0$, $k \gets 0$}
    \label{line-LMGD:initialization}
    \Repeat
    \Comment{outer loop}
      \State{$\lambda \gets M \|F_k\|$, $x_{k,0} \gets x_k$, $t \gets 0$}
      \label{line-LMGD:set_lam}
      \Repeat
      \Comment{inner loop for approximately solving subproblem~\eqref{eq:subproblem}}
      \label{line-LMGD:inner loop}
        \State{$y \gets \proj_{\mathcal C}(x_{k,t} - \frac{1}{\eta} \nabla m^k_\lambda(x_{k,t}))$}
        \label{line-LMGD:gradient-descent-projection}
        \If{$m^k_\lambda(y) \leq m^k_\lambda(x_{k,t}) + \inner{\nabla m^k_\lambda(x_{k,t})}{y - x_{k,t}} + \frac{\eta}{2}\|y - x_{k,t}\|^2$}
        \label{line-LMGD:check-eta}
          \State{$x_{k,t+1} \gets y$, $t \gets t+1$}
          \label{line-LMGD:update_x_y}
        \Else
        \label{line-LMGD:bad_eta}
          \State{$\eta \gets \alpha_{\mathrm{in}} \eta$}
          \label{line-LMGD:reset_eta}
        \EndIf
      \Until{($t = T$) \textbf{or} ($t \geq 1$ \textbf{and} $x_{k,t}$ is a $(\cstSIL \lambda \|F_k\|)$-stationary point of subproblem~\eqref{eq:subproblem})}
      \label{line-LMGD:inner loop_until}
      \State{$x \gets x_{k,t}$}
      \label{line-LMGD:set_x_xkt}
      \If{$f(x) \leq m^k_\lambda(x)$}
      \label{line-LMGD:check_M}
        \State{$(x_{k+1}, \lambda_k) \coloneqq (x, \lambda)$, $k \gets k+1$}
        \Comment{successful}
        \label{line-LMGD:successful}
      \Else
        \label{line-LMGD:small_M}
        \State{$M \gets \alpha M$}
        \Comment{unsuccessful}
        \label{line-LMGD:unsuccessful}
      \EndIf
    \Until{a solution with a desired accuracy is obtained}
  \end{algorithmic}
\end{algorithm}

Based on \cref{lem: majorization LM}, we propose an LM method that solves problem~\eqref{eq:mainproblem}.
The proposed LM is formally described in \cref{alg: proposed LM-GD} and is outlined below.
First, in \cref{line-LMGD:initialization}, three parameters are initialized: an estimate $M$ of the Lipschitz constant $L$ of $J$, a parameter $\eta$ used for solving subproblems, and the iteration counter $k$.
\cref{line-LMGD:set_lam} sets $\lambda$ using $M$ as an estimate of $L$ based on \cref{eq: lam >= L||F_k||}.
Then, the inner loop of \cref{line-LMGD:inner loop,line-LMGD:gradient-descent-projection,line-LMGD:check-eta,line-LMGD:update_x_y,line-LMGD:bad_eta,line-LMGD:reset_eta,line-LMGD:inner loop_until} solves subproblem~\eqref{eq:subproblem} approximately by a projected gradient method.
The details of the inner loop will be described later.
\cref{line-LMGD:check_M,line-LMGD:successful,line-LMGD:small_M,line-LMGD:unsuccessful} check if the current $\lambda$ and the computed solution $x$ are acceptable.
If $\lambda$ and $x$ satisfy \cref{eq: f(x_k+u) <= m(u)}, they are accepted as $\lambda_k$ and $x_{k+1}$.
Otherwise, the current value of $M$ is judged to be small as an estimate of $L$ in light of \cref{lem: majorization LM} and is increased.
We refer to the former case as a ``successful'' iteration and the latter as an ``unsuccessful'' iteration.
Note that $k$ represents not the number of outer iterations but that of only successful iterations.
As shown later in \cref{item-lem: unsuccessful upper bound global,item-thm: L_hat bound local}, the number of unsuccessful iterations is upper-bounded by a constant under certain assumptions.

\myparagraph{Inner loop for subproblem}
In the inner loop of \cref{line-LMGD:inner loop,line-LMGD:gradient-descent-projection,line-LMGD:check-eta,line-LMGD:update_x_y,line-LMGD:bad_eta,line-LMGD:reset_eta,line-LMGD:inner loop_until}, subproblem~\eqref{eq:subproblem} is solved approximately by the projected gradient method.
Here, the operator $\proj_{\mathcal C}$ in \cref{line-LMGD:gradient-descent-projection} is the projection operator defined by 
\begin{align}
  \proj_{\mathcal C}(x) \coloneqq \argmin_{y \in \mathcal C}\|y - x\|.
\end{align}
The parameter $t$ is the inner iteration counter, and the parameter $\eta$ is the inverse step-size that is adaptively chosen by a standard backtracking technique in \cref{line-LMGD:check-eta,line-LMGD:update_x_y,line-LMGD:bad_eta,line-LMGD:reset_eta}.
As shown in \cref{item-lem: inner unsuccessful upper bound global} later, \cref{line-LMGD:reset_eta} is executed a finite number of times under certain standard assumptions.
Hence, the inner loop must stop after a finite number of iterations.

\myparagraph{Input parameters}
\cref{alg: proposed LM-GD} has several input parameters.
The parameters $M_0$ and $\alpha$ are used to estimate the Lipschitz constant of the Jacobian $J$, and the parameters $\eta_0$ and $\alpha_{\mathrm{in}}$ are used to control the step-size in the inner loop.
The parameters $T$ and $c$ control how accurately the subproblems are solved through the stopping criteria of the inner loop.
Here, note that we allow for $T = \infty$.
As we will prove in \cref{sec: global convergence analysis}, the algorithm has an iteration complexity bound for an $\epsilon$-stationary point regardless of the choice of the input parameters.
However, to obtain an overall complexity bound or local quadratic convergence, there are restrictions on the choice of $T$, as explained in the next paragraph.

\myparagraph{Stopping criteria for inner loop}
There are two types of stopping criteria as in \cref{line-LMGD:inner loop_until}, and the inner loop terminates when at least one of them is satisfied.
If $T < \infty$, the projected gradient method stops after executing \cref{line-LMGD:update_x_y} at most $T$ times, and then the overall complexity for an $\epsilon$-stationary point is guaranteed to be $O(\epsilon^{-2})$.
If $T = \infty$, we have to solve subproblems more accurately to find a $(\cstSIL \lambda \|F_k\|)$-stationary point of the subproblem, and then \cref{alg: proposed LM-GD} achieves local quadratic convergence.




\begin{remark}
  \label{rem: additional input parameter}
  To make the algorithm more practical, we can introduce other parameters $0 < \beta < 1$ and $M_{\min} > 0$, and update $M \gets \max\set{\beta M, M_{\min}}$ after every successful iteration. As with the gradient descent method, such an operation prevents the estimate $M$ from being too large and eliminates the need to choose $M_0$ carefully.
  Inserting this operation never deteriorates the complexity bounds described in \cref{sec: global convergence analysis} and the local quadratic convergence in \cref{sec: local convergence analysis}.
\end{remark}

\begin{remark}
  Some methods (e.g., \citep{ueda2010global,zhao2016global}) use the condition
  \begin{align}
    \frac{m^k_\lambda(x_k) - f(x)}{m^k_\lambda(x_k) - m^k_\lambda(x)}
    \geq
    \theta
    \label{eq:trustregion_cond_theta}
  \end{align}
  with some $0 < \theta < 1$ to determine whether the computed solution $x$ to the subproblem is acceptable.
  Our acceptance condition~\cref{eq: f(x_k+u) <= m(u)} is stronger than the classical one since \cref{eq: f(x_k+u) <= m(u)} is equivalent to
  \begin{align}
    \frac{m^k_\lambda(x_k) - f(x)}{m^k_\lambda(x_k) - m^k_\lambda(x)}
    \geq
    1
  \end{align}
  under condition~\cref{eq: m^k(x) <= m^k(x_k)}.
  Therefore, \cref{lem: majorization LM} is stronger than the classical statement that condition~\cref{eq:trustregion_cond_theta} holds if $\lambda$ is sufficiently large.
\end{remark}

\section{Iteration complexity and overall complexity}
\label{sec: global convergence analysis}

We will prove that \cref{alg: proposed LM-GD} finds an $\epsilon$-stationary point 
of problem~\eqref{eq:mainproblem} within $O(\epsilon^{-2})$ outer iterations.
Futhermore, we will prove that under $T < \infty$, the overall complexity for an $\epsilon$-stationary point is also $O(\epsilon^{-2})$.
Throughout this section, $(x_k)$ and $(\lambda_k)$ denote the sequences generated by the algorithm.


\subsection{Assumptions}
\label{sec: assumption global}
We make the following assumptions to derive the complexity bound.
Recall that the sublevel set $\mathcal S(a)$ and the line segment $\mathcal L(a, b)$ are defined in \cref{eq: def of S,eq: def of L} and that $x_0 \in \mathcal C$ denotes the starting point of \cref{alg: proposed LM-GD}.

\begin{assumption}
 \label{asm: for global convergence}
For some constants $\sigma, L > 0$,
 \leavevmode
 \begin{enuminasm}
 \item
 \label{item-asm: J bounded on S}
 $\|J(x)\| \leq \sigma$, \
 $\forall x \in \mathcal C \cap \mathcal S(x_0)$,
 \item
 \label{item-asm: J Lip on S}
 $\|J(y) - J(x)\| \leq L\|y - x\|$, \
 $\forall x,y \in \mathcal C$ \ s.t.~\
 $\mathcal L(x, y) \subseteq \mathcal S(x_0)$.
 \end{enuminasm}
\end{assumption}

\cref{item-asm: J bounded on S} means the $\sigma$-boundedness of $J$ on $\mathcal C \cap \mathcal S(x_0)$.
\cref{item-asm: J Lip on S} is similar to the $L$-Lipschitz continuity of $J$ on $\mathcal C \cap \mathcal S(x_0)$ but weaker due to the condition of $\mathcal L(x, y) \subseteq \mathcal S(x_0)$.
\cref{asm: for global convergence} is milder than the assumptions in the previous work that discussed the iteration complexity, even when $\mathcal C = \R^d$.
For example, the analysis in \citep{zhao2016global} assumes $f$ and $J$ to be Lipschitz continuous on $\R^d$, which implies the boundedness of $J$ on $\R^d$.

\subsection{Approximate stationary point}
\label{sec:appropriate_stationary_point}
Before analyzing the algorithm, we define an $\epsilon$-stationary point for constrained optimization problems.
Let $\iota_{\mathcal C}: \R^d \to \R \cup \set{+\infty}$ be the indicator function of the closed convex set $\mathcal C \subseteq \R^d$.
For a convex function $g: \R^d \to \R \cup \set{+\infty}$, its subdifferential at $x \in \R^d$ is the set defined by $\partial g(x) \coloneqq \Set{p \in \R^d}{g(y) \geq g(x) + \inner{p}{y - x}, \ \forall y \in \R^d}$.

\begin{definition}[see, e.g., Definition 1 in \citep{nesterov2013gradient}]
  \label{def:eps-stationary-point}
  For $\epsilon > 0$, a point $x \in \mathcal C$ is said to be an \emph{$\epsilon$-stationary point} of the problem $\min_{x \in \mathcal C} f(x)$ if
  \begin{equation}
    \min_{p \in \partial \iota_{\mathcal C}(x)} \norm*{\nabla f(x) + p}
    \leq \epsilon.
    \label{eq:def_eps-stationary-point}
  \end{equation}
\end{definition}
This definition is consistent with the unconstrained case; the above inequalities are equivalent to $\|\nabla f(x)\| \leq \epsilon$ when $\mathcal C = \R^d$.
There is another equivalent definition of an $\epsilon$-stationary point, which we will also use.
\begin{lemma}
  For $x \in \mathcal C$ and $\epsilon > 0$, condition \cref{eq:def_eps-stationary-point} is equivalent to
  \begin{align}
    \inner{\nabla f(x)}{y - x} \geq - \epsilon \norm*{y - x},
    \quad
    \forall y \in \mathcal C.
    \label{eq:property_eps-stationary-point}
  \end{align}
\end{lemma}
\begin{proof}
  The tangent cone $\mathcal T(x)$ of $\mathcal C$ at $x \in \mathcal C$ is defined by
  \begin{align}
    \mathcal T(x)
    &\coloneqq
    \Set{\beta (y - x)}{y \in \mathcal C, \ \beta \geq 0}.
    \label{eq:def_tangent_cone}
  \end{align}
  Note that
  \begin{align}
    \mathcal T(x) &= \Set{z \in \R^d}{\inner{y}{z} \leq 0,\ \forall y \in \partial \iota_{\mathcal C}(x)}
    \label{eq:duality_TandN}
  \end{align}
  because $\mathcal C$ is a closed convex set and $\partial \iota_C(x)$ is the normal cone of $\mathcal C$.
  We have
  \begin{alignat}{2}
    \min_{p \in \partial \iota_{\mathcal C}(x)} \norm*{\nabla f(x) + p}
    &=
    \min_{p \in \partial \iota_{\mathcal C}(x)} \max_{u: \norm*{u} \leq 1}
    \inner*{- \nabla f(x) - p}{u}\\
    &=
    \max_{u: \norm*{u} \leq 1} \inf_{p \in \partial \iota_{\mathcal C}(x)}
    \set[\Big]{
      \inner*{- \nabla f(x)}{u} - \inner*{p}{u}
    }
    &\quad&\by{a minimax theorem}\\
    &=
    \max_{u \in \mathcal T(x),\,\norm*{u} \leq 1} \inner*{- \nabla f(x)}{u}
    &\quad&\by{\cref{eq:duality_TandN}}\\
    &=
    \sup_{y \in \mathcal C \setminus \set{x}} 
    \frac{\inner*{- \nabla f(x)}{y - x}}{\norm*{y - x}}
    &\quad&\by{\cref{eq:def_tangent_cone}}.
  \end{alignat}
  Therefore, condition~\cref{eq:def_eps-stationary-point} is equivalent to
  \begin{align}
    \sup_{y \in \mathcal C \setminus \set{x}} 
    \frac{\inner*{- \nabla f(x)}{y - x}}{\norm*{y - x}}
    \leq \epsilon,
  \end{align}
  which is also equivalent to \cref{eq:property_eps-stationary-point}.
\end{proof}

A useful tool for deriving iteration complexity bounds is \emph{gradient mapping} (see, e.g., \citep{nesterov2004introductory}), also known as \emph{projected gradient}~\citep{lan2013iteration} or \emph{reduced gradient}~\citep{nesterov2018lectures}.
For $\eta > 0$, the projected gradient operator $\mathcal P_\eta: \mathcal C \to \mathcal C$ and the gradient mapping $\mathcal G_\eta : \mathcal C \to \R^d$ for problem~\eqref{eq:mainproblem} are defined by
\begin{align}
  \mathcal P_\eta(x)
  &\coloneqq
  \argmin_{y \in \mathcal C} \Big\{
  \inner{\nabla f(x)}{y - x} + \frac{\eta}{2}\|y - x\|^2
  \Big\}
  =
  \proj_{\mathcal C} \Big( x - \frac{1}{\eta} \nabla f(x) \Big),
  \label{eq: def of P_eta}\\
  \mathcal G_\eta(x)
  &\coloneqq
  \eta(x - \mathcal P_\eta(x)).
  \label{eq: def of G_eta}
\end{align}
The following lemma shows the relationship between an $\epsilon$-stationary point and the gradient mapping.
\begin{lemma}
  \label{lem: small grad map implies approx stationary point}
  Suppose that \cref{asm: for global convergence} holds, and let
  \begin{equation}
    L_{f} \coloneqq \sigma^2 + L \|F_0\|.
    \label{eq: def L_f}
  \end{equation}
  Then, for any $x \in \mathcal C \cap \mathcal S(x_0)$ and $\eta \geq L_f$, the point $\mathcal P_\eta(x)$ is a $(2 \|\mathcal G_\eta(x)\|)$-stationary point of problem~\cref{eq:mainproblem}.
\end{lemma}
The proof is given in \cref{sec: proof of lemma GMtoASP}. 
This lemma will be used for the proof of \cref{item-thm:successful_iter_global}.

Although \cref{lem: small grad map implies approx stationary point} looks quite similar to \citep[Corollary 1]{nesterov2013gradient}, there exists a significant difference in their assumptions.
Indeed, \cref{lem: small grad map implies approx stationary point} assumes the boundedness and the Lipschitz property of $J$ only on a (possibly) nonconvex set $\mathcal C \cap \mathcal S(x_0)$, whereas \citep[Corollary 1]{nesterov2013gradient} assumes the Lipschitz continuity on the whole space $\R^d$.
This makes our proof more complicated than in \cite[Corollary 1]{nesterov2013gradient}.

\subsection{Preliminary lemmas}
First, we bound the decrease in the model function value due to the inner loop.
For $\eta > 0$, we define the function $\mathcal D_\eta: \mathcal C \to \R$ by
\begin{equation}
  \mathcal D_\eta(x)
  \coloneqq
  - \min_{y \in \mathcal C} \Big\{
  \inner{\nabla f(x)}{y - x} + \frac{\eta}{2}\|y - x\|^2
  \Big\}.
  \label{eq: def of D_eta}
\end{equation}
We see that $
\mathcal D_\eta(x)
\geq -\inner{\nabla f(x)}{x - x} - \frac{\eta}{2}\|x - x\|^2
= 0$
for all $x \in \mathcal C$.
In addition, $\mathcal D_\eta(x)$ is decreasing with respect to $\eta$.
\begin{lemma}
  \label{lem:inner_decrease}
  The solution $x$ obtained in \cref{line-LMGD:set_x_xkt} of \cref{alg: proposed LM-GD} satisfies
  \begin{equation}
    m^k_{\lambda} (x)
    \leq
    m^k_{\lambda} (x_k) - \mathcal D_\eta(x_k)
    \leq 
    m^k_{\lambda} (x_k),
    \label{eq: asm sufficient decrease}
  \end{equation}
  where $k$, $\lambda$, and $\eta$ are parameters in \cref{alg: proposed LM-GD}.
\end{lemma}
\begin{proof}
  The second inequality in \cref{eq: asm sufficient decrease} follows from the nonnegativity of $\mathcal D_\eta(x)$, and therefore we will prove the first one.
  Let $T'$ denote the value of $t$ when the inner loop is completed, and for each $0 \leq t \leq T'$, let $\eta_{k,t}$ denote the values of $\eta$ when $x_{k,t}$ is obtained through \cref{line-LMGD:update_x_y}.
  Our aim is to prove the first inequality in \cref{eq: asm sufficient decrease} with $(x, \eta)=(x_{k,T'},\eta_{k,T'})$. 
  We have
  \begin{alignat}{2}
    m^k_{\lambda}(x_{k,1})
    &\leq
    m^k_\lambda(x_k) + \inner{\nabla m^k_\lambda(x_k)}{x_{k,1} - x_k} + \frac{\eta_{k,1}}{2} \norm{x_{k,1} - x_k}^2
    &\quad&\by{\cref{line-LMGD:check-eta}}\\
    &=
    m^k_\lambda(x_k) + \min_{z \in \mathcal C} \set[\Big]{
      \inner{\nabla m^k_\lambda(x_k)}{z - x_k} + \frac{\eta_{k,1}}{2}\norm{z - x_k}^2
    }
    &\quad&\by{the definition of $x_{k,1}$}\\
    &= m^k_\lambda(x_k) - \mathcal D_{\eta_{k,1}}(x_k)
    &\quad&\by{$\nabla m^k_\lambda(x_k) = \nabla f(x_k)$}.
  \end{alignat}
  Since $\mathcal D_\eta(x_k)$ is decreasing in $\eta$ and $\eta_{k,1} \leq \eta_{k,2} \leq \dots \leq \eta_{k,T'}$, we have $\mathcal D_{\eta_{k,1}}(x_k) \geq \mathcal D_{\eta_{k,T'}}(x_k)$.
  On the other hand, we have $m^k_\lambda(x_{k,1}) \geq \dots \ge  m^k_\lambda(x_{k, T'})$.
  Combining these inequalities, we obtain the desired result.
\end{proof}

From the above lemma and \cref{line-LMGD:check_M}, it follows that for all $k$,
\begin{equation}
  f(x_{k+1})
  \leq m^k_{\lambda_k}(x_{k+1})
  \leq m^k_{\lambda_k}(x_k)
  = f(x_k).
  \label{eq: f(x_(k+1)) <= m(x_(k+1)) <= m(x_k) = f(x_k)}
\end{equation}
This monotonicity of $f(x_k)$ in $k$ is an important property of the majorization-minimization and will be used in our analysis.

The following two lemmas show that the parameters $M$ and $\eta$ in the algorithm are upper-bounded, and hence \cref{line-LMGD:reset_eta,line-LMGD:unsuccessful} are executed only a finite number of times per single run.
\begin{lemma}
  \label{lem: unsuccessful global}
  Suppose that \cref{item-asm: J Lip on S} holds, and let
  \begin{equation}
    \bar M \coloneqq \max\{M_0, \alpha L\},
    \label{eq: def of bar L global}
  \end{equation}
  where $M_0$ and $\alpha$ are the inputs of \cref{alg: proposed LM-GD}. Then,
  \begin{enuminlem}
    \item
    \label{item-lem: hat L upper bound global}
    the parameter $M$ in \cref{alg: proposed LM-GD} always satisfies $M \leq \bar M$;
    \item
    \label{item-lem: unsuccessful upper bound global}
    throughout the algorithm, the number of unsuccessful iterations is at most $\ceil{\log_\alpha (\bar M / M_0)} = O(1)$.
  \end{enuminlem}
\end{lemma}
\begin{proof}
  We have $\mathcal S(x_k) \subseteq \mathcal S(x_0)$ from \cref{eq: f(x_(k+1)) <= m(x_(k+1)) <= m(x_k) = f(x_k)}, and therefore \cref{item-asm: J Lip on S} implies \cref{eq: J Lip on X} with $\mathcal X = \mathcal C$.
  On the other hand, \cref{eq: asm sufficient decrease} directly implies \cref{eq: m^k(x) <= m^k(x_k)}.
  Hence, by \cref{lem: majorization LM} with $\mathcal X = \mathcal C$ and \cref{lem:inner_decrease}, if $M \geq L$ holds at \cref{line-LMGD:set_lam}, the condition in \cref{line-LMGD:check_M} must be true.
  Therefore, if $M_0 \geq L$, no unsuccessful iterations occur and the parameter $M$ always satisfies $M = M_0$. Otherwise, there exists an integer $l \geq 1$ such that $L \leq \alpha^l M_0 < \alpha L$. Since $M = \alpha^l M_0$ after $l$ unsuccessful iterations, the parameter $M$ always satisfies $M < \alpha L$. Consequently, we obtain the first result, and the second follows from the first.
\end{proof}

\begin{lemma}
  \label{lem:inner_global}
  Suppose that \cref{asm: for global convergence} holds, and let 
  \begin{align}
    \bar \eta \coloneqq \max\set{\eta_0, \alpha_{\mathrm{in}} (\sigma^2 + \bar M \|F_0\|) },
    \label{eq:def_bar_eta}
  \end{align}
  where $\eta_0$ and $\alpha_{\mathrm{in}}$ are the inputs of \cref{alg: proposed LM-GD} and $\bar M$ is defined in \cref{eq: def of bar L global}.
  Then,
  \begin{enuminlem}
    \item
    \label{item-lem: eta upper bound global}
    the parameter $\eta$ in \cref{alg: proposed LM-GD} always satisfies $\eta \leq \bar \eta$;
    \item
    \label{item-lem: inner unsuccessful upper bound global}
    throughout the algorithm, \cref{line-LMGD:reset_eta} will be executed at most $\ceil{\log_{\alpha_\mathrm{in}} (\bar \eta / \eta_0)} = O(1)$ times.
  \end{enuminlem}
\end{lemma}
\begin{proof}
  Since the function $m^k_\lambda$ defined by \cref{eq: def of model} has the $(\|J_k\|^2 + \lambda)$-Lipschitz continuous gradient, we have
  \begin{align}
    m^k_\lambda(y)
    \leq
    m^k_\lambda(x)
    + \inner{\nabla m^k_\lambda(x)}{y - x}
    + \frac{\norm*{J_k}^2 + \lambda}{2} \norm*{y - x}^2,\quad
    \forall x, y \in \R^d
  \end{align}
  (see, e.g., \citep[Eq.~(2.1.9)]{nesterov2018lectures}).
  We also have $\|J_k\|^2 + \lambda \leq \sigma^2 + \bar M \|F_0\|$ from \cref{item-asm: J bounded on S,lem: unsuccessful global}.
  Therefore, the inequality in \cref{line-LMGD:check-eta} must hold if $\eta \geq \sigma^2 + \bar M \|F_0\|$.
  With the same arguments as in \cref{lem: unsuccessful global}, we obtain the desired results.
\end{proof}

As we can see from the proofs of \cref{lem: unsuccessful global,lem:inner_global}, if $M_0 \geq L$ and $\eta_0 \geq \sigma^2 + M_0 \norm{F_0}$, then no unsuccessful iterations occur in both outer and inner loops.
Adjusting $M$ and $\eta$ adaptively as in the presented algorithm avoids a too small step-size in practice.

\subsection{Iteration complexity and overall complexity}
\label{sec: global sublinear convergence}
We use the following lemma for the analysis.
\begin{lemma}
  \label{lem: property of D and G}
  \begin{equation}
  \mathcal D_\eta(x)
  \geq
  \frac{1}{2\eta}\|\mathcal G_\eta(x)\|^2,\quad
  \forall x \in \mathcal C,
  \quad
  \eta > 0.
  \end{equation}
\end{lemma}
\begin{proof}
  By the first-order optimality condition on \cref{eq: def of P_eta} and the convexity of $\mathcal C$, we have
  \begin{equation}
    \inner{\nabla f(x) + \eta (\mathcal P_\eta(x) - x)}{y - \mathcal P_\eta(x)} \geq 0,
    \quad
    \forall y \in \mathcal C.
    \label{eq: opt cond on def of P}
  \end{equation}
  Using this inequality, we obtain
  \begin{alignat}{2}
    \mathcal D_\eta(x)
    &=
    \inner{\nabla f(x)}{x - \mathcal P_\eta(x)}
    - \frac{\eta}{2}\|x - \mathcal P_\eta(x)\|^2
    &\quad&\text{(by \cref{eq: def of P_eta,eq: def of D_eta})}\\
    &\geq
    \frac{\eta}{2}\|x - \mathcal P_\eta(x)\|^2
    &\quad&\text{(by \cref{eq: opt cond on def of P} with $y = x$)}\\
    &=
    \frac{1}{2\eta}\|\mathcal G_\eta(x)\|^2
    &\quad&\text{(by \cref{eq: def of G_eta})}.
    \end{alignat}
\end{proof}

We show the asymptotic global convergence and the iteration complexity bound of \cref{alg: proposed LM-GD}.
\begin{theorem}
  Suppose that \cref{asm: for global convergence} holds, and define $\bar \eta$ by \cref{eq:def_bar_eta}.
  Then,
  \begin{enuminthm}
    \item
    $\displaystyle \lim_{k \to \infty} \norm{\mathcal G_{\bar \eta}(x_k)} = 0$, and therefore, any accumulation point of $(x_k)$ is a stationary point of problem~\cref{eq:mainproblem};
    \item
    \label{item-thm:successful_iter_global}
    $\mathcal P_{\bar \eta}(x_k)$ is an $\epsilon$-stationary point of problem~\cref{eq:mainproblem} for some $k = O(\epsilon^{-2})$.
  \end{enuminthm}
\end{theorem}
\begin{proof}
  We have
  \begin{alignat}{2}
    f(x_{k+1}) - f(x_k)
    \label{eq: f(x_k+1) - f(x_k) <= -||G||^2}
    &\leq
    m^k_{\lambda_k}(x_{k+1}) - m^k_{\lambda_k}(x_k)
    &\quad&\by{\cref{line-LMGD:check_M} and $m^k_{\lambda_k}(x_k) = f(x_k)$}\\
    &\leq
    - \mathcal D_{\bar \eta}(x_k)
    &\quad&\text{(by \cref{lem:inner_decrease,item-lem: eta upper bound global})}\\
    &\leq
    - \frac{1}{2 \bar \eta}\|\mathcal G_{\bar \eta}(x_k)\|^2
    &\quad&\text{(by \cref{lem: property of D and G})}.
  \end{alignat}
  Summing up this inequality for $k=0,1,\dots,K-1$, we obtain
  \begin{equation}
    \sum_{k=0}^{K-1} \|\mathcal G_{\bar \eta}(x_k)\|^2
    \leq
    2 \bar \eta (f(x_0) - f(x_K))
    \leq 2 \bar \eta f(x_0)
    \label{eq:bound_sum_gradmap}
  \end{equation}
  for all $K \geq 0$.
  Therefore, we also have $\sum_{k=0}^{\infty} \|\mathcal G_{\bar \eta}(x_k)\|^2 \leq 2 \bar \eta f(x_0)$,
  yielding $\lim_{k \to \infty} \norm{\mathcal G_{\bar \eta}(x_k)} = 0$, the first result.

  Combining \eqref{eq:bound_sum_gradmap} with $\min_{0 \leq k < K} \|\mathcal G_{\bar \eta}(x_k)\|^2 \leq \frac{1}{K} \sum_{k=0}^{K-1} \|\mathcal G_{\bar \eta}(x_k)\|^2$, we have $\|\mathcal G_{\bar \eta}(x_k)\| \leq \epsilon / 2$ for some $k = O(\epsilon^{-2})$.
  For such $x_k$, the point $\mathcal P_{\bar \eta}(x_k)$ is an $\epsilon$-stationary point from \cref{lem: small grad map implies approx stationary point} and $\bar \eta \geq L_f$.
  Thus, we have obtained the second result.
\end{proof}

From \cref{item-lem: unsuccessful upper bound global,item-thm:successful_iter_global}, we obtain the iteration complexity bound of our algorithm as follows.
\begin{corollary}
  \label{cor:iteration_complexity}
  Under \cref{asm: for global convergence}, \cref{alg: proposed LM-GD} finds an $\epsilon$-stationary point within $O(\epsilon^{-2})$ outer iterations, namely,
  $O(\epsilon^{-2})$ successful and unsuccessful iterations.
\end{corollary}

From this iteration complexity bound and \cref{item-lem: inner unsuccessful upper bound global}, we also obtain the overall complexity bound.
\begin{corollary}\label{cor: global overall complexity}
  Suppose that \cref{asm: for global convergence} holds.
  Then, \cref{alg: proposed LM-GD} with $T < \infty$ finds an $\epsilon$-stationary point after $O(\epsilon^{-2} T)$ basic operations.
\end{corollary}
We use the term \emph{basic operations} to refer to evaluation of $F(x)$, Jacobian-vector multiplications $J(x) u$ and $J(x)^\top v$, and projection onto $\mathcal C$ as in \cref{sec:oracle_model}.

In order to compute an $\epsilon$-stationary point based on \cref{item-thm:successful_iter_global}, knowledge of the value of $\bar \eta$ is required. However, this requirement can be circumvented with a slight modification of the algorithm.
See \cref{subsec:without_etabar} for the details.

\section{Local quadratic convergence}
\label{sec: local convergence analysis}
For \emph{zero-residual} problems, we will prove that the sequence $(x_k)$ 
generated by \cref{alg: proposed LM-GD} with $T = \infty$ converges locally quadratically to an optimal solution. 
Let us denote the set of optimal solutions to problem \cref{eq:mainproblem} by $\mathcal X^* \coloneqq \Set{x \in \mathcal C}{F(x) = \0}$ and the distance between $x \in \R^d$ and $\mathcal X^*$ simply by $\dist(x) \coloneqq \min_{y \in \mathcal X^*} \|y - x\|$.
Throughout this section, we fix a point $x^* \in \mathcal X^*$ and denote a neighborhood of $x^*$ by $\mathcal B(r) \coloneqq \Set{x \in \R^d }{\|x - x^*\| \leq r}$ for $r > 0$.\footnote{
  If $x^*$ is an interior point of the constraint $\mathcal C$, the problem can be regarded as an unconstrained one, and the quadratic convergence is easier to prove.
  We do not assume this, i.e., $x^*$ may be on the boundary of $\mathcal C$.
}
As in the previous section, we denote the sequences generated by \cref{alg: proposed LM-GD} with $T = \infty$ by $(x_k)$ and $(\lambda_k)$.

\subsection{Assumptions}
\label{sec: assumption local}
We make the following assumptions to prove local quadratic convergence.

\begin{assumption}
  \label{asm: local convergence}
  \leavevmode
  \begin{enuminasm}
    \item
    \label{item-asm: zero-residual}
    There exists $x \in \mathcal C$ such that $F(x) = \0$.
    \end{enuminasm}
    For some constants $\cstEB, L, r > 0$,
    \begin{enuminasm}
    \setcounter{enuminasmi}{1}
    \item
    \label{item-asm: EB on N(r)}
    $\cstEB \dist(x) \leq \|F(x)\|$, \ \
    $\forall x \in \mathcal C \cap \mathcal B(r)$,
    \item
    \label{item-asm: J Lip on N(r)}
    $\|J(y) - J(x)\| \leq L\|y - x\|$, \ \
    $\forall x,y \in \mathcal C \cap \mathcal B(r)$.
  \end{enuminasm}
\end{assumption}

\cref{item-asm: zero-residual} requires the problem to be zero-residual, 
\cref{item-asm: EB on N(r)} is called \tred{a} local error bound condition, and 
\cref{item-asm: J Lip on N(r)} is the local Lipschitz continuity of $J$. 
These assumptions are used in the previous analyses of LM methods~\citep{bellavia2015strong, yamashita2001rate, kanzow2004levenberg, dan2002convergence, fan2003modified, fan2004inexact, fan2006convergence, fan2005quadratic, fan2001convergence, behling2012unified, facchinei2013family}.

\subsection{Fundamental inequalities for analysis}

Since $\mathcal C \cap \mathcal B(r)$ is compact, there exists a constant $\sigma > 0$ such that
\begin{equation}
  \|J(x)\| \leq \sigma,
  \quad
  \forall x \in \mathcal C \cap \mathcal B(r),
  \label{eq:J_bounded_local}
\end{equation}
which implies
\begin{equation}
  \|F(y) - F(x)\| \leq \sigma \|y - x\|,
  \quad
  \forall x, y \in \mathcal C \cap \mathcal B(r).
  \label{eq:F_Lip_local}
\end{equation}
Let $\sigma$ denote such a constant in the rest of this section.

For a point $x \in \R^d$, let $\tilde x \in \mathcal X^*$ denote
an optimal solution closest to $x$; $\|\tilde x - x\| = \dist(x)$. 
In particular, $\tilde x_k $ denotes one of the closest solutions to $x_k$ for each $k \geq 0$. 
Since $\|\tilde a - x^*\| \leq \|\tilde a - a\| + \|a - x^*\| \leq 2\|a - x^*\|$, we have
\begin{equation}
  a \in \mathcal B(r/2)
  \quad \Longrightarrow\quad
  \tilde a \in \mathcal B(r).
  \label{eq:a_in_r/2_implies_sol_in_r}
\end{equation}
Therefore, \cref{eq:F_Lip_local} with $y \coloneqq \tilde x$ implies
\begin{equation}
  \|F(x)\| \leq \sigma \|x - \tilde x\| = \sigma \dist (x),
  \quad
  \forall x \in \mathcal C \cap \mathcal B(r/2).
  \label{eq:F_<=_dist}
\end{equation}

From the stopping criterion in \cref{line-LMGD:inner loop_until} of \cref{alg: proposed LM-GD} with $T = \infty$ and \cref{def:eps-stationary-point}, the solution $x$ obtained in \cref{line-LMGD:set_x_xkt} satisfies
\begin{equation}
  \inner{\nabla m^k_\lambda(x)}{y - x} \geq - \cstSIL \lambda \|F_k\| \|y - x\|,
  \quad
  \forall y \in \mathcal C.
  \label{eq:asm_local_subproblem_y_mu}
\end{equation}
From the definition of $x_{k+1}$ and $\lambda_k$, we also have the inequality with $(x, \lambda) = (x_{k+1}, \lambda_k)$, i.e.,
\begin{equation}
  \inner{\nabla m^k_{\lambda_k}(x_{k+1})}{y - x_{k+1}} \geq - \cstSIL \lambda_k \|F_k\| \|y - x_{k+1}\|,
  \quad
  \forall y \in \mathcal C.
  \label{eq:asm_local_subproblem_x_lam}
\end{equation}

\subsection{Preliminary lemma}
\begin{lemma}
  \label{lem:local_induction_step}
  Suppose that \cref{asm: local convergence} holds, and define $\bar M$ by \cref{eq: def of bar L global}.
  Define the constants $C_1, C_2, \delta > 0$ by
  \begin{subequations}
    \begin{align}
      C_1 &\coloneqq \sqrt{1 + \cstSIL^2 \sigma^2 + \frac{L^2 r}{16 \cstEB M_0}},
      \label{eq:def_C1}\\
      C_2
      &\coloneqq \frac{1}{c^2}
      \bigg(
      \sigma^2 \Big( \cstSIL \bar M + \frac{L}{2\cstEB} \Big)
      + \frac{L\sigma C_1^2}{2} + (L + \bar M)\sigma C_1
      \bigg),
      \label{eq:def_C2}\\
      \delta &\coloneqq \frac{r}{2(1 + C_1)},
      \label{eq:def_delta}
    \end{align}
  \end{subequations}
  where $M_0$ and $c$ are the inputs of \cref{alg: proposed LM-GD}.
  Assume that $x_k \in \mathcal B(\delta)$ and $M \leq \bar M$ hold at \cref{line-LMGD:set_lam}.
  Then,
  \begin{enuminlem}
    \item
    \label{item-lem:x_xk_upper_local}
    the solution $x$ obtained in \cref{line-LMGD:set_x_xkt} satisfies
    \begin{align}
      \|x - x_k\| \leq C_1 \dist(x_k);
      \label{eq:y-x_dist_bound}
    \end{align}
    \item
    \label{item-lem:M_upper_local}
    $M \leq \bar M$ holds when $x_{k+1}$ is obtained;
    \item
    \label{item-lem:dist_rec_local}
    the following hold:
    \begin{align}
      \|x_{k+1} - x_k\| &\leq C_1 \dist(x_k),
      \label{eq:x_diff_dist_bound}\\
      \dist(x_{k+1}) &\leq C_2 \dist(x_k)^2.
      \label{eq: xxx implies dist quadratic}
    \end{align}
  \end{enuminlem}
\end{lemma}

\begin{proof}[Proof of \cref{item-lem:x_xk_upper_local}]
  From $x_k \in \mathcal B(\delta)$, $\delta \leq r/2$, and \cref{eq:a_in_r/2_implies_sol_in_r}, we have
  \begin{align}
    x_k \in \mathcal B(r/2)
    \quad\text{and}\quad
    \tilde x_k \in \mathcal B(r).
    \label{eq:local_proof_xk_in_ball}
  \end{align}
  Moreover,
  we have from $\nabla m^k_{\lambda}(x) = J_k^\top(F_k + J_k(x - x_k)) + \lambda(x - x_k)$ that
  \[
  \underbrace{\inner{\nabla m^k_{\lambda}(x)}{x - \tilde x_k}}_{\text{(A)}}
  =
  \underbrace{\inner{F_k + J_k(x - x_k)}{J_k(x - \tilde x_k)}}_{\text{(B)}}
  + \lambda \underbrace{\inner{x - x_k}{x - \tilde x_k}}_{\text{(C)}}.
  \]
  We bound the terms (A)--(C) as follows:
  \begin{align}
    \text{(A)} \leq
    \cstSIL \lambda \|F_k\| \|x - \tilde x_k\|
    &\leq
    \cstSIL \sigma \lambda\|x_k - \tilde x_k\| \|x - \tilde x_k\|\\
    &\leq
    \frac{\cstSIL^2 \sigma^2 \lambda}{2} \|x_k - \tilde x_k\|^2 + \frac{\lambda}{2} \|x - \tilde x_k\|^2,
  \end{align}
  where the first and second inequalities follow from 
  \cref{eq:asm_local_subproblem_y_mu} and \cref{eq:F_<=_dist}, respectively,
  and the last inequality follows from the arithmetic and geometric means;
  \[
    \text{(B)}
    \geq - \frac{1}{4} \|F_k + J_k(\tilde x_k - x_k)\|^2
    \geq - \frac{L^2}{16}\|\tilde x_k - x_k\|^4,
  \]
  where the first inequality follows from $4\inner{a}{b} = \|a+b\|^2-\|a-b\|^2 \geq -\|a-b\|^2$
  and the second inequality from \cref{item-lem: Lip J implies xxx,eq:local_proof_xk_in_ball,item-asm: J Lip on N(r)}; and
  \begin{equation}
    \text{(C)}
    = \frac{1}{2} \Big( 
    \|x - x_k\|^2 + \|x - \tilde x_k\|^2 - \|\tilde x_k - x_k\|^2
    \Big).
  \end{equation}
  Combining these bounds and rearranging terms yield
  \begin{equation}
    \|x - x_k\|^2
    \leq
    (1 + \cstSIL^2 \sigma^2) \|\tilde x_k - x_k\|^2 + \frac{L^2}{8\lambda}\|\tilde x_k - x_k\|^4.
    \label{eq:bound_yk_xk}
  \end{equation}
  From \cref{eq:local_proof_xk_in_ball}, \cref{item-asm: EB on N(r)}, and $\lambda = M \|F_k\| \geq M_0 \|F_k\|$, we have
  \begin{equation}
    \|\tilde x_k - x_k\|^2
    \leq
    \frac{r}{2} \times \frac{\|F_k\|}{\rho}
    \leq
    \frac{r \lambda}{2 \cstEB M_0}.
  \end{equation}
  Applying this bound to the second term on the right-hand side of \cref{eq:bound_yk_xk},
  we obtain the desired result \eqref{eq:y-x_dist_bound}.
\end{proof}

\begin{proof}[Proof of \cref{item-lem:M_upper_local}]
  As in \cref{item-lem:x_xk_upper_local}, let $x$ denote the $x$ obtained in \cref{line-LMGD:set_x_xkt}.
  By \cref{eq:def_delta}, \cref{eq:y-x_dist_bound}, and $x_k \in \mathcal B(\delta)$, we have 
  \begin{align}
    \|x - x^*\|
    &\leq
    \|x_k - x^*\| + \|x - x_k\|\\
    &\leq
    \|x_k - x^*\| + C_1 \dist(x_k)\\
    &\leq
    (1 + C_1) \|x_k - x^*\|
    \leq
    (1 + C_1) \delta
    = r/2,
  \end{align}
  i.e.,
  \begin{align}
    x \in \mathcal B(r/2).
    \label{eq:x_in_r/2}
  \end{align}
  We now have $x_k, x \in \mathcal C \cap \mathcal B(r)$.
  As in the proof of \cref{item-lem: hat L upper bound global}, by using \cref{lem: majorization LM} with $\mathcal X \coloneqq \mathcal C \cap \mathcal B(r)$, we see that if $M \geq L$ holds at \cref{line-LMGD:set_lam}, the outer iteration must be successful.
  This leads to the desired result.
\end{proof}

\begin{proof}[Proof of \cref{item-lem:dist_rec_local}]
  Eq.~\cref{eq:x_diff_dist_bound} follows from \cref{item-lem:x_xk_upper_local,item-lem:M_upper_local}.
  We prove \cref{eq: xxx implies dist quadratic} below.
  From \cref{eq:a_in_r/2_implies_sol_in_r,eq:x_in_r/2}, we have $x_{k+1}, \tilde x_{k+1} \in \mathcal B(r)$. Moreover, we have 
  \begin{align}
    &\mathInd
    \|F_{k+1}\|^2 - 
    \overbrace{\inner{\nabla m^k_{\lambda_k}(x_{k+1})}{x_{k+1} - \tilde x_{k+1}}}^{\text{(D)}}\\
    &=
    \inner{F_{k+1}}{F_{k+1} + J_{k+1}(\tilde x_{k+1} - x_{k+1})}
    + \inner{J_{k+1}^\top F_{k+1} - \nabla m^k_{\lambda_k}(x_{k+1})}{x_{k+1} - \tilde x_{k+1}}
    \\
    &\leq
    \underbrace{\|F_{k+1}\| \|F_{k+1} + J_{k+1}(\tilde x_{k+1} - x_{k+1})\|}_{\text{(E)}}
    +
    \underbrace{\|J_{k+1}^\top F_{k+1} - \nabla m^k_{\lambda_k}(x_{k+1})\|}_{\text{(F)}}
    \dist(x_{k+1})
  \end{align}
  and bound the terms (D)--(F) as follows:
  \begin{equation}
    \text{(D)}
    \leq
    \cstSIL \lambda_k \|F_k\| \dist (x_{k+1})
    \leq
    \cstSIL \bar M \|F_k\|^2 \dist (x_{k+1})
  \end{equation}
  by \cref{eq:asm_local_subproblem_x_lam,item-lem:M_upper_local};
  \begin{equation}
    \text{(E)}
    \leq
    \frac{L}{2} \|F_{k+1}\| \dist(x_{k+1})^2
    \leq
    \frac{L}{2\cstEB} \|F_{k+1}\|^2 \dist(x_{k+1})
    \leq
    \frac{L}{2\cstEB} \|F_k\|^2 \dist(x_{k+1})
  \end{equation}
  by \cref{item-lem: Lip J implies xxx}, \cref{item-asm: EB on N(r)}, and $\|F_{k+1}\| \leq \|F_k\|$ from \cref{eq: f(x_(k+1)) <= m(x_(k+1)) <= m(x_k) = f(x_k)}; and
  \begin{alignat}{2}
    \text{(F)}
    &= \| J_{k+1}^\top F_{k+1} - J_k^\top (F_k + J_ku) - \lambda_k u\|
    &\quad&\by{letting $u \coloneqq x_{k+1} - x_k$}\\
    &\leq \|J_k^\top (F_{k+1} - F_k - J_ku)\|\\
    &\qquad + \|(J_{k+1} - J_k)^\top F_{k+1}\| + \lambda_k \|u\|\\
    &\leq \frac{L\sigma}{2} \|u\|^2 + L\|F_{k+1}\|\|u\| + \lambda_k \|u\|
    &\quad&\by{\cref{eq:J_bounded_local}, \cref{item-lem: Lip J implies xxx},\\
    and \cref{item-asm: J Lip on N(r)}}\\
    &\leq \frac{L\sigma}{2} \|u\|^2 + (L + \bar M) \|F_k\| \|u\|
    &\quad&\by{$\|F_{k+1}\| \leq \|F_k\|$ and \cref{item-lem:M_upper_local}}\\
    &\leq \Big( \frac{L\sigma C_1^2}{2} + (L + \bar M)\sigma C_1 \Big) \dist(x_k)^2
    &\quad&\by[.]{\cref{eq:F_<=_dist,eq:x_diff_dist_bound}}
  \end{alignat}
  Combining these bounds yields
  \[
    \|F_{k+1}\|^2
    \leq
    \bigg(
    \Big( \cstSIL \bar M + \frac{L}{2\cstEB} \Big) \|F_k\|^2
    + \Big( \frac{L\sigma C_1^2}{2} + (L + \bar M)\sigma C_1 \Big) \dist(x_k)^2
    \bigg) \dist(x_{k+1}).
  \]
  We bound $\|F_k\|$ and $\|F_{k+1}\|$ in the above inequality by using \cref{item-asm: EB on N(r),eq:F_<=_dist} and obtain
  \[
    \cstEB^2 \dist(x_{k+1})^2
    \leq
    \bigg(
    \sigma^2 \Big( \cstSIL \bar M + \frac{L}{2\cstEB} \Big)
    + \frac{L\sigma C_1^2}{2} + (L + \bar M)\sigma C_1
    \bigg) \dist(x_k)^2\dist(x_{k+1}),
  \]
  which implies the desired result \eqref{eq: xxx implies dist quadratic}.
\end{proof}

\subsection{Local quadratic convergence}
Let us state the local quadratic convergence result of \cref{alg: proposed LM-GD}.
\begin{theorem}
  \label{thm: local superliner convergence and unsuccessful}
  Suppose that \cref{asm: local convergence} holds, and define $\bar M$ by \cref{eq: def of bar L global}. Set $x_0 \in \mathcal B(\delta_0)$ 
  for a sufficiently small constant $\delta_0 > 0$ such that
  \begin{equation}
    C_2 \delta_0 < 1,\quad
    \delta_0 + \frac{C_1 \delta_0}{1 - C_2 \delta_0} \leq \delta,
    \label{eq:def_delta0}
  \end{equation}
  where $C_1$, $C_2$, and $\delta$ are the constants defined in \cref{eq:def_C1,eq:def_C2,eq:def_delta}. Then,
  \begin{enuminthm}
  \item
  \label{item-thm: L_hat bound local}
  the number of unsuccessful iterations is at most $\ceil{\log_\alpha (\bar M / M_0)} = O(1)$, and
  \item
  \label{item-thm:quadratic_convergence}
  the sequence $(x_k)$ converges quadratically to an optimal solution $\hat x \in \mathcal X^*$.
  \end{enuminthm}
\end{theorem}

\begin{proof}[Proof of \cref{item-thm: L_hat bound local}]
  First, we will prove that
  \begin{subequations}
    \begin{align}
      &\text{$x_k \in \mathcal B(\delta)$, and}
      \label{eq:local_proof_induction1}\\
      &\text{$M \leq \bar M$ holds when $x_k$ is obtained}
      \label{eq:local_proof_induction2}
    \end{align}
  \end{subequations}
  for all $k \geq 0$ by induction.
  For $k=0$, \cref{eq:local_proof_induction1,eq:local_proof_induction2} are obvious.
  For a fixed $K \geq 0$, assume \cref{eq:local_proof_induction1,eq:local_proof_induction2} for all $k \leq K$.
  We then have \cref{eq:x_diff_dist_bound,eq: xxx implies dist quadratic,eq:local_proof_induction2} for $k \leq K+1$ by \cref{lem:local_induction_step}.
  To complete the induction, we prove \cref{eq:local_proof_induction1} for $k = K+1$.
  Solving the recursion of \cref{eq: xxx implies dist quadratic} and using $\dist(x_0) \leq \delta_0$, we have
  \begin{equation}
    \dist(x_k)
    \leq \dist(x_0) (C_2 \dist(x_0))^{2^k - 1}
    \leq \delta_0 (C_2 \delta_0)^{2^k - 1}
    \leq \delta_0 (C_2 \delta_0)^k
    \label{eq:dist_linear_convergence}
  \end{equation}
  for all $k \leq K+1$.
  We obtain \cref{eq:local_proof_induction1} for $k = K+1$ as follows:
  \begin{alignat}{2}
    \|x_{K+1} - x^*\|
    &\leq \|x_0 - x^*\| + \sum_{k=0}^K \|x_{k+1} - x_k\|
    &\quad&\by{the triangle inequality}\\
    &\leq \delta_0 + C_1 \sum_{k=0}^K \dist(x_k)
    &\quad&\by{\cref{eq:x_diff_dist_bound}}\\
    &\leq \delta_0 + \frac{C_1 \delta_0}{1 - C_2 \delta_0}
    \leq \delta
    &\quad&\by{\cref{eq:dist_linear_convergence,eq:def_delta0}}.
  \end{alignat}
  Now, we have proved \cref{eq:local_proof_induction1,eq:local_proof_induction2} for all $k \geq 0$.
\end{proof}

\begin{proof}[Proof of \cref{item-thm:quadratic_convergence}]
  Note that we have proved \cref{eq: xxx implies dist quadratic,eq:dist_linear_convergence} for all $k \geq 0$ in the proof of \cref{item-thm: L_hat bound local}.
  By \cref{eq:dist_linear_convergence} and $C_2\delta_0<1$ in \cref{eq:def_delta0}, we have
  \begin{equation}
    \label{eq:dist x_k to 0}
    \lim_{k \to \infty}\dist(x_k) = 0.
  \end{equation}
  As with \cref{eq:dist_linear_convergence}, we have for $i \geq k$,
  \begin{equation}
    \dist(x_i) \leq \dist(x_k) (C_2 \dist(x_k))^{2^{i-k} - 1}
    \leq \dist(x_k) (C_2 \delta_0)^{i-k}.
    \label{eq:dist_linear_convergence_offset}
  \end{equation}
  Using this bound and \cref{eq:y-x_dist_bound}, we obtain
  \begin{equation}
    \|x_k - x_l\|
    \leq \sum_{i=k}^{l-1} \|x_{i+1} - x_i\|
    \leq C_1 \sum_{i=k}^{l-1} \dist(x_i)
    \leq \frac{C_1}{1 - C_2\delta_0} \dist(x_k)
    \label{eq:Cauchy_seq_xk}
  \end{equation}
  for all $k, l$ such that $0 \leq k < l$.
  \Cref{eq:Cauchy_seq_xk,eq:dist x_k to 0} imply that $(x_k)$ is a Cauchy sequence.
  Accordingly, the sequence $(x_k)$ converges to a point $\hat x \in \mathcal X^*$ by \cref{eq:dist x_k to 0}.
  Thus, we obtain
  \begin{alignat}{2}
    \|x_{k+1} - \hat x\|
    &=
    \lim_{l \to \infty} \|x_{k+1} - x_l\|
    &\quad&\text{(by the continuity of $\|\cdot\|$)}\\
    &\leq
    \frac{C_1}{1 - C_2\delta_0} \dist(x_{k+1})
    &\quad&\text{(by \cref{eq:Cauchy_seq_xk})}\\
    &\leq
    \frac{C_1C_2}{1 - C_2\delta_0} \dist(x_k)^2
    &\quad&\text{(by \cref{eq: xxx implies dist quadratic})}\\
    &\leq
    \frac{C_1C_2}{1 - C_2\delta_0} \|x_k - \hat x\|^2
    &\quad&\text{(by $\hat x \in \mathcal X^*$)},
  \end{alignat}
  which implies \cref{item-thm:quadratic_convergence}.
\end{proof}


\section{Practical variant of the proposed method}
\label{sec:generalization}
We present a more practical variant (\cref{alg:proposed_experiment}) of \cref{alg: proposed LM-GD}, which also achieves the theoretical guarantees given for \cref{alg: proposed LM-GD} in \cref{sec: global convergence analysis,sec: local convergence analysis}.

\subsection{Generalized version of \cref{alg: proposed LM-GD}}
To obtain the practical variant, we first present a generalized framework of \cref{alg: proposed LM-GD}.
\cref{alg: proposed LM-GD} runs the vanilla projected gradient (PG) method in the inner loop. 
This PG can be replaced with other algorithms keeping $O(\epsilon^{-2})$ iteration complexity and quadratic convergence that were gained for \cref{alg: proposed LM-GD}.
Indeed, these theoretical results rely on the fact that the $x$ obtained in \cref{line-LMGD:set_x_xkt} of \cref{alg: proposed LM-GD} satisfies the following conditions:
\begin{condition}[for $O(\epsilon^{-2})$ iteration complexity bound]
  \label{cond:iteration_complexity}
  There exists a constant $\gamma > 0$ such that for all $k$,
  \begin{equation}
    m^k_{\lambda} (x) - m^k_{\lambda} (x_k)
    \leq
    - \mathcal D_\gamma(x_k).
    \label{eq:sufficient_decrease_PG_general}
  \end{equation}
\end{condition}
\begin{condition}[for local quadratic convergence]
  \label{cond:local_convergence}
  Both of the following hold:
  \begin{enumincond}
    \item
    $m^k_\lambda(x) \leq m^k_\lambda(x_k)$ for all $k$;
    \item 
    there exists a constant $c > 0$ such that $x$ is a $(\cstSIL \lambda \|F_k\|)$-stationary point of subproblem~\eqref{eq:subproblem} for all $k$.
  \end{enumincond}
\end{condition}
This fact yields a general algorithmic framework that achieves the $O(\epsilon^{-2})$ iteration complexity bound together with the quadratic convergence as in \cref{alg: proposed LM}.

\begin{algorithm}[t]
  \small
  \caption{
    Generalized version of \cref{alg: proposed LM-GD}
  }
  \label{alg: proposed LM}
  \begin{algorithmic}[1]
    \Require{
      $x_0 \in \mathcal C$, $M_0 > 0$, $\alpha > 1$
    }
    \State{$M \gets M_0$, $k \gets 0$}
    \Comment{initialization}
    \label{line-alg: initialization}
    \Repeat
      \State{$\lambda \gets M\|F_k\|$}
      \label{line-alg: def of lam}
      \State{%
        \pbox[t]{\linewidth}{%
          Compute an approximate solution $x \in \mathcal C$ to subproblem~\cref{eq:subproblem} that satisfies\\ \cref{cond:iteration_complexity} or \ref{cond:local_convergence} or both.
        }
      }
      \label{line-alg: subproblem}
      \If{$f(x) \leq m^k_\lambda(x)$}
        \State{$(x_{k+1}, \lambda_k) \coloneqq (x, \lambda)$, $k \gets k+1$}
        \Comment{successful}
        \label{line-alg: successful}
      \Else
        \State{$M \gets \alpha M$}
        \Comment{unsuccessful}
        \label{line-alg: unsuccessful}
      \EndIf
    \Until{a solution with a desired accuracy is obtained}
  \end{algorithmic}
\end{algorithm}

\begin{algorithm}[t]
  \small
  \caption{Proposed LM method using APG with adaptive restart}
  \label{alg:proposed_experiment}
  \begin{algorithmic}[1]
    \Require{
      \pbox[t]{40em}{
        $x_0 \in \mathcal C$;
        $M_0, \eta_0 > 0$;
        $\alpha, \alpha_{\mathrm{in}} > 1$;
        $0 < \beta, \beta_{\mathrm{in}} < 1$;
        $0 < M_{\mathrm{min}} \leq M_0$;
        $T \in \Z_{> 0} \cup \set{\infty}$;
        $c > 0$
      }
      \vspace{0.2\baselineskip}
    }
    \State{$M \gets M_0$, $\eta \gets \eta_0$, $k \gets 0$}
    \Repeat
    \Comment{outer loop}
      \State{$\lambda \gets M \|F_k\|$}
      \State{$\eta \gets \max \set{\eta, \lambda}$}
      \State{%
        $x_{k, -1} \gets x_k$,
        $x_{k, 0} \gets x_k$,
        $\theta_{-1} \gets 1$,
        $t \gets 0$
      }
      \Repeat
      \Comment{inner loop (APG)}
        \State{$\theta_t \gets \sqrt{\lambda / \eta}$}
        \State{$y \gets x_{k, t} + \frac{\theta_t (1 - \theta_{t-1})}{\theta_{t-1} (1 + \theta_t)} (x_{k,t} - x_{k,t-1})$}
        \State{$z \gets \proj_{\mathcal C}(y - \frac{1}{\eta} \nabla m^k_\lambda(y))$}
        \If{$m^k_\lambda(z) \leq m^k_\lambda(y) + \inner{\nabla m^k_\lambda(y)}{z - y} + \frac{\eta}{2}\|z - y\|^2$}
          \If{$m^k_\lambda(z) \leq m^k_\lambda(x_{k, t})$}
            \State{$x_{k, t+1} \gets z$, $t \gets t+1$}
            \State{$\eta \gets \max \set{\beta_{\mathrm{in}} \eta, \lambda}$}
          \Else 
            \State{%
              $x_{k, t-1} \gets x_{k, t}$,
              $\theta_{t-1} \gets 1$
            }
            \Comment{restart of APG}
          \EndIf
        \Else
          \State{$\eta \gets \alpha_{\mathrm{in}} \eta$}
        \EndIf
      \Until{($t = T$) \textbf{or} ($t \geq 1$ \textbf{and} $x_{k,t}$ is a $(\cstSIL \lambda \|F_k\|)$-stationary point of subproblem~\eqref{eq:subproblem})}
      \label{line-alg:check_xkt}
      \State{$x \gets x_{k, t}$}
      \label{line-LMAPG:x_gets_xkt}
      \If{$f(x) \leq m^k_\lambda(x)$}
        \State{$(x_{k+1}, \lambda_k) \coloneqq (x, \lambda)$, $k \gets k+1$}
        \Comment{successful}
        \State{$M \gets \max \set{\beta M, M_{\min}}$}
      \Else
        \State{$M \gets \alpha M$}
        \Comment{unsuccessful}
      \EndIf
    \Until{a solution with a desired accuracy is obtained}
  \end{algorithmic}
\end{algorithm}

In \cref{line-alg: subproblem} of \cref{alg: proposed LM}, any globally convergent algorithm for subproblem~\eqref{eq:subproblem} can be employed.
For example, we may use (block) coordinate descent methods, Frank–Wolfe methods, interior point methods, active set methods, or augmented Lagrangian methods.
For unconstrained cases, since the subproblem reduces to solving a system of linear equations, we may use conjugate gradient methods or direct methods, including Gaussian elimination.

\subsection{Proposed method with an accelerated projected gradient}
A practical example of \cref{alg: proposed LM} is presented in \cref{alg:proposed_experiment}.
This algorithm employs the \emph{accelerated} projected gradient (APG) method~\citep[Algorithm~1]{lin2014adaptive} with the adaptive restarting technique~\citep[Section 3.2]{odonoghue2015adaptive} to solve subproblems and adopts the additional parameters mentioned in \cref{rem: additional input parameter}.
Since the solution $x$ obtained in \cref{line-LMAPG:x_gets_xkt} of \cref{alg:proposed_experiment} satisfies \cref{cond:iteration_complexity}, this algorithm enjoys the $O(\epsilon^{-2})$ iteration complexity bound.
In addition, it also achieves the $O(\epsilon^{-2})$ overall complexity bound if $T < \infty$ as with \cref{cor: global overall complexity}, and it achieves local quadratic convergence if $T = \infty$.
\cref{alg:proposed_experiment} will be used for the numerical experiments in the next section.

\section{Numerical experiments}
\label{sec:experiments}
We examine the practical performance of the proposed method.
We implemented all methods in Python with SciPy~\citep{virtanen2020scipy} and JAX~\citep{jax2018github} and executed them on a computer with Apple M1 Chip (8 cores, 3.2 GHz) and 16 GB RAM.

\subsection{Problem setting}
We consider three types of instances: (i) compressed sensing with quadratic measurement, (ii) nonnegative matrix factorization with missing values, and (iii) autoencoder with MNIST dataset.
\subsubsection{Compressed sensing with quadratic measurement}
Given $A_1,\dots,A_n \in \R^{r \times d}$, $b_1,\dots,b_n \in \R^d$, and $c_1,\dots,c_n, R \in \R$, we consider the following problem:
\begin{align}
  \min_{x \in \R^d}\ 
  \sum_{i=1}^n \prn[\Big]{
    \frac{1}{2r} \norm{A_i x}^2 + \inner{b_i}{x} - c_i
  }^2
  \quad\subjectto\quad
  \norm{x}_1 \leq R,
  \label{eq:problem_cs}
\end{align}
where $\norm{\cdot}_1$ denotes $\ell_1$-norm.
Problem \cref{eq:problem_cs} formulates the situation where a sparse vector $x^* \in \R^d$ is recovered from a small number (i.e., $n < d$) of quadratic observations, $\frac{1}{2r} \norm{A_i x^*}^2 + \inner{b_i}{x^*}$ for $i = 1, \dots, n$.
Such a problem arises in the context of compressed sensing~\citep{li2013sparse,blumensath2013compressed} and phase retrieval~\citep{candes2015phase,zhang2017nonconvex}.
Problem \cref{eq:problem_cs} can be transformed into the form of problem~\cref{eq:mainproblem}.

\myparagraph{Generating instances}
First, we generate the optimal solution $x^* \in \R^d$ with only $d_{\mathrm{nnz}} \,(< d)$ nonzero entries.
The indexes of the nonzero entries are chosen uniformly randomly, and the value of those elements are independently drawn from the uniform distribution on $[-x_{\mathrm{max}}, x_{\mathrm{max}}]$.
Each entry of $A_i$'s and $b_i$'s is drawn independently from the standard normal distribution $\mathcal N(0, 1)$.
Then, we set $R = \norm{x^*}_1$ and $c_i = \frac{1}{2r} \norm{A_i x^*}^2 + \inner{b_i}{x^*}$ for all $i$.
We fix $d = 200$, $r = 10$, and $n = 50$, and set $d_{\mathrm{nnz}} \in \set{5, 10, 20}$ and $x_{\mathrm{max}} \in \set{0.1, 1}$.
We set the starting point for each algorithm as $x_0 = \0$.

\subsubsection{Nonnegative matrix factorization with missing values}
Given $A \in \R^{m \times n}$ and $H \in \set{0, 1}^{m \times n}$, we consider the following problem:
\begin{align}
  \min_{X \in \R^{m \times r},\ Y \in \R^{n \times r}}
  \norm{ H \odot (X Y^\top - A) }_{\mathrm{F}}^2
  \quad\subjectto\quad
  X \geq O,\quad
  Y \geq O,
  \label{eq:problem_nmf}
\end{align}
where $\odot$ denotes the elementwise product, $X \geq O$ and $Y \geq O$ denote elementwise inequalities, and $\norm{\cdot}_{\mathrm{F}}$ denotes the Frobenius norm.
Problem~\cref{eq:problem_nmf} formulates the situation where a data matrix $A$ with some missing entries is approximated by the product $X Y^\top$ of two nonnegative matrices.
Such a problem is called nonnegative matrix factorization (NMF) with missing values and is widely used for nonnegative data analysis, especially for collaborative filtering~\citep{zhang2006learning,luo2014efficient}.
For more information on NMF, see~\citep{wang2012nonnegative,berry2007algorithms} and the references therein.
Problem~\cref{eq:problem_nmf} can also be written as problem~\cref{eq:mainproblem}.


\myparagraph{Generating instances}
To generate $A$ and $H$, we introduce two parameters: $\gamma \geq 1$ and $0 < p \leq 1$. 
The parameters $\gamma$ and $p$ control the condition number of $A$ and the number of $1$'s in $H$, respectively.
Let $l \coloneqq \min \set{m, n}$.
First, a matrix $\tilde A \in \R^{m \times n}$ is generated by $\tilde A = U D V^\top$, and then the matrix $A$ is obtained by normalizing $\tilde A = (\tilde a_{ij})_{i,j}$ as $A = \tilde A / \max_{i, j} \tilde a_{ij}$.
Here, each entry of $U \in \R^{m \times l}$ and $V \in \R^{n \times l}$ follows independently the uniform distribution on $[0, 1]$, and $D = \diag(\gamma^{0}, \gamma^{-1/l}, \gamma^{-2/l},\dots, \gamma^{-(l-1)/l} ) \in \R^{l \times l}$ is a diagonal matrix.
$H$ is a random matrix whose entries follow independently the Bernoulli distribution with parameter $p$, i.e., each entry of $H$ is $1$ with probability $p$.
We fix $m = n = 50$ and $\gamma = 10^5$, and set $r \in \set{10, 40}$ and $p \in \set{0.02, 0.1, 0.5}$.
Since $(X, Y) = (O, O)$ is a stationary point of problem~\cref{eq:problem_nmf}, we set the starting point to random matrices whose entries independently follow the uniform distribution on $[0, 10^{-3}]$.

\subsubsection{Autoencoder with MNIST dataset}
The third instance is highly nonlinear and large-scale.
In machine learning, autoencoders (see, e.g., \citep[Section 14]{goodfellow2016deep}) are a popular model to compress real-world data, represented as high-dimensional vectors, into low-dimensional vectors.
Given $p$-dimensional data $a_1, \dots, a_N \in \R^p$, autoencoders try to learn an encoder $\phi_x^{\mathrm{enc}}: \R^p \to \R^q$ and a decoder $\phi_y^{\mathrm{dec}}: \R^q \to \R^p$, where $q < p$.
Here, $x$ and $y$ are parameters to be learned by solving the following optimization problem:
\begin{align}
  \min_{x, y}\ 
  \sum_{i=1}^N \norm*{a_i - \phi_y^{\mathrm{dec}}(\phi_x^{\mathrm{enc}}(a_i))}^2.
  \label{eq:problem_autoencoder}
\end{align}
As we see from the optimization problem above, the autoencoder aims to extract latent features that can be used to reconstruct the original data.



For this experiment, we use the MNIST hand-written digit dataset.
Each data is a $28 \times 28$ pixel grayscale image, which is represented as a vector $a_i \in [0, 1]^p$ with $p = 28 \times 28 = 728$.
The dataset contains 60,000 training data, of which $N = 1000$ were randomly chosen for use.
We set $q = 16$; our model encodes $728$-dimensional data into $16$ dimensions.
Both encoder and decoder are two-layer neural networks with a hidden layer of size $64$ and logistic sigmoid activation functions.
Specifically, the encoder $\phi_x^{\mathrm{enc}}$ is written as
\begin{align}
  \phi_x^{\mathrm{enc}}(a)
  =
  \phi_x^2 \circ
  \phi_x^1 (a),
  \quad\text{where}\quad
  \phi_x^i(a)
  &\coloneqq
  S( W_i a + b_i ).
\end{align}
Here, $S$ is the elementwise logistic sigmoid function, and $W_1 \in \R^{728 \times 64}$, $b_1 \in \R^{64}$, $W_2 \in \R^{64 \times 16}$, and $b_2 \in \R^{16}$ are parameters of the network; $x = ((W_i, b_i))_{i=1}^2$.
The decoder $\phi_y^{\mathrm{dec}}$ is formulated in a similar way.
When we rewrite problem~\cref{eq:problem_autoencoder} in the form of \cref{eq:mainproblem}, the dimension of the function $F: \R^d \to \R^n$ is $d = 96{,}104$
and $n = Np = 728{,}000$.

\subsection{Algorithms and implementation}
We compare the proposed method with six existing methods.
The details are below.

\myparagraph{Proposed (\cref{alg:proposed_experiment}) and Proposed-NA (\cref{alg: proposed LM-GD}) method}
To see the effect of acceleration for subproblems, we implemented both \cref{alg: proposed LM-GD,alg:proposed_experiment}; \cref{alg:proposed_experiment} is expected to be faster, of course.
In \cref{line-LMGD:inner loop_until} of \cref{alg: proposed LM-GD} and \cref{line-alg:check_xkt} of \cref{alg:proposed_experiment}, we have to check if $x_{k,t}$ is a $(\cstSIL \lambda \|F_k\|)$-stationary point, but it is not very easy.
We thus replace the criterion with one using gradient mapping, i.e., check if $\eta \norm{x_{k,t} - y} \leq \cstSIL \lambda \|F_k\|$.
The input parameters of \cref{alg: proposed LM-GD,alg:proposed_experiment} are set to $M_0 = \eta_0 = 1$, $\alpha = \alpha_{\mathrm{in}} = 2$, $\beta = \beta_{\mathrm{in}} = 0.9$, $M_{\min} = 10^{-10}$, $T = 100$, and $c = 1$.

\myparagraph{Fan method \citep[Algorithm 2.1]{fan2013levenberg} and KYF method \citep[Algorithm 2.12]{kanzow2004levenberg}}
The Fan and KYF methods are constrained LM methods with a global convergence guarantee.
To solve subproblem~\eqref{eq:subproblem}, an APG method is used as well as \cref{alg:proposed_experiment} for a fair comparison.
The difference from the APG in \cref{alg:proposed_experiment} is in the stopping criterion; the condition ``$x_{k,t}$ is a $(\cstSIL \lambda \|F_k\|)$-stationary point'' in \cref{line-alg:check_xkt} of \cref{alg:proposed_experiment} is replaced with $\eta \norm{x_{k,t} - y} \leq 10^{-9}$.\footnote{
  Since $\cstSIL \lambda \|F_k\|$ in \cref{alg:proposed_experiment} derives from our update rule of $\lambda$ and our analysis, it does not seem appropriate to use a criterion with $\cstSIL \lambda \|F_k\|$ directly in another algorithm.
  We thus use $\eta \norm{x_{k,t} - y} \leq 10^{-9}$ instead of $\eta \norm{x_{k,t} - y} \leq \cstSIL \lambda \|F_k\|$ here.
}
The input parameters in \citep{kanzow2004levenberg,fan2013levenberg} are set to $\mu = 10^{-4}$, $\beta = 0.9$, $\sigma = 10^{-4}$, $\gamma= 0.99995$, and $\delta \in \set{1, 2}$, following the recommendations of \citep{kanzow2004levenberg,fan2013levenberg}.\footnote{
  $\delta = 1$ and $\delta = 2$ correspond to the Fan method and the KYF method, respectively.
}

\myparagraph{Facchinei method \citep[Algorithm 3]{facchinei2013family}}
This is a constrained LM method that allows subproblems to be solved inexactly.
We solve the subproblems in almost the same way as the Fan and KYF methods.
The input parameters in \citep{facchinei2013family} are set to $\gamma_0 = 1$ and $S = 2$.

\myparagraph{GGO method \citep[Algorithm G-LMA-IP]{goncalves2021inexact}}
This is an LM-type method that requires the solution of a linear system at each iteration.
The main advantage of this algorithm is that it does not require exact projection and can be applied to problems with a complex feasible region.
Still, it is reported to perform well even when the projection is easy to compute exactly~\citep{goncalves2021inexact}.
The linear systems are solved via QR decomposition (\texttt{scipy.linalg.qr}~\citep{virtanen2020scipy}) and the input parameters in \citep{goncalves2021inexact} are set to $M \in \set{1, 15}$, $\eta_1 = 10^{-4}$, $\eta_2 = 10^{-2}$, $\eta_3 = 10^{10}$, $\gamma = 10^{-3}$, $\beta = 1/2$, and $\theta_k = 0$, following~\citep{goncalves2021inexact}.\footnote{
  Because the algorithm with $M = 1$ outperformed $M = 15$ in our experiments, we omit the results with $M = 15$.
}

\myparagraph{Projected gradient (PG) method}
The PG method is one of the most standard first-order methods for problem~\cref{eq:mainproblem}.
The step-size is adaptively chosen in a similar way to the APG in \cref{alg:proposed_experiment} with $\eta_0 = 1$, $\alpha_{\mathrm{in}} = 2$, and $\beta_{\mathrm{in}} = 0.9$.

\myparagraph{Trust-region reflective (TRF) method}
This is an interior trust-region method for box-constrained nonlinear optimization.
It was proposed in \citep{branch1999subspace} and is implemented in SciPy~\citep{virtanen2020scipy} with several improvements.
For the TRF method, we call \texttt{scipy.optimize.least\_squares}~\citep{virtanen2020scipy} with a \texttt{gtol=1e-5} option to avoid the long execution time caused by searching for too precise a solution.

\myparagraph{Other information}
As mentioned in \cref{sec:oracle_model}, there are two ways to handle Jacobian matrices: explicitly computing $J_k \coloneqq J(x_k)$ or using Jacobian-vector products $J_k u$ and $J_k^\top v$.
In our experiments, the latter implementation outperformed the former, so we adopted the latter if possible (i.e., for Proposed, Proposed-NA, Fan, KYF, Facchinei, and PG).\footnote{
  When using the Jacobian-vector products, i.e., not computing the Jacobian explicitly, almost all of the algorithm's runtime is spent solving subproblems.
}
We note that GGO is based on QR decomposition, which is probably impossible to implement using Jacobian-vector products.

For projection onto the feasible region of problem \cref{eq:problem_cs}, we employ \citep[Algorithm 1]{condat2016fast}, whose time complexity is $O(d \log d)$.

\begin{figure}[thbp]
  \centering
  \includegraphics[width=0.9\linewidth]{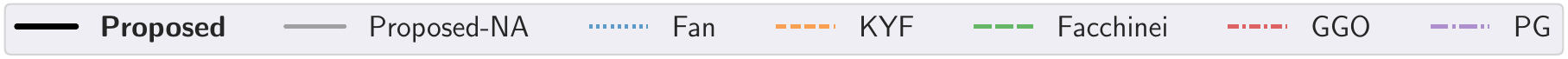}\par\medskip
  \subcaptionbox{$x_{\mathrm{max}} = 0.1$, $d_{\mathrm{nnz}} = 5$}[0.32\linewidth]{
    \includegraphics[width=\linewidth]{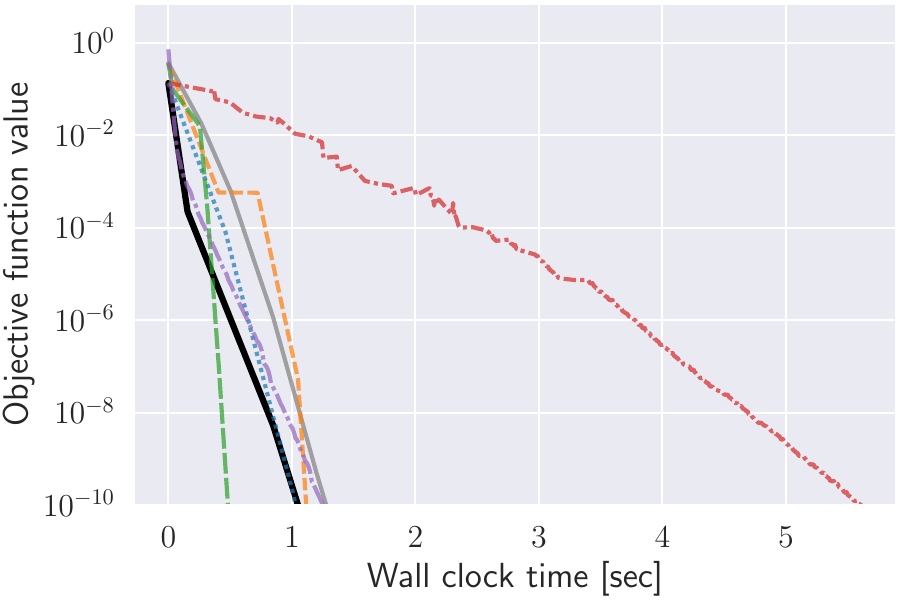}\par
    \includegraphics[width=\linewidth]{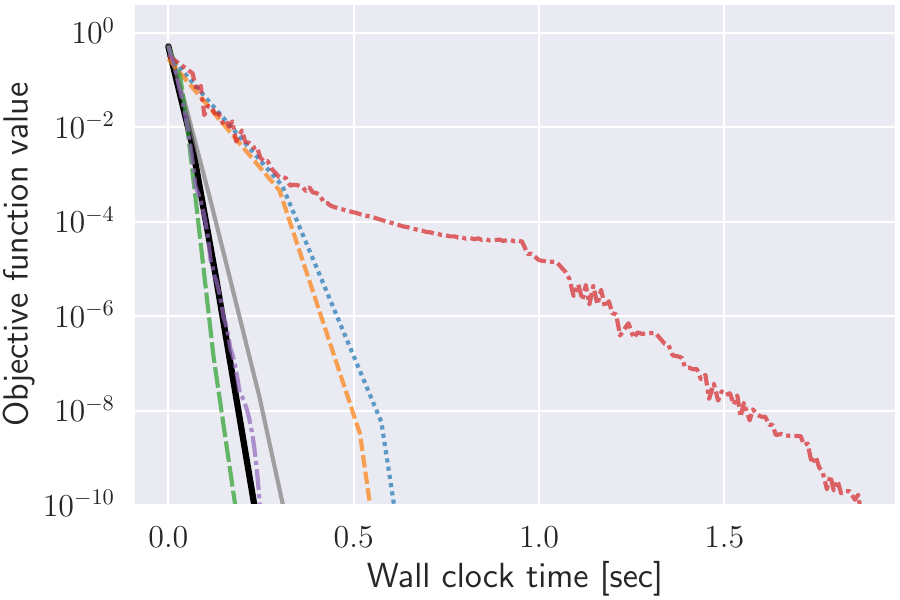}\par
  }\hfill
  \subcaptionbox{$x_{\mathrm{max}} = 0.1$, $d_{\mathrm{nnz}} = 10$\label{fig:exp_cs_01_10}}[0.32\linewidth]{
    \includegraphics[width=\linewidth]{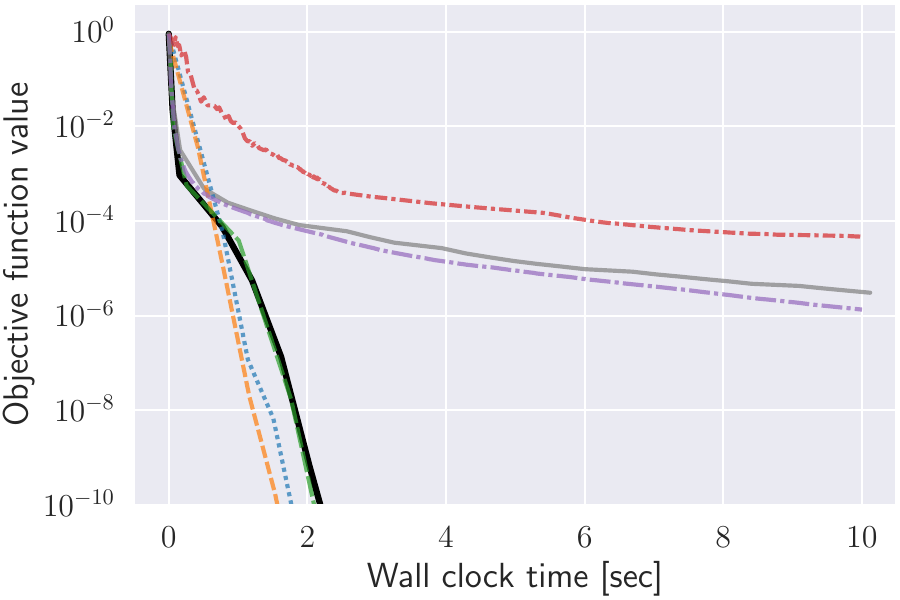}\par
    \includegraphics[width=\linewidth]{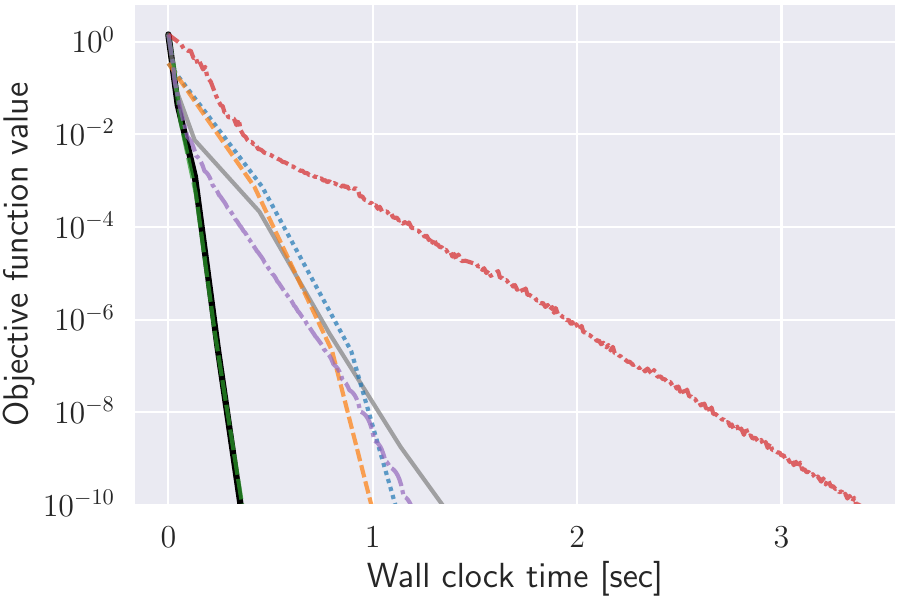}\par
  }\hfill
  \subcaptionbox{$x_{\mathrm{max}} = 0.1$, $d_{\mathrm{nnz}} = 20$}[0.32\linewidth]{
    \includegraphics[width=\linewidth]{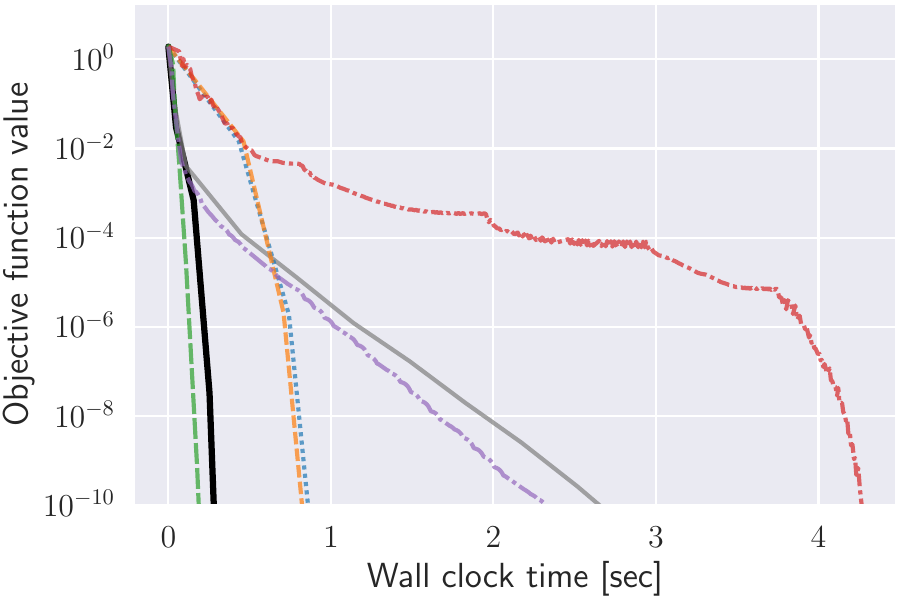}\par
    \includegraphics[width=\linewidth]{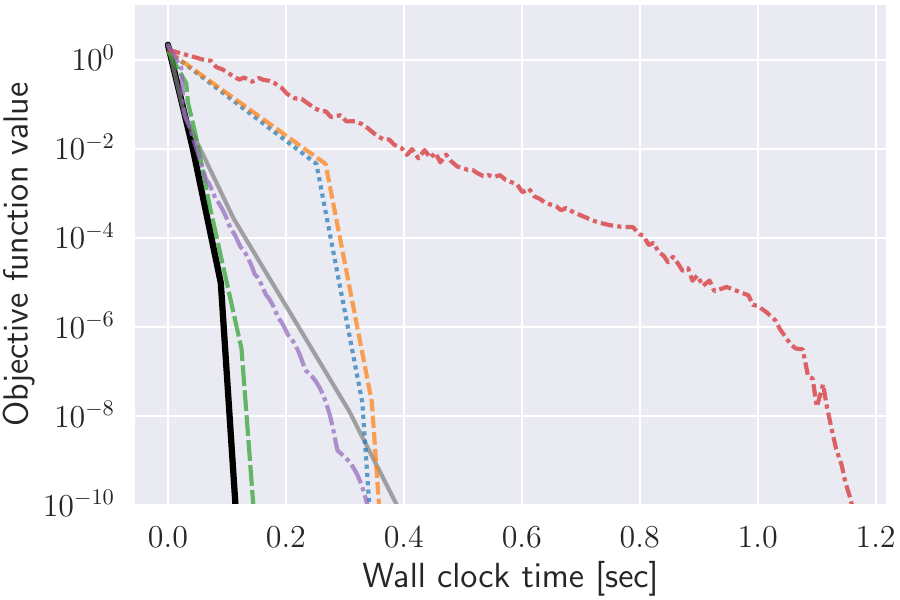}\par
  }\par\bigskip
  \subcaptionbox{$x_{\mathrm{max}} = 1$, $d_{\mathrm{nnz}} = 5$}[0.32\linewidth]{
    \includegraphics[width=\linewidth]{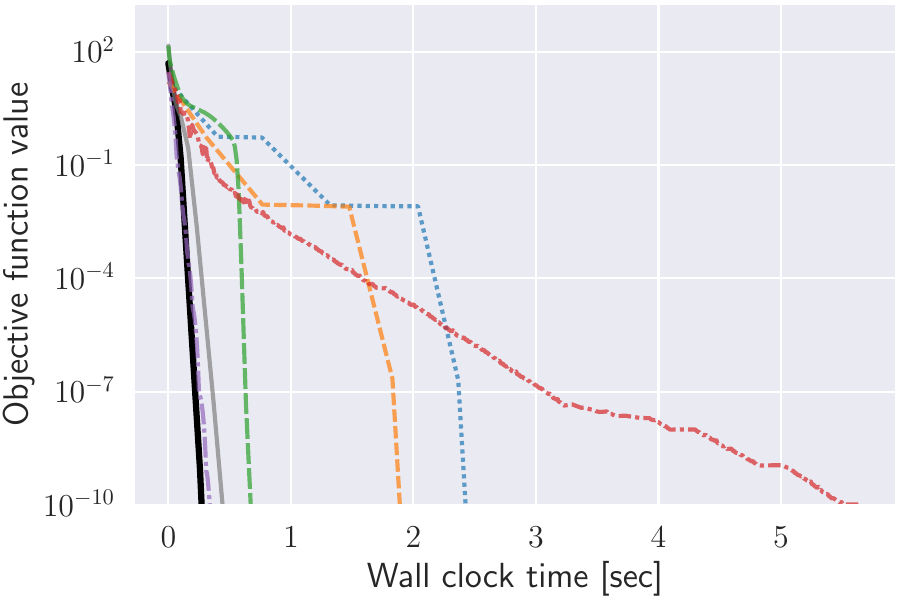}\par
    \includegraphics[width=\linewidth]{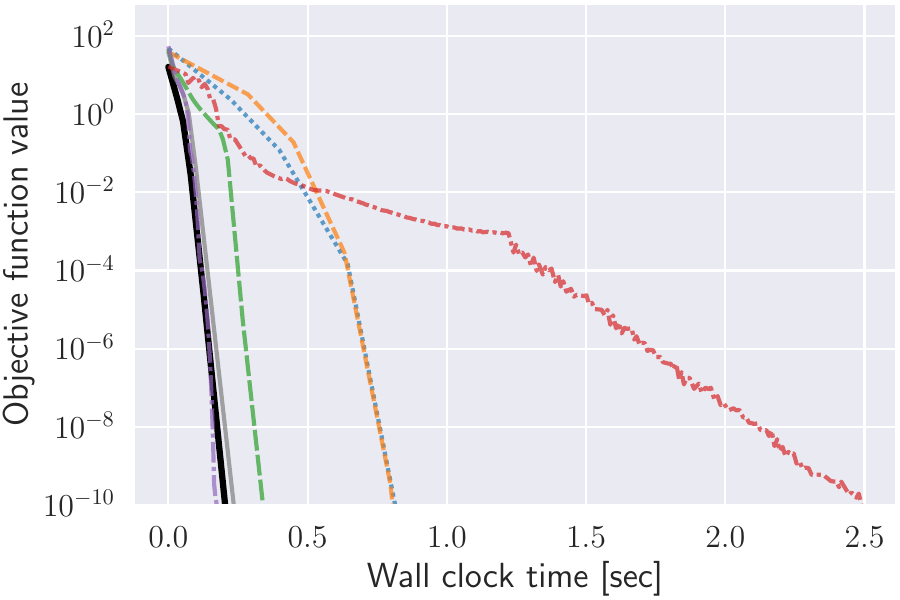}\par
  }\hfill
  \subcaptionbox{$x_{\mathrm{max}} = 1$, $d_{\mathrm{nnz}} = 10$}[0.32\linewidth]{
    \includegraphics[width=\linewidth]{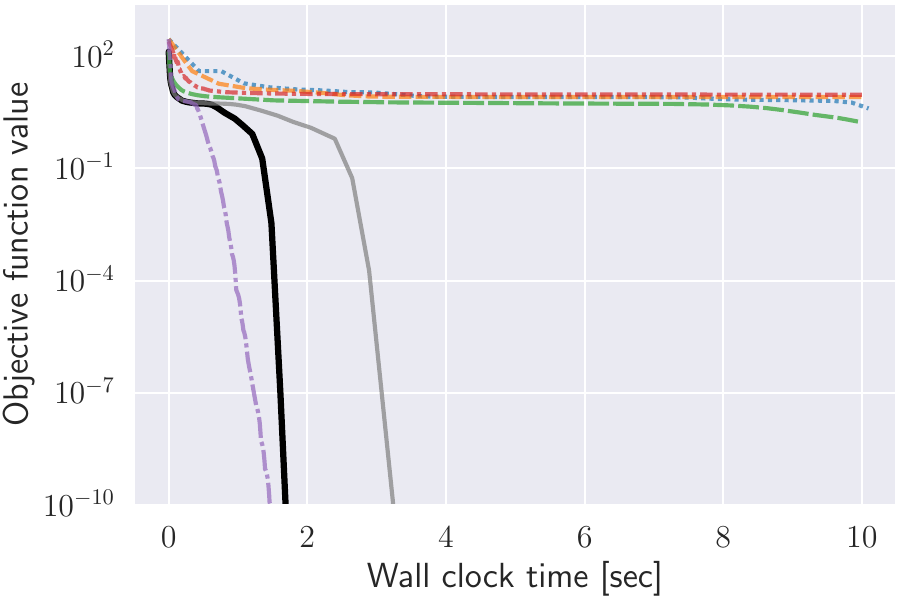}\par
    \includegraphics[width=\linewidth]{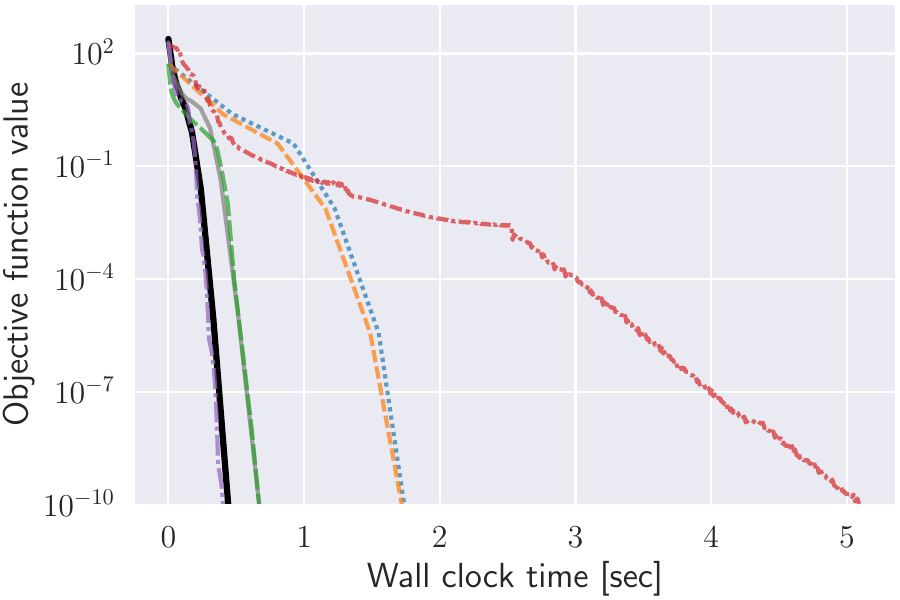}\par
  }\hfill
  \subcaptionbox{$x_{\mathrm{max}} = 1$, $d_{\mathrm{nnz}} = 20$\label{fig:exp_cs_1_20}}[0.32\linewidth]{
    \includegraphics[width=\linewidth]{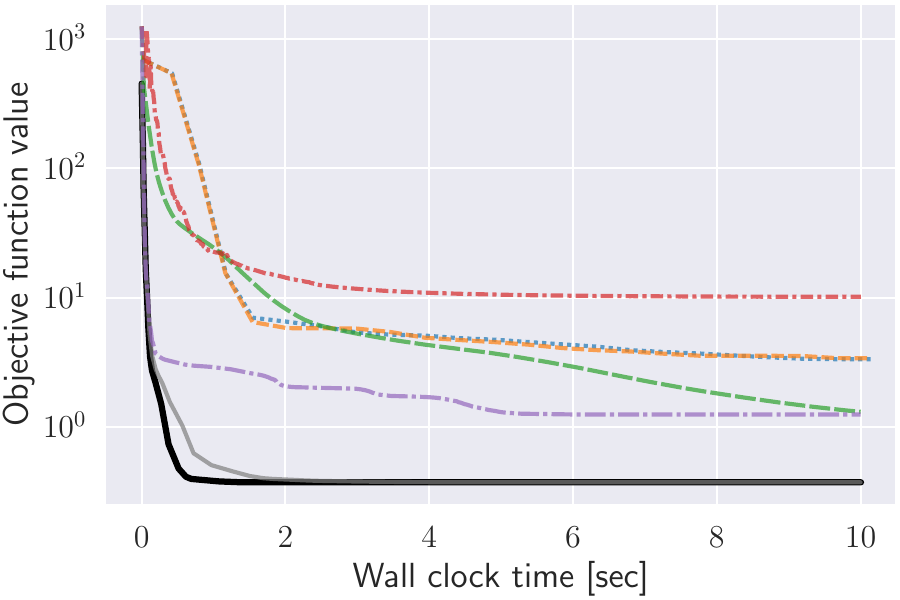}\par
    \includegraphics[width=\linewidth]{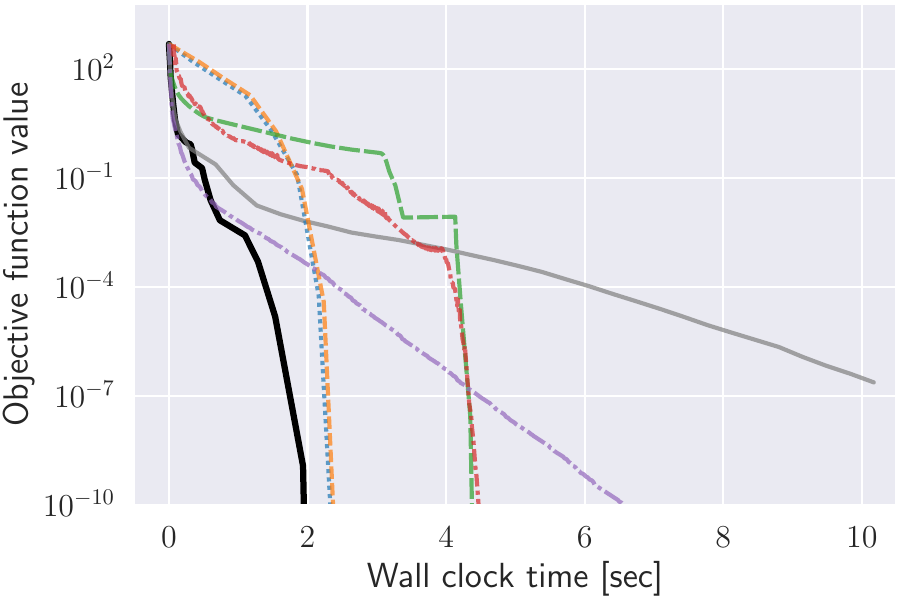}\par
  }\par\medskip
  \caption{Results of compressed sensing (problem~\cref{eq:problem_cs}).}
  \label{fig:exp_cs}
\end{figure}

\begin{figure}[thp]
  \centering
  \includegraphics[width=0.9\linewidth]{figure/legend.pdf}\par\medskip
  \subcaptionbox{$r = 10$, $p = 0.02$\label{fig:exp_nmf_10_002}}[0.32\linewidth]{
    \includegraphics[width=\linewidth]{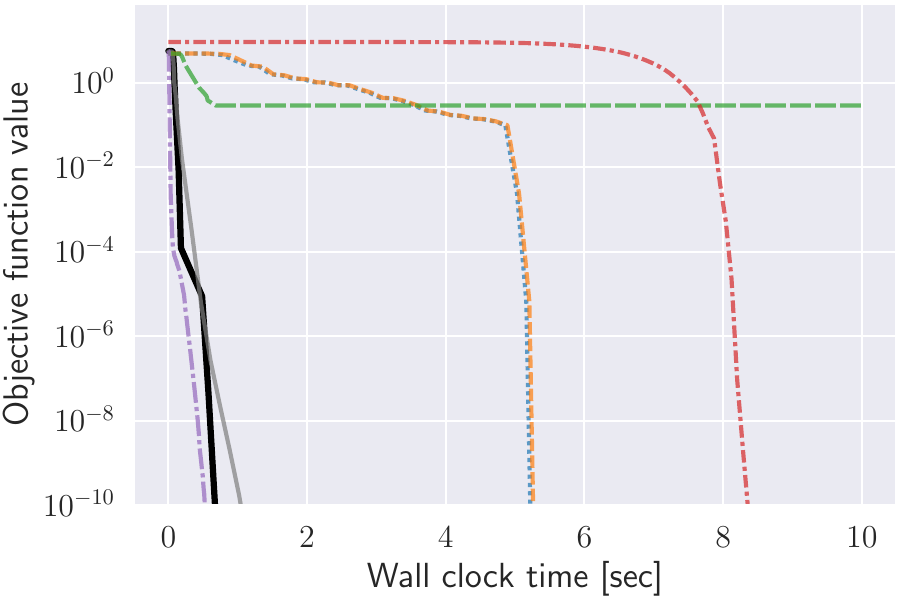}\par
    \includegraphics[width=\linewidth]{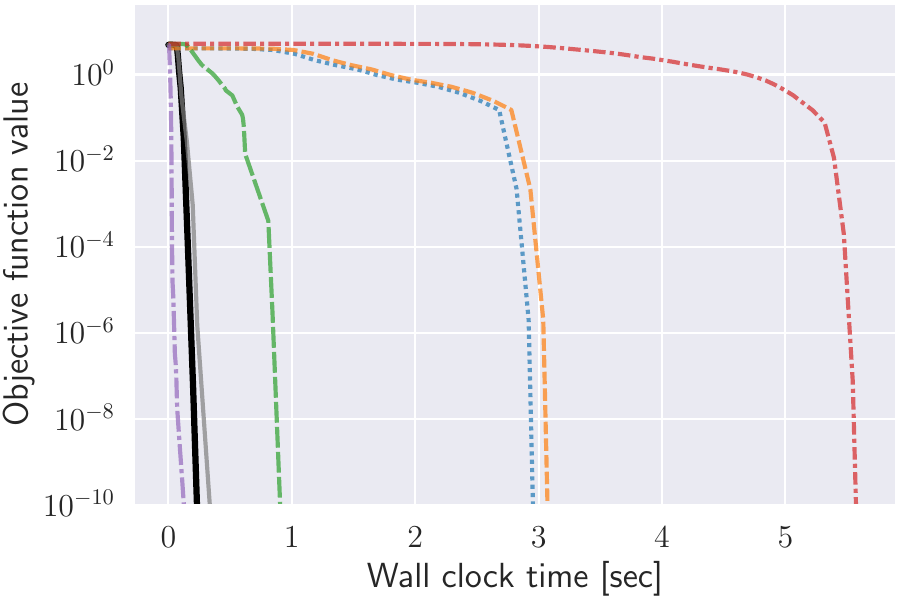}\par
  }\hfill
  \subcaptionbox{$r = 10$, $p = 0.1$\label{fig:exp_nmf_10_01}}[0.32\linewidth]{
    \includegraphics[width=\linewidth]{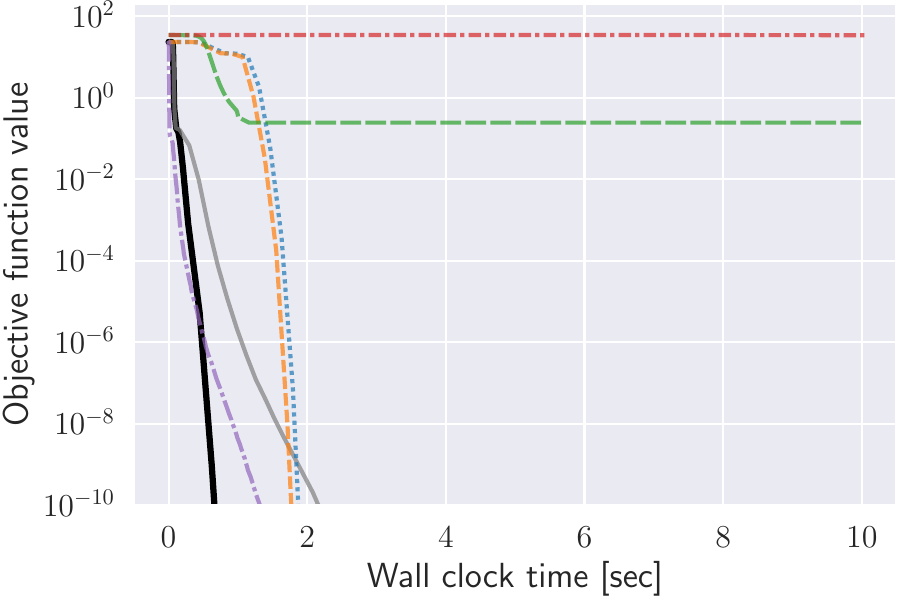}\par
    \includegraphics[width=\linewidth]{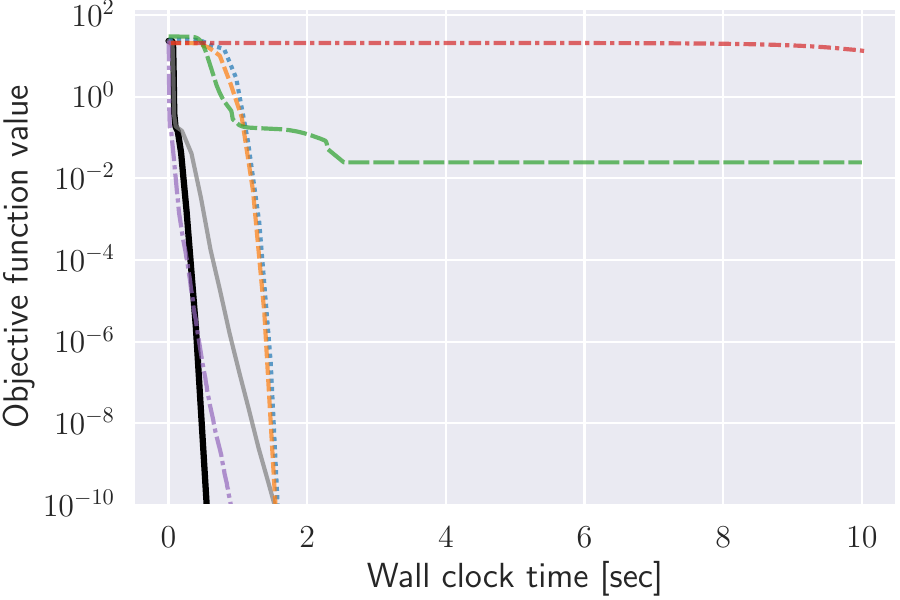}\par
  }\hfill
  \subcaptionbox{$r = 10$, $p = 0.5$\label{fig:exp_nmf_10_05}}[0.32\linewidth]{
    \includegraphics[width=\linewidth]{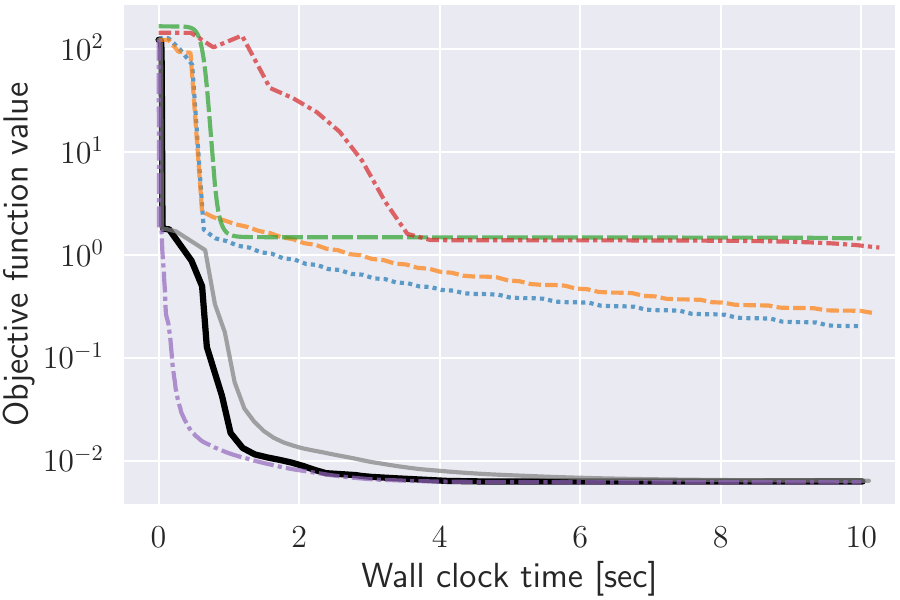}\par
    \includegraphics[width=\linewidth]{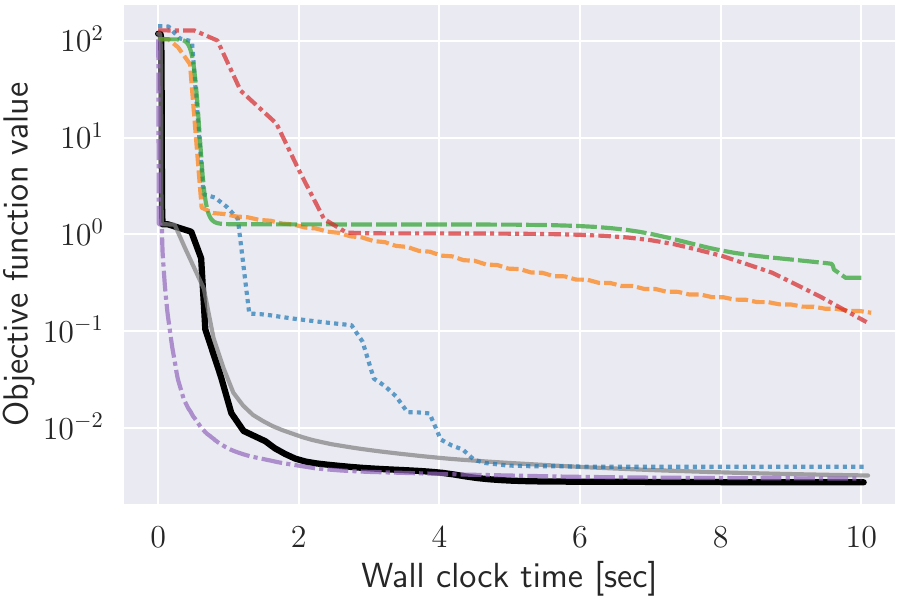}\par
  }\par\bigskip
  \subcaptionbox{$r = 40$, $p = 0.02$\label{fig:exp_nmf_40_002}}[0.32\linewidth]{
    \includegraphics[width=\linewidth]{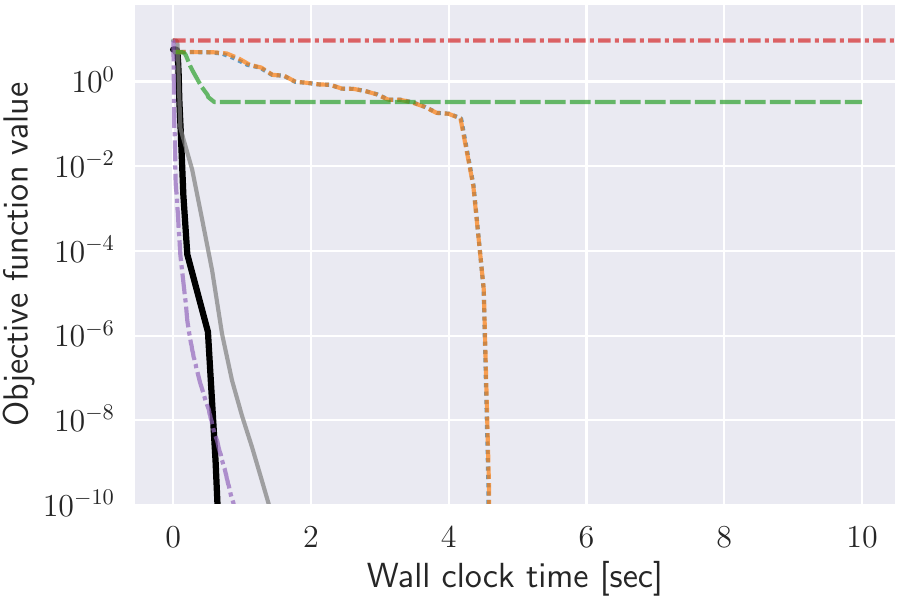}\par
    \includegraphics[width=\linewidth]{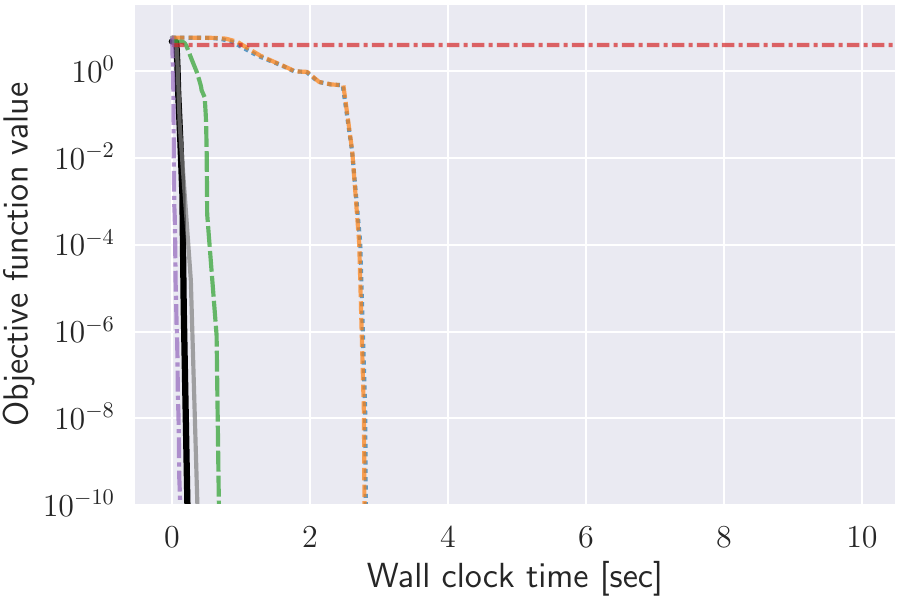}\par
  }\hfill
  \subcaptionbox{$r = 40$, $p = 0.1$\label{fig:exp_nmf_40_01}}[0.32\linewidth]{
    \includegraphics[width=\linewidth]{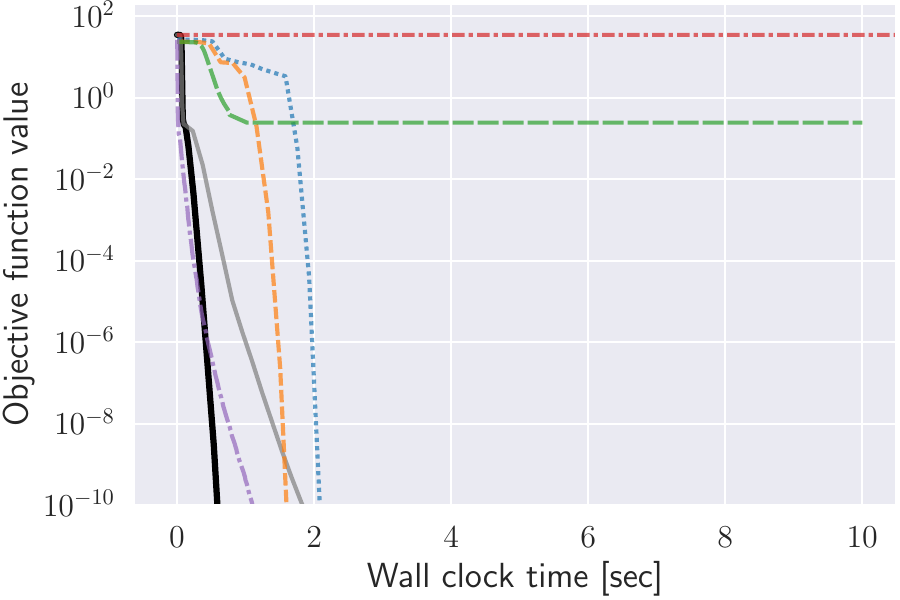}\par
    \includegraphics[width=\linewidth]{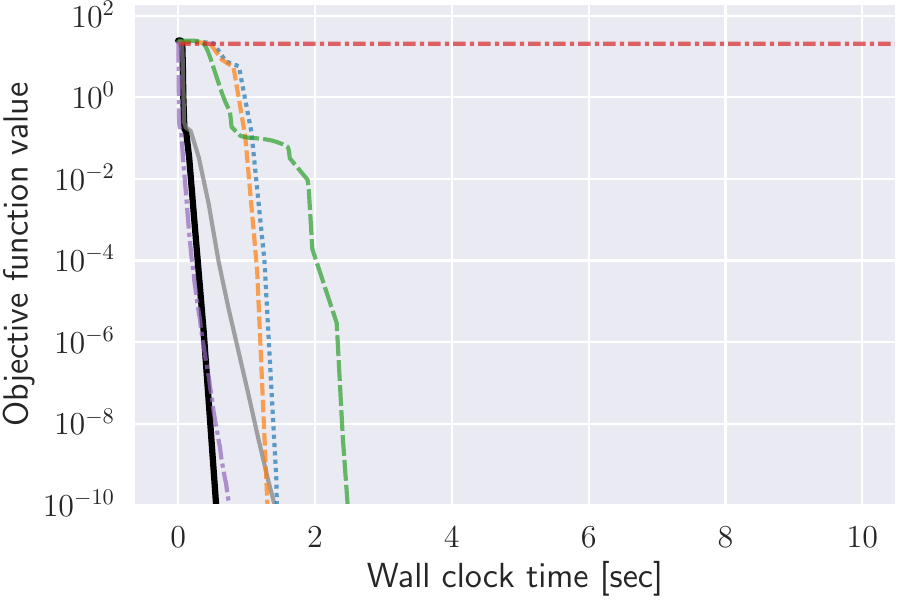}\par
  }\hfill
  \subcaptionbox{$r = 40$, $p = 0.5$\label{fig:exp_nmf_40_05}}[0.32\linewidth]{
    \includegraphics[width=\linewidth]{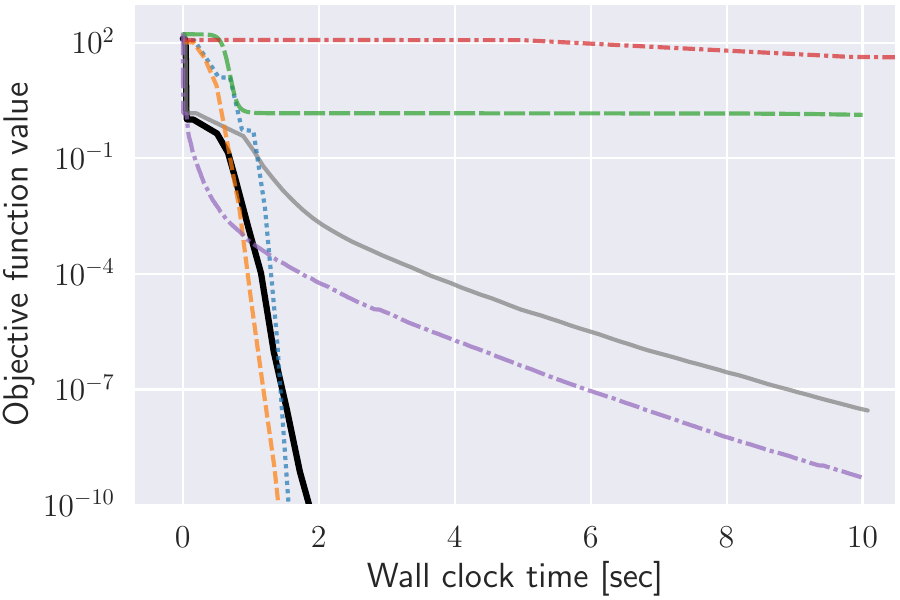}\par
    \includegraphics[width=\linewidth]{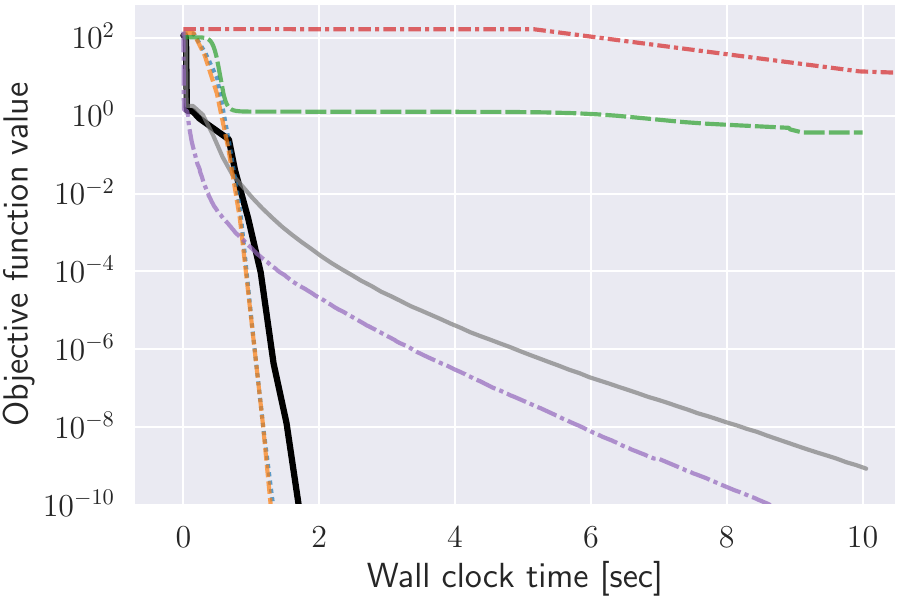}\par
  }\par\medskip
  \caption{Results of NMF (problem~\cref{eq:problem_nmf}).}
  \label{fig:exp_nmf}
\end{figure}

\sisetup{
  table-format = <1.1e+3,
  table-align-comparator = false,
}
\begin{table}
  \caption{Results of compressed sensing (problem~\cref{eq:problem_cs}).}
  \label{table:exp_stats_cs}
  \begin{subtable}[b]{\textwidth}
    \subcaption{$x_{\mathrm{max}} = 0.1$, $d_{\mathrm{nnz}} = 5$}
    \label{table:exp_stats_cs_nnz5_xmax1e-01}
    \centering\footnotesize
    \begin{tabular}{@{}lSSrrrrrr@{}}\toprule
      & {\multirow{2}{*}{objective}} & {\multirow{2}{*}{GM norm}} & \multirow{2}{*}{\pbox[t]{2.1em}{\centering time (sec)}} & \multirow{2}{*}{\#iter} & \multicolumn{3}{c}{\#evaluations} & \multirow{2}{*}{\pbox[t]{3.1em}{\centering success (\%)}}\\\cmidrule(lr){6-8}
      & & & & & $F$ & JVP & $\proj_{\mathcal C}$\\\midrule
      \textbf{Proposed} & 3.5e-14 & 1.4e-06 & 0.55 & 3.4 & 7.8 & 343.2 & 117.8 & 100\\
Proposed-NA & 1.7e-12 & 1.6e-06 & 1.06 & 4.3 & 9.6 & 787.9 & 295.8 & 100\\
Fan & 1.0e-17 & 2.9e-08 & 1.03 & 3.1 & 7.2 & 832.5 & 280.6 & 100\\
KYF & 1.7e-17 & 3.9e-08 & 0.95 & 3.2 & 7.4 & 861.0 & 290.2 & 100\\
Facchinei & 6.7e-14 & 1.3e-06 & 0.38 & 4.9 & 10.8 & 278.7 & 97.8 & 100\\
GGO & 2.2e-13 & 9.4e-06 & 3.31 & 258.5 & 1461.8 & 0.0 & 775.5 & 100\\
PG & 1.1e-11 & 9.1e-06 & 0.59 & 181.2 & 397.5 & 181.2 & 396.5 & 100\\
\bottomrule
    \end{tabular}
  \end{subtable}\par\medskip
  \begin{subtable}[b]{\textwidth}
    \subcaption{$x_{\mathrm{max}} = 0.1$, $d_{\mathrm{nnz}} = 10$}
    \label{table:exp_stats_cs_nnz10_xmax1e-01}
    \centering\footnotesize
    \begin{tabular}{@{}lSSrrrrrr@{}}\toprule
      & {\multirow{2}{*}{objective}} & {\multirow{2}{*}{GM norm}} & \multirow{2}{*}{\pbox[t]{2.1em}{\centering time (sec)}} & \multirow{2}{*}{\#iter} & \multicolumn{3}{c}{\#evaluations} & \multirow{2}{*}{\pbox[t]{3.1em}{\centering success (\%)}}\\\cmidrule(lr){6-8}
      & & & & & $F$ & JVP & $\proj_{\mathcal C}$\\\midrule
      \textbf{Proposed} & 7.6e-13 & 2.1e-06 & 1.11 & 5.0 & 11.1 & 978.6 & 331.2 & 100\\
Proposed-NA & 3.2e-07 & 7.0e-05 & 4.35 & 13.6 & 28.6 & 3934.7 & 1455.3 & 80\\
Fan & 7.2e-14 & 4.3e-07 & 1.69 & 4.2 & 9.4 & 1497.0 & 503.2 & 100\\
KYF & 6.2e-12 & 1.1e-06 & 1.55 & 4.1 & 9.2 & 1469.4 & 493.9 & 100\\
Facchinei & 2.2e-13 & 2.0e-06 & 1.29 & 6.4 & 13.8 & 1167.0 & 395.4 & 100\\
GGO & 4.9e-06 & 1.5e-02 & 8.87 & 594.1 & 3807.1 & 0.0 & 1782.3 & 40\\
PG & 1.4e-07 & 5.4e-05 & 3.93 & 1244.6 & 2686.0 & 1244.6 & 2685.0 & 80\\
\bottomrule
    \end{tabular}
  \end{subtable}\par\medskip
  \begin{subtable}[b]{\textwidth}
    \subcaption{$x_{\mathrm{max}} = 0.1$, $d_{\mathrm{nnz}} = 20$}
    \label{table:exp_stats_cs_nnz20_xmax1e-01}
    \centering\footnotesize
    \begin{tabular}{@{}lSSrrrrrr@{}}\toprule
      & {\multirow{2}{*}{objective}} & {\multirow{2}{*}{GM norm}} & \multirow{2}{*}{\pbox[t]{2.1em}{\centering time (sec)}} & \multirow{2}{*}{\#iter} & \multicolumn{3}{c}{\#evaluations} & \multirow{2}{*}{\pbox[t]{3.1em}{\centering success (\%)}}\\\cmidrule(lr){6-8}
      & & & & & $F$ & JVP & $\proj_{\mathcal C}$\\\midrule
      \textbf{Proposed} & 3.9e-14 & 1.6e-06 & 0.25 & 3.8 & 8.6 & 286.2 & 99.2 & 100\\
Proposed-NA & 4.5e-12 & 3.1e-06 & 1.38 & 5.8 & 12.6 & 1246.6 & 464.8 & 100\\
Fan & 7.0e-15 & 8.1e-07 & 0.64 & 3.1 & 7.2 & 914.4 & 307.9 & 100\\
KYF & 1.9e-15 & 4.8e-07 & 0.66 & 3.1 & 7.2 & 940.2 & 316.5 & 100\\
Facchinei & 1.9e-14 & 8.5e-07 & 0.24 & 6.9 & 14.8 & 239.1 & 86.6 & 100\\
GGO & 1.3e-13 & 8.5e-06 & 2.07 & 168.9 & 1158.6 & 0.0 & 506.7 & 100\\
PG & 1.6e-11 & 9.7e-06 & 1.05 & 329.8 & 717.7 & 329.8 & 716.7 & 100\\
\bottomrule
    \end{tabular}
  \end{subtable}
\end{table}

\begin{table}
  \ContinuedFloat
  \caption{Results of compressed sensing (problem~\cref{eq:problem_cs}).}
  \begin{subtable}[b]{\textwidth}
    \subcaption{$x_{\mathrm{max}} = 1$, $d_{\mathrm{nnz}} = 5$}
    \label{table:exp_stats_cs_nnz5_xmax1e+00}
    \centering\footnotesize
    \begin{tabular}{@{}lSSrrrrrr@{}}\toprule
      & {\multirow{2}{*}{objective}} & {\multirow{2}{*}{GM norm}} & \multirow{2}{*}{\pbox[t]{2.1em}{\centering time (sec)}} & \multirow{2}{*}{\#iter} & \multicolumn{3}{c}{\#evaluations} & \multirow{2}{*}{\pbox[t]{3.1em}{\centering success (\%)}}\\\cmidrule(lr){6-8}
      & & & & & $F$ & JVP & $\proj_{\mathcal C}$\\\midrule
      \textbf{Proposed} & 1.1e-14 & 8.6e-07 & 0.35 & 8.7 & 18.4 & 310.8 & 112.3 & 100\\
Proposed-NA & 3.6e-14 & 1.4e-06 & 0.47 & 9.6 & 20.2 & 409.5 & 160.1 & 100\\
Fan & 2.1e-14 & 1.0e-06 & 1.41 & 5.5 & 16.1 & 1227.4 & 418.6 & 100\\
KYF & 1.9e-14 & 7.5e-07 & 1.41 & 5.7 & 20.6 & 1295.3 & 445.4 & 100\\
Facchinei & 3.0e-14 & 7.1e-07 & 0.66 & 78.1 & 157.2 & 449.4 & 227.9 & 100\\
GGO & 1.6e-13 & 9.7e-06 & 5.62 & 469.4 & 2897.9 & 0.0 & 1408.2 & 100\\
PG & 3.3e-12 & 9.2e-06 & 0.28 & 82.4 & 184.3 & 82.4 & 183.3 & 100\\
\bottomrule
    \end{tabular}
  \end{subtable}\par\medskip
  \begin{subtable}[b]{\textwidth}
    \subcaption{$x_{\mathrm{max}} = 1$, $d_{\mathrm{nnz}} = 10$}
    \label{table:exp_stats_cs_nnz10_xmax1e+00}
    \centering\footnotesize
    \begin{tabular}{@{}lSSrrrrrr@{}}\toprule
      & {\multirow{2}{*}{objective}} & {\multirow{2}{*}{GM norm}} & \multirow{2}{*}{\pbox[t]{2.1em}{\centering time (sec)}} & \multirow{2}{*}{\#iter} & \multicolumn{3}{c}{\#evaluations} & \multirow{2}{*}{\pbox[t]{3.1em}{\centering success (\%)}}\\\cmidrule(lr){6-8}
      & & & & & $F$ & JVP & $\proj_{\mathcal C}$\\\midrule
      \textbf{Proposed} & 9.4e-15 & 1.2e-06 & 0.91 & 23.4 & 48.1 & 789.0 & 286.4 & 100\\
Proposed-NA & 1.2e-12 & 2.4e-06 & 1.80 & 25.5 & 52.3 & 1626.4 & 619.3 & 100\\
Fan & 4.0e-01 & 8.4e-01 & 3.74 & 12.2 & 25.4 & 3465.3 & 1167.3 & 90\\
KYF & 8.0e-01 & 2.5e-01 & 3.74 & 14.9 & 241.5 & 3190.4 & 1283.2 & 90\\
Facchinei & 1.7e-01 & 2.6e-01 & 3.25 & 675.5 & 1352.0 & 2371.2 & 1465.9 & 90\\
GGO & 1.4e+00 & 1.2e+01 & 8.68 & 673.8 & 4538.9 & 0.0 & 2021.4 & 50\\
PG & 7.7e-12 & 9.5e-06 & 0.98 & 305.2 & 665.2 & 305.2 & 664.2 & 100\\
\bottomrule
    \end{tabular}
  \end{subtable}\par\medskip
  \begin{subtable}[b]{\textwidth}
    \subcaption{$x_{\mathrm{max}} = 1$, $d_{\mathrm{nnz}} = 20$}
    \label{table:exp_stats_cs_nnz20_xmax1e+00}
    \centering\footnotesize
    \begin{tabular}{@{}lSSrrrrrr@{}}\toprule
      & {\multirow{2}{*}{objective}} & {\multirow{2}{*}{GM norm}} & \multirow{2}{*}{\pbox[t]{2.1em}{\centering time (sec)}} & \multirow{2}{*}{\#iter} & \multicolumn{3}{c}{\#evaluations} & \multirow{2}{*}{\pbox[t]{3.1em}{\centering success (\%)}}\\\cmidrule(lr){6-8}
      & & & & & $F$ & JVP & $\proj_{\mathcal C}$\\\midrule
      \textbf{Proposed} & 1.1e-01 & 4.9e-05 & 5.97 & 94.2 & 195.0 & 5390.1 & 1890.9 & 80\\
Proposed-NA & 1.1e-01 & 1.6e-02 & 9.73 & 59.5 & 120.7 & 8925.8 & 3319.3 & 20\\
Fan & 4.6e-01 & 2.3e+00 & 9.47 & 23.3 & 479.6 & 8310.6 & 3214.7 & 10\\
KYF & 4.7e-01 & 2.4e+00 & 9.24 & 22.6 & 484.0 & 8091.7 & 3146.7 & 20\\
Facchinei & 2.6e-01 & 4.1e-01 & 9.24 & 1806.5 & 3614.0 & 6916.2 & 4111.9 & 20\\
GGO & 1.4e+00 & 1.8e+01 & 9.47 & 820.6 & 6427.6 & 0.0 & 2461.8 & 10\\
PG & 1.9e-01 & 6.2e-03 & 8.53 & 2737.1 & 5899.1 & 2737.1 & 5898.1 & 30\\
\bottomrule
    \end{tabular}
  \end{subtable}
\end{table}

\begin{table}
  \caption{Results of NMF (problem~\cref{eq:problem_nmf}).}
  \label{table:exp_stats_nmf}
  \begin{subtable}[b]{\textwidth}
    \subcaption{$r = 10$, $p = 0.02$}
    \label{table:exp_stats_nmf_r10_p2e-02}
    \centering\footnotesize
    \begin{tabular}{@{}lSSrrrrrr@{}}\toprule
      & {\multirow{2}{*}{objective}} & {\multirow{2}{*}{GM norm}} & \multirow{2}{*}{\pbox[t]{2.1em}{\centering time (sec)}} & \multirow{2}{*}{\#iter} & \multicolumn{3}{c}{\#evaluations} & \multirow{2}{*}{\pbox[t]{3.1em}{\centering success (\%)}}\\\cmidrule(lr){6-8}
      & & & & & $F$ & JVP & $\proj_{\mathcal C}$\\\midrule
      \textbf{Proposed} & 3.1e-13 & 1.6e-07 & 0.46 & 36.4 & 75.0 & 908.4 & 339.2 & 100\\
Proposed-NA & 3.6e-10 & 3.7e-06 & 0.65 & 36.8 & 75.3 & 1414.3 & 549.1 & 100\\
Fan & 2.3e-11 & 3.0e-06 & 3.77 & 25.6 & 116.3 & 8170.9 & 2785.2 & 100\\
KYF & 2.4e-11 & 3.1e-06 & 3.84 & 25.6 & 116.3 & 8528.8 & 2904.5 & 100\\
Facchinei & 1.2e-01 & 1.6e-01 & 6.08 & 2650.4 & 5301.8 & 8645.7 & 5532.3 & 50\\
GGO & 2.1e-11 & 2.8e-06 & 6.90 & 85.0 & 86.0 & 0.0 & 170.0 & 100\\
PG & 1.3e-09 & 9.8e-06 & 0.25 & 256.2 & 552.8 & 256.2 & 551.8 & 100\\
\bottomrule
    \end{tabular}
  \end{subtable}\par\medskip
  \begin{subtable}[b]{\textwidth}
    \subcaption{$r = 10$, $p = 0.1$}
    \label{table:exp_stats_nmf_r10_p1e-01}
    \centering\footnotesize
    \begin{tabular}{@{}lSSrrrrrr@{}}\toprule
      & {\multirow{2}{*}{objective}} & {\multirow{2}{*}{GM norm}} & \multirow{2}{*}{\pbox[t]{2.1em}{\centering time (sec)}} & \multirow{2}{*}{\#iter} & \multicolumn{3}{c}{\#evaluations} & \multirow{2}{*}{\pbox[t]{3.1em}{\centering success (\%)}}\\\cmidrule(lr){6-8}
      & & & & & $F$ & JVP & $\proj_{\mathcal C}$\\\midrule
      \textbf{Proposed} & 7.8e-12 & 1.3e-06 & 0.69 & 36.0 & 73.1 & 1383.9 & 497.3 & 100\\
Proposed-NA & 1.3e-09 & 6.7e-06 & 1.68 & 42.2 & 85.4 & 3785.0 & 1421.6 & 100\\
Fan & 1.8e-11 & 2.1e-06 & 1.72 & 10.4 & 35.1 & 3619.9 & 1223.4 & 100\\
KYF & 9.1e-12 & 1.2e-06 & 1.62 & 10.4 & 34.4 & 3579.9 & 1209.5 & 100\\
Facchinei & 1.7e-01 & 1.8e-01 & 10.00 & 4520.2 & 9041.4 & 13965.0 & 9175.2 & 0\\
GGO & 2.2e+01 & 3.7e+00 & 10.06 & 89.2 & 90.2 & 0.0 & 178.4 & 0\\
PG & 2.8e-09 & 9.9e-06 & 0.85 & 883.1 & 1903.2 & 883.1 & 1902.2 & 100\\
\bottomrule
    \end{tabular}
  \end{subtable}\par\medskip
  \begin{subtable}[b]{\textwidth}
    \subcaption{$r = 10$, $p = 0.5$}
    \label{table:exp_stats_nmf_r10_p5e-01}
    \centering\footnotesize
    \begin{tabular}{@{}lSSrrrrrr@{}}\toprule
      & {\multirow{2}{*}{objective}} & {\multirow{2}{*}{GM norm}} & \multirow{2}{*}{\pbox[t]{2.1em}{\centering time (sec)}} & \multirow{2}{*}{\#iter} & \multicolumn{3}{c}{\#evaluations} & \multirow{2}{*}{\pbox[t]{3.1em}{\centering success (\%)}}\\\cmidrule(lr){6-8}
      & & & & & $F$ & JVP & $\proj_{\mathcal C}$\\\midrule
      \textbf{Proposed} & 3.8e-03 & 1.0e-04 & 9.71 & 174.7 & 368.2 & 20065.2 & 6863.1 & 20\\
Proposed-NA & 4.0e-03 & 8.3e-04 & 10.07 & 120.2 & 242.4 & 22780.5 & 8440.7 & 0\\
Fan & 9.4e-02 & 4.6e-01 & 10.07 & 62.6 & 1380.9 & 21825.5 & 8522.6 & 0\\
KYF & 2.2e-01 & 9.0e-01 & 10.12 & 62.8 & 1646.0 & 21834.9 & 8778.1 & 0\\
Facchinei & 8.9e-01 & 3.8e-01 & 10.00 & 4587.8 & 9176.6 & 13808.1 & 9190.5 & 0\\
GGO & 6.4e-01 & 4.5e-01 & 10.17 & 28.0 & 32.7 & 0.0 & 59.9 & 0\\
PG & 4.0e-03 & 3.7e-04 & 10.00 & 10338.3 & 22252.0 & 10338.3 & 22251.0 & 0\\
\bottomrule
    \end{tabular}
  \end{subtable}
\end{table}

\begin{table}
  \ContinuedFloat
  \caption{Results of NMF (problem~\cref{eq:problem_nmf}).}
  \begin{subtable}[b]{\textwidth}
    \subcaption{$r = 40$, $p = 0.02$}
    \label{table:exp_stats_nmf_r40_p2e-02}
    \centering\footnotesize
    \begin{tabular}{@{}lSSrrrrrr@{}}\toprule
      & {\multirow{2}{*}{objective}} & {\multirow{2}{*}{GM norm}} & \multirow{2}{*}{\pbox[t]{2.1em}{\centering time (sec)}} & \multirow{2}{*}{\#iter} & \multicolumn{3}{c}{\#evaluations} & \multirow{2}{*}{\pbox[t]{3.1em}{\centering success (\%)}}\\\cmidrule(lr){6-8}
      & & & & & $F$ & JVP & $\proj_{\mathcal C}$\\\midrule
      \textbf{Proposed} & 1.3e-10 & 2.3e-06 & 0.41 & 33.8 & 69.7 & 765.3 & 288.9 & 100\\
Proposed-NA & 3.5e-10 & 3.0e-06 & 0.67 & 34.5 & 70.5 & 1409.4 & 545.7 & 100\\
Fan & 1.8e-11 & 1.6e-06 & 3.58 & 22.9 & 99.3 & 7309.7 & 2487.7 & 100\\
KYF & 1.8e-11 & 1.6e-06 & 3.59 & 22.9 & 99.3 & 7534.4 & 2562.6 & 100\\
Facchinei & 1.6e-01 & 2.2e-01 & 7.29 & 3119.5 & 6240.0 & 9808.5 & 6389.0 & 30\\
GGO & 5.5e+00 & 3.1e-02 & 11.30 & 4.0 & 5.0 & 0.0 & 8.0 & 0\\
PG & 1.6e-09 & 9.6e-06 & 0.26 & 254.8 & 549.9 & 254.8 & 548.9 & 100\\
\bottomrule
    \end{tabular}
  \end{subtable}\par\medskip
  \begin{subtable}[b]{\textwidth}
    \subcaption{$r = 40$, $p = 0.1$}
    \label{table:exp_stats_nmf_r40_p1e-01}
    \centering\footnotesize
    \begin{tabular}{@{}lSSrrrrrr@{}}\toprule
      & {\multirow{2}{*}{objective}} & {\multirow{2}{*}{GM norm}} & \multirow{2}{*}{\pbox[t]{2.1em}{\centering time (sec)}} & \multirow{2}{*}{\#iter} & \multicolumn{3}{c}{\#evaluations} & \multirow{2}{*}{\pbox[t]{3.1em}{\centering success (\%)}}\\\cmidrule(lr){6-8}
      & & & & & $F$ & JVP & $\proj_{\mathcal C}$\\\midrule
      \textbf{Proposed} & 1.7e-11 & 1.4e-06 & 0.63 & 33.3 & 67.7 & 1187.4 & 429.1 & 100\\
Proposed-NA & 1.3e-09 & 7.1e-06 & 1.45 & 38.0 & 77.0 & 3102.2 & 1168.2 & 100\\
Fan & 2.8e-12 & 1.2e-06 & 1.69 & 9.7 & 33.8 & 3356.5 & 1135.4 & 100\\
KYF & 5.4e-12 & 1.3e-06 & 1.49 & 9.1 & 29.1 & 3105.2 & 1048.2 & 100\\
Facchinei & 1.8e-01 & 1.8e-01 & 9.27 & 3971.2 & 7943.4 & 12461.4 & 8125.0 & 10\\
GGO & 2.6e+01 & 9.4e-02 & 12.53 & 3.9 & 4.9 & 0.0 & 7.8 & 0\\
PG & 2.2e-09 & 9.9e-06 & 0.73 & 715.5 & 1542.4 & 715.5 & 1541.4 & 100\\
\bottomrule
    \end{tabular}
  \end{subtable}\par\medskip
  \begin{subtable}[b]{\textwidth}
    \subcaption{$r = 40$, $p = 0.5$}
    \label{table:exp_stats_nmf_r40_p5e-01}
    \centering\footnotesize
    \begin{tabular}{@{}lSSrrrrrr@{}}\toprule
      & {\multirow{2}{*}{objective}} & {\multirow{2}{*}{GM norm}} & \multirow{2}{*}{\pbox[t]{2.1em}{\centering time (sec)}} & \multirow{2}{*}{\#iter} & \multicolumn{3}{c}{\#evaluations} & \multirow{2}{*}{\pbox[t]{3.1em}{\centering success (\%)}}\\\cmidrule(lr){6-8}
      & & & & & $F$ & JVP & $\proj_{\mathcal C}$\\\midrule
      \textbf{Proposed} & 5.9e-10 & 5.4e-06 & 1.67 & 56.0 & 115.8 & 3066.0 & 1078.0 & 100\\
Proposed-NA & 8.0e-09 & 1.1e-05 & 9.51 & 109.6 & 220.5 & 20514.7 & 7602.7 & 80\\
Fan & 4.6e-11 & 2.6e-06 & 1.38 & 8.1 & 22.3 & 2863.2 & 964.8 & 100\\
KYF & 4.5e-11 & 1.6e-06 & 1.36 & 8.1 & 21.2 & 2835.5 & 954.6 & 100\\
Facchinei & 7.8e-01 & 4.1e-01 & 10.00 & 4402.4 & 8805.8 & 13251.9 & 8819.7 & 0\\
GGO & 2.7e+01 & 2.3e+01 & 13.74 & 2.8 & 3.8 & 0.0 & 7.6 & 0\\
PG & 7.3e-09 & 1.0e-05 & 6.83 & 6750.1 & 14531.0 & 6750.1 & 14530.0 & 100\\
\bottomrule
    \end{tabular}
  \end{subtable}
\end{table}

\begin{figure}[thp]
  \centering
  \subcaptionbox{$r = 10$, $p = 0.02$}[0.32\linewidth]{
    \includegraphics[width=\linewidth]{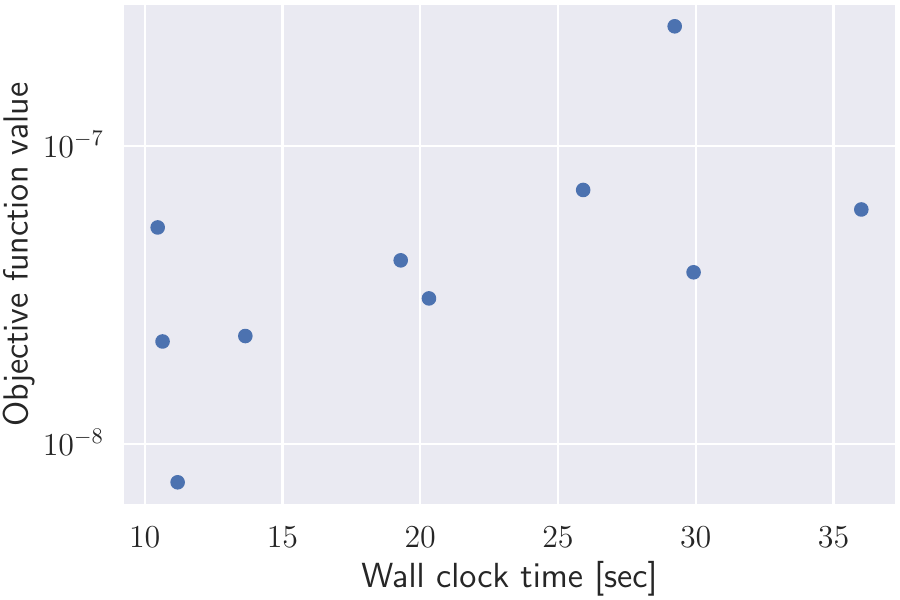}\par
  }\hfill
  \subcaptionbox{$r = 10$, $p = 0.1$}[0.32\linewidth]{
    \includegraphics[width=\linewidth]{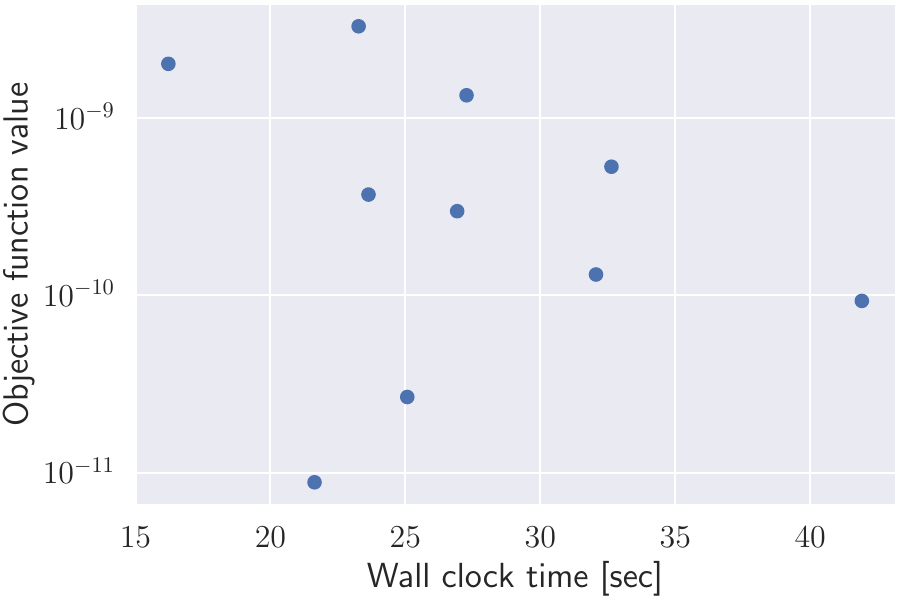}\par
  }\hfill
  \subcaptionbox{$r = 10$, $p = 0.5$}[0.32\linewidth]{
    \includegraphics[width=\linewidth]{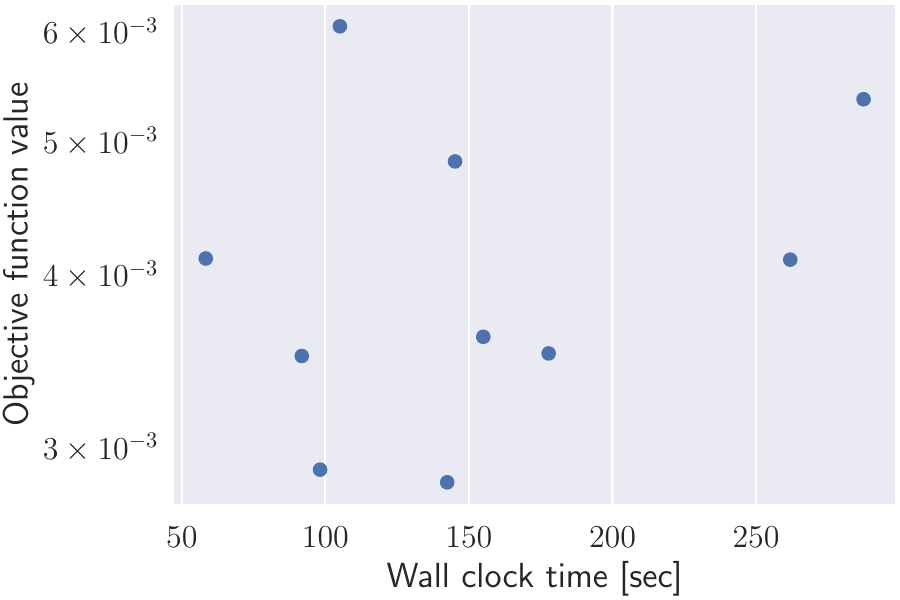}\par
  }\par\bigskip
  \subcaptionbox{$r = 40$, $p = 0.02$\label{fig:exp_nmf_trf_40_002}}[0.32\linewidth]{
    \includegraphics[width=\linewidth]{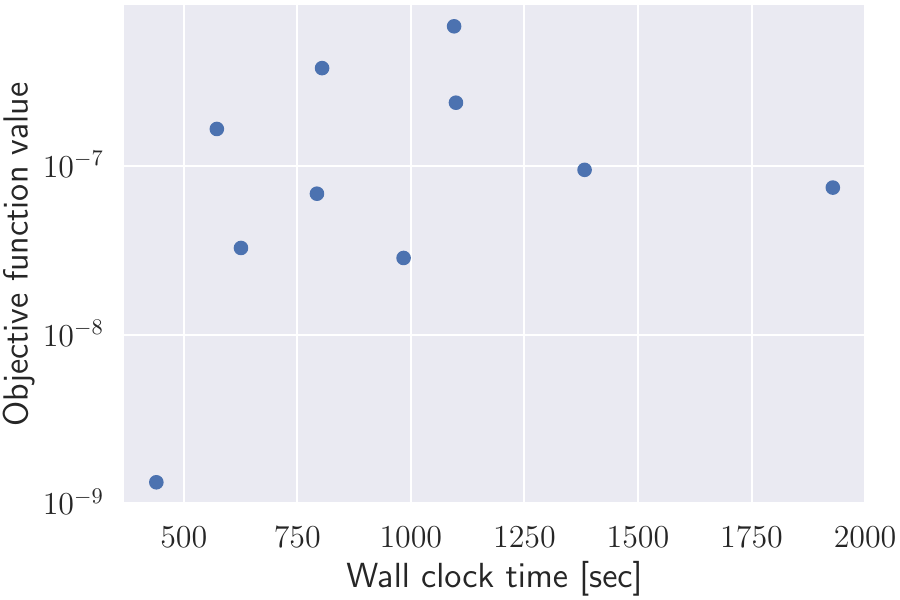}\par
  }\hfill
  \subcaptionbox{$r = 40$, $p = 0.1$\label{fig:exp_nmf_trf_40_01}}[0.32\linewidth]{
    \includegraphics[width=\linewidth]{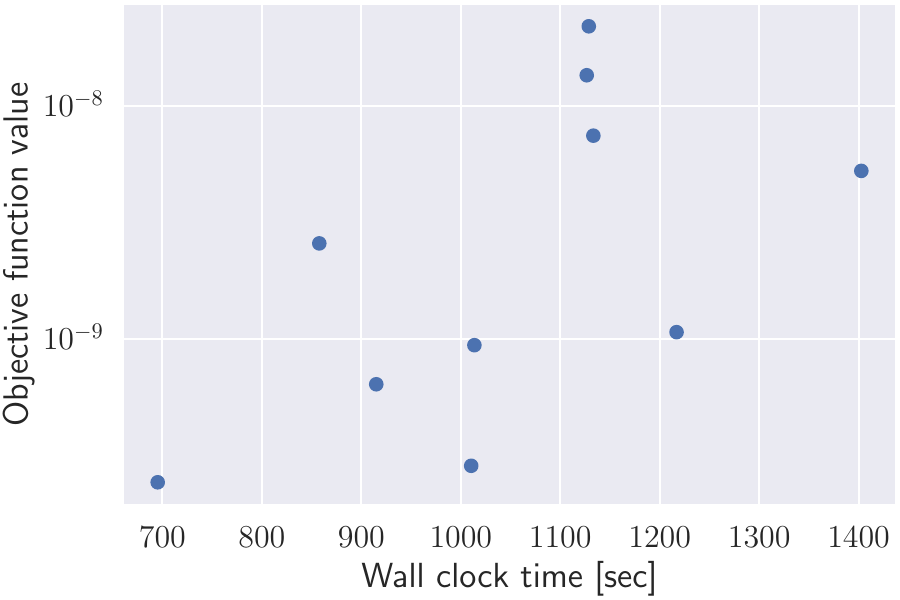}\par
  }\hfill
  \subcaptionbox{$r = 40$, $p = 0.5$\label{fig:exp_nmf_trf_40_05}}[0.32\linewidth]{
    \includegraphics[width=\linewidth]{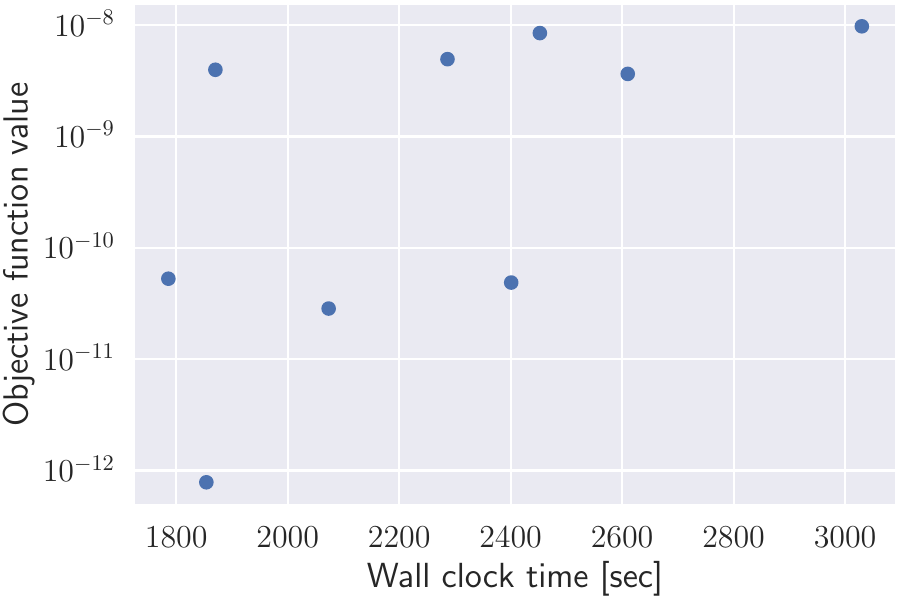}\par
  }\par\medskip
  \caption{Results of NMF (problem~\cref{eq:problem_nmf}) by the TRF method.}
  \label{fig:exp_nmf_trf}
\end{figure}

\subsection{Results}
\subsubsection{Compressed sensing and NMF}
\cref{fig:exp_cs,fig:exp_nmf} show the results of compressed sensing in \cref{eq:problem_cs} and NMF in \cref{eq:problem_nmf}.
Each figure consists of six subfigures, and they consist of two plots; the upper one shows the worst case among ten randomly generated instances, and the lower one shows the best case.\footnote{
  Here, for each method, we determine the best and worst cases out of ten instances as follows.
  Each algorithm is stopped when either of the following conditions is fulfilled: (i) the objective function value falls below $10^{-10}$; (ii) the execution time exceeds 10 seconds.
  First, we define that the case stopped by condition (i) is better than that stopped by condition (ii).
  Then, among the cases stopped by condition (i), the case with a shorter execution time is defined as better.
  Similarly, among the cases stopped by condition (ii), the case with a smaller objective value is defined as better.
  Note that from the above definition, the instances corresponding to plots in the same figure may be distinct.
}

\cref{table:exp_stats_cs,table:exp_stats_nmf} provide more detailed information.
For the tables, each algorithm is stopped when either of the following conditions is fulfilled: (i) the algorithm finds a point where the norm of the gradient mapping is less than $10^{-5}$; (ii) the execution time exceeds 10 seconds.
The ``success'' column indicates the percentage of instances (out of 10) that ended up satisfying condition (i).
The other columns show the averages of the following values: the objective function value reached, the gradient-mapping norm, the execution time, the number of iterations, and the number of basic operations.
JVP stands for Jacobian-vector products.

A remarkable feature of our method is its stability in addition to fast convergence.
For example, while the Fan and KYF methods perform well in most cases, they sometimes do not converge fast, as shown in \cref{table:exp_stats_cs_nnz10_xmax1e+00,table:exp_stats_cs_nnz20_xmax1e+00}.
The proposed method shows the best or comparable performance in all our settings compared to the other methods.
This suggests that our method is stable without careful parameter tuning.


As seen from \cref{table:exp_stats_nmf_r40_p2e-02,table:exp_stats_nmf_r40_p1e-01,table:exp_stats_nmf_r40_p5e-01}, the Facchinei and GGO methods do not work well in some cases.
As for the Facchinei method, the reason is presumably that the method does not guarantee global convergence.
For GGO, it is observed from the tables that the number of iterations performed within the time limit is small, say 3 or 4.
It is because the method at each iteration computes a Jacobian explicitly and solves a linear system, resulting in a high cost per iteration for large-scale problems.
Our method guarantees global convergence and repeats relatively low-cost iterations without Jacobian computation, which also leads to a stable performance.


\cref{fig:exp_nmf_trf} shows the results of the TRF method.
Since this method can only handle box constraints, the results only of problem~\cref{eq:problem_nmf} are presented.
One marker corresponds to one instance, representing the elapsed time and the obtained objective value.\footnote{
  We ran the TRF method for ten instances for each $(r, p)$, but \cref{fig:exp_nmf_trf_40_05} has only nine markers because the algorithm stopped with the error ``\texttt{SVD did not converge}'' for one instance.
}
TRF takes more time to converge than the proposed method; in particular, comparing \cref{table:exp_stats_nmf_r40_p5e-01,fig:exp_nmf_trf_40_05}, we see that the elapsed time is about 1000 times longer than ours.
This result may be due to the difference in how TRF and ours handle the constraint.
When the optimal solution or a stationary point is at the boundary of the constraint set, our method can reach the boundary in a finite number of iterations.
However, TRF does not, as it is an interior point method.

\begin{figure}
  \centering
  \includegraphics[width=0.9\linewidth]{figure/legend.pdf}\par\medskip
  \includegraphics[width=0.5\linewidth]{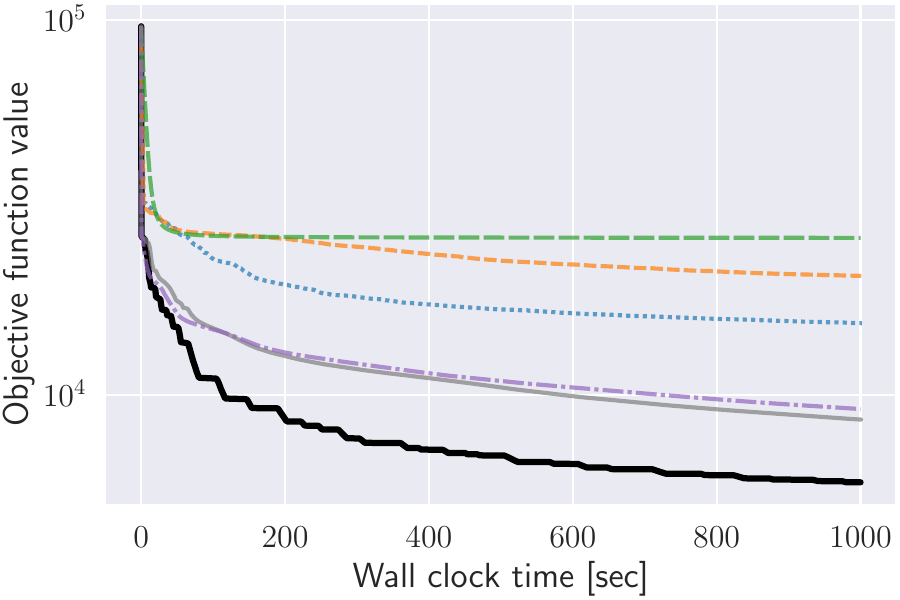}
  \caption{
    Results of autoencoder with MNIST (problem~\cref{eq:problem_autoencoder}).
  }
  \label{fig:exp_ae_obj}
\end{figure}

\subsubsection{Autoencoder with MNIST}
\cref{fig:exp_ae_obj} shows the results of problem~\cref{eq:problem_autoencoder}.
The results of the GGO method are omitted because the method explicitly computes the Jacobian, but it was infeasible in this large-scale setting, where $d = 96{,}104$ and $n = 728{,}000$.
Among the existing methods, the PG method converges the fastest, but the proposed method converges about five times faster than PG.
This result suggests that our method is also effective for large-scale and highly nonlinear problems.

\section{Conclusion and future work}
\label{sec:conclusion}

We proposed an LM method for solving constrained least-squares problems.
Our method finds an $\epsilon$-stationary point of (possibly) nonzero-residual problems
after $O(\epsilon^{-2})$ computation,
and also achieves local quadratic convergence for zero-residual problems.
There are few LM methods having both overall complexity bounds 
and local quadratic convergence even for unconstrained problems;
in fact, our investigation yielded only one such algorithm~\citep{bergou2020convergence}.
The key to our analysis is a simple update rule for $(\lambda_k)$ and the majorization lemma (\cref{lem: majorization LM}).

We may be able to extend the convergence analysis shown in this paper 
to different problem settings. 
For example, it would be interesting to derive an overall complexity bound of LM methods for a nonsmooth function $F$.
It would be also interesting to integrate a stochastic technique into our LM method against problems with $F$ of a huge size.
Finally, in recent years, studies on local convergence analysis for non-zero residual problems are progressive~\citep{ipsen2011rank,behling2019local,bergou2020convergence}.
It is important to research our LM method further in this line.


\appendix
\section{Lemmas and proofs}
\subsection{Lemma on Lipschitz-like properties}
Recall that the line segment $\mathcal L(a, b)$ is defined in \cref{eq: def of L}.
\begin{lemma}\label{lem: Lip implies xxx}
Let $\mathcal X \subseteq \R^d$ be any (possibly nonconvex) set. 
For some constants $\sigma, L > 0$, consider the following two sets of conditions:
 \begin{subequations}
 \begin{align}
 \label{item-cond: J bounded on X}
 &\|J(x)\| \leq \sigma, \quad \forall x \in \mathcal X,\\
 \label{item-cond: J Lip on X}
 &\|J(y) - J(x)\| \leq L\|y - x\|,
 \quad
 \forall x,y \in \mathcal X \ \ \text{s.t.} \ \ \mathcal L(x, y) \subseteq \mathcal X,
 \end{align}
 \end{subequations}
 and
 \begin{subequations}
 \begin{align}
 &
 \label{eq: F Lip-like}
 \|F(y) - F(x)\| \leq \sigma \|y - x\|,
 \quad
 \forall x,y \in \mathcal X \ \ \text{s.t.} \ \ \mathcal L(x, y) \subseteq \mathcal X,\\
 &
 \label{eq: F smooth-like}
 \|F(y) - F(x) - J(x)(y - x)\|
 \leq
 \frac{L}{2}\|y - x\|^2,
 \quad
 \forall x,y \in \mathcal X \ \ \text{s.t.} \ \ \mathcal L(x, y) \subseteq \mathcal X,\\
 &
 \label{eq: grad_f Lip-like}
 \|\nabla f(y) - \nabla f(x)\|
 \leq
 (\sigma^2 + L\|F(x)\|) \|y - x\|,
 \quad
 \forall x,y \in \mathcal X \ \ \text{s.t.} \ \ \mathcal L(x, y) \subseteq \mathcal X.
 \end{align}
 \end{subequations}
 Then,
 \begin{enuminlem}
 \item
 \label{item-lem: bounded J implies xxx}
 \cref{item-cond: J bounded on X}%
 \quad$\Longrightarrow$\quad%
 \cref{eq: F Lip-like},
 
 \item
 \label{item-lem: Lip J implies xxx}
 \cref{item-cond: J Lip on X}%
 \quad$\Longrightarrow$\quad%
 \cref{eq: F smooth-like},
 
 \item
 \label{item-lem: bounded Lip J implies xxx}
 \cref{item-cond: J bounded on X,item-cond: J Lip on X}%
 \quad$\Longrightarrow$\quad%
 \cref{eq: grad_f Lip-like}.
 \end{enuminlem}
\end{lemma}

\begin{proof}
By applying the multivariate mean value theorem, i.e.,
 \begin{equation}
 F(y) - F(x)
 =
 \int_0^1 J((1-\theta) x + \theta y)(y - x)\,d\theta,
 \label{eq: multivariate mean value}
 \end{equation}
we can easily obtain \cref{item-lem: bounded J implies xxx,item-lem: Lip J implies xxx}. \cref{item-lem: bounded Lip J implies xxx} is obtained as follows:
 \begin{alignat}{2}
 \|\nabla f(y) - \nabla f(x)\|
 &=
 \|J(y)^\top F(y) - J(x)^\top F(x)\|\\
 &\leq
 \|J(y)^\top F(y) - J(y)^\top F(x)\| + \|J(y)^\top F(x) - J(x)^\top F(x)\|\\
 &\leq
 \|J(y)\| \|F(y) - F(x)\| + \|J(y) - J(x)\| \|F(x)\|\\
 &\leq
 (\sigma^2 + L\|F(x)\|) \|y - x\|.
 \end{alignat}
The last inequality follows from \cref{item-cond: J bounded on X}, \cref{item-cond: J Lip on X}, and \cref{item-lem: bounded J implies xxx}.
\end{proof}

\begin{remark}
 \label{rem: Lip grad f implies xxx}
By replacing $(F, J)$ with $(f, \nabla f)$ in \cref{item-lem: Lip J implies xxx},
 \begin{equation}
 \|\nabla f(y) - \nabla f(x)\| \leq L_f \|y - x\|,
 \quad
 \forall x,y \in \mathcal X \ \ \text{s.t.} \ \ \mathcal L(x, y) \subseteq \mathcal X
 \end{equation}
implies
 \begin{equation}
 |f(y) - f(x) - \inner{\nabla f(x)}{y - x}| \leq \frac{L_f}{2}\|y - x\|^2,
 \quad
 \forall x,y \in \mathcal X \ \ \text{s.t.} \ \ \mathcal L(x, y) \subseteq \mathcal X.
 \end{equation}
\end{remark}

\subsection{Proof of \texorpdfstring{\cref{lem: majorization LM}}{}}
\label{sec: proof of MM lemma}
The proof requires the following lemma, 
which is useful for deriving the majorization lemma for (general) MM-based algorithms
under the assumption of Lipschitz continuity only on a sublevel set. 
We will use this lemma to prove \cref{lem: majorization LM} as well as \cref{prop: PG decrease}.

\begin{lemma}
 \label{lem: MM technical}
Let $\mathcal X \subseteq \R^d$ be any convex set.
Fix a point $x_k \in \mathcal X$,
and a strictly convex function $\tilde m: \R^d \to \R$ 
such that $\tilde m(x_k) = f(x_k)$ and $\nabla \tilde m(x_k) = \nabla f(x_k)$.
We consider three subsets of $\mathcal X$:
\begin{equation}
  \label{eq: def of R1R2R3}
  \begin{aligned}
    \mathcal R_1 &\coloneqq \Set{x \in \mathcal X}{\tilde m(x) \leq \tilde m(x_k)},\\
    \mathcal R_2 &\coloneqq \Set{x \in \mathcal X}{\mathcal L(x_k, x) \subseteq \mathcal S(x_k)},\\
    \mathcal R_3 &\coloneqq \Set{x \in \mathcal X}{f(x) \leq \tilde m(x)}.
  \end{aligned}
\end{equation}
If
\begin{equation}
  \label{eq: R1 cap R2 in R3}
  (\mathcal R_1 \cap \mathcal R_2) \subseteq \mathcal R_3,
\end{equation}
then $\mathcal R_1 \subseteq \mathcal R_2$,
and therefore $\mathcal R_1 \subseteq \mathcal R_2 \subseteq \mathcal R_3$.
\end{lemma}

\begin{proof}
We fix $x \in \mathcal R_1$ arbitrarily and will show $x \in \mathcal R_2$.
This is obvious if $x = x_k$, and thus, let $x \neq x_k$ below.
Accordingly, we have
\begin{alignat}{2}
  \inner{\nabla f(x_k)}{x - x_k}
  < 0
  \label{eq:first_coef_taylor_f_negative}
\end{alignat}
since
\begin{alignat}{2}
  \inner{\nabla f(x_k)}{x - x_k}
  &=
  \inner{\nabla \tilde m(x_k)}{x - x_k}
  &\quad&\text{(by $\nabla \tilde m(x_k) = \nabla f(x_k)$)}\\
  &<
  \tilde m(x) - \tilde m(x_k)
  &\quad&\text{(by the strictly convexity of $\tilde m$ and $x \neq x_k$)}\\
  &\leq 0
  &\quad&\text{(by $x \in \mathcal R_1$)}.
\end{alignat}
By the Taylor expansion
\begin{alignat}{2}
  f((1 - \theta)x_k + \theta x)
  &=
  f(x_k)
  + \theta \inner{\nabla f(x_k)}{x - x_k}
  + o(\theta)
\end{alignat}
together with \cref{eq:first_coef_taylor_f_negative},
there exists $\theta_1 > 0$ such that
\begin{equation}
  \label{eq:exists_theta2}
  f((1 - \theta)x_k + \theta x) < f(x_k),\quad
  \text{for all $\theta \in (0, \theta_1]$.}
\end{equation}
We will prove $x \in \mathcal R_2$ by contradiction; 
assume $x \notin \mathcal R_2$, i.e., there exists $\theta_2 \in [0, 1]$ such that
\begin{equation}
  \label{eq:exists_theta1}
  f((1 - \theta_2)x_k + \theta_2 x) > f(x_k).
\end{equation}
Combining \cref{eq:exists_theta1,eq:exists_theta2} with the intermediate value theorem
yields that there exists $\theta_3 \in (\theta_1, \theta_2)$ such that
\begin{subequations}
 \begin{align}
 &
 f((1 - \theta)x_k + \theta x) = f(x_k),
 \quad
 \text{for $\theta = \theta_3$}
 \label{eq: t2 such that xxx1}\\
 &
 f((1 - \theta)x_k + \theta x) \leq f(x_k),
 \quad
 \text{for all $\theta \in [0, \theta_3]$}.
 \label{eq: t2 such that xxx2}
 \end{align}
 \end{subequations}
Note that \cref{eq: t2 such that xxx2} is equivalent to $(1 - \theta_3)x_k + \theta_3 x \in \mathcal R_2$.
On the other hand, we also have $(1 - \theta_3)x_k + \theta_3 x \in \mathcal R_1$ 
by the convexity of $\mathcal R_1$ and $x_k, x \in \mathcal R_1$.
Thus, we have
 \begin{equation}
 (1 - \theta_3)x_k + \theta_3 x \in \mathcal R_3
 \label{eq: f(x + theta_3u) <= m(theta_3u)}
 \end{equation}
by \cref{eq: R1 cap R2 in R3}.
Therefore, we obtain
\begin{alignat}{2}
  f(x_k)
  &=
  f((1 - \theta_3)x_k + \theta_3 x)
  &\quad&\text{(by \cref{eq: t2 such that xxx1})}\\
  &\leq \tilde m((1 - \theta_3)x_k + \theta_3 x)
  &\quad&\text{(by \cref{eq: f(x + theta_3u) <= m(theta_3u)})}\\
  &< (1-\theta_3) \tilde m(x_k) + \theta_3 \tilde m(x)
  &\quad&\text{(by the strictly convexity of $\tilde m$)}\\
  &\leq \tilde m(x_k)
  &\quad&\text{(by $x \in \mathcal R_1$)}\\
  &=
  f(x_k)
  &\quad&\text{(by the assumption on $\tilde m$)},
\end{alignat}
which is a contradiction.
\end{proof}

Now, we prove \cref{lem: majorization LM}.

\begin{proof}[Proof of \cref{lem: majorization LM}]
The model function $m^k_\lambda$ is strictly convex and 
satisfies $m^k_\lambda(x_k) = f(x_k)$ and $\nabla m^k_\lambda(x_k) = \nabla f(x_k)$.
We use \cref{lem: MM technical} with $\tilde m \coloneqq m^k_\lambda$.
Note that 
\cref{eq: m^k(x) <= m^k(x_k),eq: f(x_k+u) <= m(u)} correspond to
$x \in \mathcal R_1$ and $x \in \mathcal R_3$, respectively,
where $\mathcal R_1$ and $\mathcal R_3$ are defined in \cref{eq: def of R1R2R3}.
Thus, by \cref{lem: MM technical}, it suffices to prove \cref{eq: R1 cap R2 in R3}.
We fix $x \in \mathcal R_1 \cap \mathcal R_2$ arbitrarily and will show $x \in \mathcal R_3$.
Let $u \coloneqq x - x_k$.
From the convexity of $\mathcal X$,
$x \in \mathcal R_2$, \cref{eq: J Lip on X}, and \cref{item-lem: Lip J implies xxx},
we have
 \begin{equation}
 \|F(x) - F_k - J_k u\|
 \leq
 \frac{L}{2}\|u\|^2.
 \label{eq: ineq from J lip}
 \end{equation}
From the inequality of arithmetic and geometric means, we have
 \begin{equation}
 \frac{\lambda}{2}\|u\|^2 + \frac{L^2}{2\lambda} \|u\|^2\|F_k + J_k u\|^2
 \geq
 L \|u\|^2\|F_k + J_k u\|.
 \label{eq: ineq from AM-GM}
 \end{equation}
Furthermore, by \cref{eq: def of model} and $x \in \mathcal R_1$, we have
 \begin{equation}
 \|F_k + J_k u\|^2 + \lambda \|u\|^2
 =
 2m^k_\lambda(x)
 \leq
 2m^k_\lambda(x_k)
 =
 2f(x_k)
 =
 \|F_k\|^2.
 \label{eq: ||F + Ju||^2 + lam||u||^2 <= ||F||^2}
 \end{equation}
Using these inequalities, we obtain $x \in \mathcal R_3$ as follows:
 \begin{alignat}{2}
 f(x) - m^k_\lambda(x)
 &=
 \frac{1}{2}\|F(x)\|^2 - m^k_\lambda(x)\\
 &\leq
 \frac{1}{2}\Big(
 \|F_k + J_k u\| + \|F(x) - F_k - J_k u\|
 \Big)^2
 - m^k_\lambda(x)\\
 &\leq
 \frac{1}{2} \Big( \|F_k + J_k u\| + \frac{L}{2}\|u\|^2 \Big)^2
 - m^k_\lambda(x)
 &\quad&\text{(by \cref{eq: ineq from J lip})}\\
 &=
 \frac{L^2}{8}\|u\|^4
 - \frac{\lambda}{2}\|u\|^2
 + \frac{L}{2}\|u\|^2\|F_k + J_k u\|
 &&\text{(by \cref{eq: def of model})}\\
 &\leq
 \frac{L^2}{8}\|u\|^4
 - \frac{\lambda}{4}\|u\|^2
 + \frac{L^2}{4\lambda} \|u\|^2\|F_k + J_k u\|^2
 &&\text{(by \cref{eq: ineq from AM-GM})}\\
 &\leq
 \frac{L^2}{8}\|u\|^4
 - \frac{\lambda}{4}\|u\|^2
 + \frac{L^2}{4\lambda} \|u\|^2 \Big( \|F_k\|^2 - \lambda \|u\|^2 \Big)
 &&\text{(by \cref{eq: ||F + Ju||^2 + lam||u||^2 <= ||F||^2})}\\
 &=
 - \frac{L^2}{8}\|u\|^4
 - \frac{1}{4\lambda} \|u\|^2 \Big( \lambda^2 - L^2 \|F_k\|^2 \Big)\\
 &\leq
 0
 &&\text{(by \cref{eq: lam >= L||F_k||})}.
 \end{alignat}

\end{proof}

The proof of \cref{lem: majorization LM} is a little complicated mainly because $x$ and $y$ 
in \cref{eq: J Lip on X} are restricted on the sublevel set $\mathcal S(x_k)$. 
If we assume the Lipschitz continuity of $J$ on the convex set $\mathcal X$ 
as in \citep{zhao2016global}, \cref{lem: MM technical} is unnecessary and the proof of \cref{lem: majorization LM} can be simplified.

\subsection{Proposition on projected gradient methods}
\label{sec: majorization lemma PG}
\begin{proposition}
 \label{prop: PG decrease}
Fix a point $x_k \in \mathcal C$. For some constant $L_f > 0$, assume that
 \begin{equation}
 \label{eq: Lip-like nabla f on S}
 \|\nabla f(y) - \nabla f(x)\| \leq L_f\|y - x\|,
 \quad
 \forall x,y \in \mathcal C \ \ \text{s.t.} \ \ \mathcal L(x, y) \subseteq \mathcal S(x_k).
 \end{equation}
Then, for $\eta \geq L_f$,
 \begin{align}
 &
 \mathcal L(x_k, \mathcal P_\eta(x_k)) \subseteq \mathcal S(x_k),
 \label{eq: (1-t)x + ty in S PG}\\
 &
 f(\mathcal P_\eta(x_k)) - f(x_k)
 \leq
 - \mathcal D_\eta(x_k).
 \label{eq: sufficient decrease PG}
 \end{align}
\end{proposition}

\begin{proof}
For $\eta \geq L_f$, we define
 \begin{equation}
 \tilde m(x)
 \coloneqq
 f(x_k) + \inner{\nabla f(x_k)}{x - x_k} + \frac{\eta}{2}\|x - x_k\|^2,
 \label{eq: def of model for PG}
 \end{equation}
and use \cref{lem: MM technical} with this function and $\mathcal X = \mathcal C$.
Note that this $\tilde m$ is strictly convex and satisfies that
$\tilde m(x_k) = f(x_k)$ and $\nabla \tilde m(x_k) = \nabla f(x_k)$.
By \cref{eq: Lip-like nabla f on S}, \cref{rem: Lip grad f implies xxx}, and $\eta \geq L_f$,
we have $\mathcal R_2 \subseteq \mathcal R_3$,
where $\mathcal R_2$ and $\mathcal R_3$ are defined in \cref{eq: def of R1R2R3},
and we therefore have \cref{eq: R1 cap R2 in R3}.
Thus, by \cref{lem: MM technical}, we obtain $\mathcal R_1 \subseteq \mathcal R_2 \subseteq \mathcal R_3$, 
which yields
$\mathcal P_\eta(x_k) = \argmin_{x \in \mathcal C} \tilde m(x) \in \mathcal R_1 \subseteq \mathcal R_2 \subseteq \mathcal R_3$.
The first result \cref{eq: (1-t)x + ty in S PG} is
equivalent to $\mathcal P_\eta(x_k) \in \mathcal R_2$,
and the second \cref{eq: sufficient decrease PG} is equivalent to $\mathcal P_\eta(x_k) \in \mathcal R_3$
since $\tilde m (\mathcal P_\eta(x_k)) = f(x_k) - \mathcal D_\eta (x_k)$.
\end{proof}


\subsection{Proof of \texorpdfstring{\cref{lem: small grad map implies approx stationary point}}{}}
\label{sec: proof of lemma GMtoASP}
To prove \cref{lem: small grad map implies approx stationary point}, we first show the following Lipschitz-like property on $\nabla f$.

\begin{lemma}
  \label{lem: asms imply Lip-like nabla f}
  Let \cref{asm: for global convergence} hold and define $L_f$ by \cref{eq: def L_f}.
  Then, for $\eta \geq L_{f}$, we have
  \begin{equation}
  \|\nabla f(\mathcal P_\eta(x)) - \nabla f(x)\|
  \leq
  \eta \|\mathcal P_\eta(x) - x\|,\quad
  \forall x \in \mathcal C \cap \mathcal S(x_0).
  \end{equation}
\end{lemma}

\begin{proof}
  Fix a point $x' \in \mathcal C \cap \mathcal S(x_0)$ arbitrarily. Since $\|F(x')\| \leq \|F(x_0)\| = \|F_0\|$, \cref{asm: for global convergence,item-lem: bounded Lip J implies xxx} with $\mathcal X = \mathcal C \cap \mathcal S(x')$ imply
  \begin{equation}
    \label{eq: Lip-like grad f on S}
    \|\nabla f(y) - \nabla f(x)\|
    \leq
    L_f \|y - x\|,\quad
    \forall x,y \in \mathcal C \ \ \text{s.t.} \ \
    \mathcal L(x, y) \subseteq \mathcal S(x').
  \end{equation}
  By \cref{prop: PG decrease,eq: Lip-like grad f on S}, we have $\mathcal L(x', \mathcal P_\eta(x')) \subseteq \mathcal S(x')$ for $\eta \geq L_f$. Therefore, by using \cref{eq: Lip-like grad f on S} again, we obtain
  \begin{equation}
    \|\nabla f(\mathcal P_\eta(x')) - \nabla f(x')\|
    \leq
    L_f \|\mathcal P_\eta(x') - x'\|
    \leq
    \eta \|\mathcal P_\eta(x') - x'\|,
  \end{equation}
  which is the desired result.
\end{proof}

Now, we prove \cref{lem: small grad map implies approx stationary point}.

\begin{proof}
  [Proof of \cref{lem: small grad map implies approx stationary point}] Since
  \begin{equation}
    \mathcal G_\eta(x)
    = \eta(x - \mathcal P_\eta(x))
    \in \nabla f(x) + \partial \iota_{\mathcal C}(\mathcal P_\eta(x))
    \label{eq: G_eta(x) in F(x) + xxx}
  \end{equation}
  from the definitions of $\mathcal P_\eta$ and $\mathcal G_\eta$ in \cref{eq: def of P_eta,eq: def of G_eta}, we obtain
  \begin{alignat}{2}
    &\mathInd
    \min_{g \in \partial \iota_{\mathcal C}(\mathcal P_\eta(x))}
    \|\nabla f(\mathcal P_\eta(x)) + g\|\\
    &\leq
    \min_{g \in \partial \iota_{\mathcal C}(\mathcal P_\eta(x))}
    \|\nabla f(x) + g\|
    + \|\nabla f(\mathcal P_\eta(x)) - \nabla f(x)\|\\
    &\leq
    \|\mathcal G_\eta(x)\| + \eta \|\mathcal P_\eta(x) - x\|
    &\quad&\by{\cref{eq: G_eta(x) in F(x) + xxx} and \cref{lem: asms imply Lip-like nabla f}}\\
    &=
    2 \|\mathcal G_\eta(x)\|
    &\quad&\text{(by \cref{eq: def of G_eta})}.
  \end{alignat}
  \par\vspace{-0.2\baselineskip}
\end{proof}


\subsection{Relaxing an assumption in \cref{item-thm:successful_iter_global}}
\label{subsec:without_etabar}

In order to compute an $\epsilon$-stationary point based on \cref{item-thm:successful_iter_global}, knowledge of the value of $\bar \eta$ is required. However, this requirement can be circumvented with a slight modification of the algorithm. We show it in this section.

Let $\eta_k$ be the value of $\eta$ when $x_k$ is obtained in \cref{alg: proposed LM-GD}.
As with the proof of \cref{item-thm:successful_iter_global}, we can show that
\begin{align}
  \norm{\mathcal G_{\eta_{k+1}}(x_k)} \leq \epsilon / 2
  \label{eq:gradmap_norm_upper}
\end{align}
for some $k = O(\epsilon^{-2})$.
If \cref{eq:gradmap_norm_upper} and $\eta_{k+1} \geq L_f$ hold, then $\mathcal P_{\eta_{k+1}}(x_k)$ is an $\epsilon$-stationary point by \cref{lem: small grad map implies approx stationary point}, but $\eta_{k+1} \geq L_f$ is not necessarily true.
To address this issue, we modify \cref{alg: proposed LM-GD} a little.

As can be seen from the proof of \cref{lem: small grad map implies approx stationary point}, even if $\eta < L_f$, the point $\mathcal P_\eta(x_k)$ is an $\epsilon$-stationary point of problem~\cref{eq:mainproblem} as long as $\norm{\mathcal G_{\eta}(x_k)} \leq \epsilon / 2$ and the following hold:
\begin{align}
  \|\nabla f(\mathcal P_{\eta}(x_k)) - \nabla f(x_k)\| \leq \eta \|\mathcal P_{\eta}(x_k) - x_k\|.
  \label{eq:liplike_gradf_xk}
\end{align}
Thus, by updating $\eta$ by $\eta \gets \alpha_{\mathrm{in}} \eta$ until \cref{eq:liplike_gradf_xk} is satisfied when $x_k$ is obtained in \cref{alg: proposed LM-GD}, we can guarantee that $\mathcal P_{\eta_{k+1}}(x_k)$ is an $\epsilon$-stationary point for some $k = O(\epsilon^{-2})$.
Since \cref{eq:liplike_gradf_xk} must hold for $\eta \geq L_f$ by \cref{lem: asms imply Lip-like nabla f}, this modification of the algorithm does not sacrifice \cref{lem:inner_global} and other convergence guarantees.
The important point here is that we can check if \cref{eq:liplike_gradf_xk} holds with no prior knowledge of constants (e.g., $\sigma$, $L$, and $L_f$) of the problem.
We have obtained the modified algorithm that does not require the knowledge of the constants.

\section*{Acknowledgments}
We gratefully acknowledge the constructive comments of the anonymous referees.
This work was supported by JSPS KAKENHI Grant Numbers 20K19748 and 19H04069,
and JST ERATO Grand Number JPMJER1903.

\section*{Data availability}
The source code used in our numerical experiments is available on \url{https://github.com/n-marumo/constrained-lm}.

\bibliography{myrefs}

\end{document}